\documentclass[12pt]{article}
\usepackage[cmtip,arrow]{xy}
\usepackage{pb-diagram, pb-xy}
\usepackage{amssymb,epsfig}
\newtheorem{defin}{Definition}
\newtheorem{prop}{Proposition}
\newtheorem{nt}{Remark}
\newtheorem{Th}{Theorem}
\newtheorem{lemma}{Lemma}
\newtheorem{cons}{Corollary}

\newtheorem{example}{Example}

\newfont{\sdbl}{msbm9}
\newfont{\dbl}{msbm10 at 12pt}
\newcommand{\eqdef}{\stackrel{\rm def}{=}}
\newcommand{\proof}{{\bf Proof\ }}
\newcommand{\Ob}{\mathop {\rm Ob}}

\newcommand{\ff}{{\cal F}}

\newcommand{\betta}{\beta}

\newcommand{\Hom}{\mathop {\rm Hom}}
\newcommand{\Aut}{\mathop {\rm Aut}}
\newcommand{\Mor}{\mathop {\rm Mor}}

\newcommand{\Lim}{\mathop {\rm lim}}

\newcommand{\Supp}{\mathop {\rm supp}}

\newcommand{\dz}{{\mbox{\dbl Z}}}

\newcommand{\Z}{\dz}

\newcommand{\veps}{\varepsilon}

\newcommand{\lto}{\longrightarrow}

\newcommand{\df}{{\mbox{\dbl F}}}
\newcommand{\sdf}{{\mbox{\sdbl F}}}
\newcommand{\cc}{{\bf C}}
\newcommand{\F}{{\bf F}}
\newcommand{\I}{{\bf I}}
\newcommand{\D}{{\cal D}}
\newcommand{\E}{{\cal E}}
\newcommand{\oc}{\otimes_{\cc}}
\newcommand{\znakneponyaten}{{\bf 1}}
\newcommand{\cO}{\hat {\cal O}}

\newcommand{\Cal}{\cal}
\newcommand{\pr}{\prod\limits}
\def\Bbb#1{{\mathbb#1}}
\def\R{{\mathbb R}}
\def\Z{{\mathbb Z}}

\def\C{{\mathbb C}}
\def\A{{\mathbb A}}
\newcommand{\ab}{{\rm ab}}

\newcommand{\Gal}{{\rm Gal}}
\newcommand{\Fq}{{\mathbb F}_q}

\begin{document}
\author{D.V. Osipov, A.N. Parshin
\footnote{Both authors are supported by RFBR (grant no.
08-01-00095-a), by a program of President of RF for supporting of
Leading Scientific Schools (grant no. 1987.2008.1), by INTAS (grant
no. 05-1000008-8118); besides the first author is supported by a
grant of Russian Science Support Foundation and by a program of
President of RF for supporting of young russian scientists (grant
no. MK-864.2008.1).} }

\title{Harmonic analysis on local fields and adelic
spaces I}
\date{}

\maketitle \abstract{ Harmonic analysis is developed on objects of
some category $C_2$ of infinite-dimensional filtered vector spaces
over a finite field. This category includes two-dimensional local
fields and adelic spaces of algebraic surfaces defined over a finite
field. The main result is the theory of the Fourier transform on
these objects and obtaining of two-dimensional Poisson formulas.}

\section{Introduction}
Local fields and  adelic groups composed from them are   well known
tools of arithmetic. They were introduced by C. Chevalley in the
1930s and were used to formulate and   solve
many different problems in number theory and algebraic geometry
(see, for example, \cite{A,W}). These constructions are associated
with  fields of algebraic numbers and  fields of algebraic
functions of one variable, i.e.  with schemes of dimension $1$. A need for
such constructions  in higher dimensions was realized by A. N.
Parshin in the  1970s. He managed to develop these constructions
in  the local case for any dimension and  in the global case for
dimension $2$, \cite{P3}. For the schemes of an arbitrary dimension
this approach was developed by A. A. Beilinson, \cite{B}.

If $X$ is a scheme of dimension $n$ and of finite type over $\Z$,
and $X_n\subset X_{n-1}\subset \dots \subset X_1\subset X_0 \subset
X$ is a flag of irreducible subschemes ($\mbox{codim}(X_i)=i$), then
one can define a ring $K_{X_0,\dots,X_{n-1}}$ associated to the
flag. In case when all the subschemes are regularly embedded, the ring is an
$n$-dimensional local field. Then one can form an adelic object
$$\Bbb A_X= {\prod}' K_{X_0,\dots,X_{n-1}}$$
where the product is taken over all the flags, but  with
certain restrictions on  the components of adeles (see \cite{B, H, P3,
PF, Osi}). For a scheme over a finite field $\Fq$ this definition is
the ultimate definition of the adelic space attached to $X$. In
general, one has to extend it to $\Bbb A_X \oplus \Bbb A_{X \otimes
\R}$  by adding archimedian components.

In dimension $1$, the adelic groups $\Bbb A_X$ and $\Bbb A_X^*$ are
locally compact groups and thus we can develope the classical
harmonic analysis on these groups. The starting point for that is
the measure theory on locally compact local fields attached to
points on schemes $X$ of dimension $1$.

This harmonic analysis can be applied to a study of zeta- and
$L$-functions of one-dimensional schemes $X$. J. Tate  and
independently K. Iwasawa have introduced an analytically defined
$L$-function
$$L(s,\chi,f)=\int_{\A^*} \chi(a)|a|^sf(a)
d^* \mu,$$ where $d^* \mu$ is a Haar measure on $\Bbb A^*_X$, the
function $f$ belongs to the Bruhat-Schwartz space of functions on
$\Bbb A_X$ and $\chi$ is an abelian character, coming from the
Galois group $\Gal(K^{\ab}/K)$ by the reciprocity map
$$\A_X^* \rightarrow \Gal(K^{\ab}/K).
$$
 For $L(s,\chi,f)$,  they have proved an
analytical continuation to the whole $s$-plane and the functional
equation
$$
L(s,\chi,f) = L(1-s,\chi^{-1},{\hat f}),
$$
using Fourier transform $f \mapsto {\hat f}$ and Poisson formula for
functions $f$ on the space $\Bbb A_X$ \cite{A,W}. If $L_X(s,\chi)$
is the arithmetically defined $L$-function of the scheme $X$, then
$L_X(s,\chi) = L(s,\chi,f)$ for an appropriate function $f$.

 For a long time, the second author had been pointing  to an essential gap in arithmetics related with
absence of such approach for schemes of  dimen\-sion $> 1$ and
suggesting that the gap should be  filled. We note that there is a
cohomological formalism defined for schemes over  $\df_q$ and
working equally well in all dimensions.

{ \bf Problem 1.} {\em Extend Tate--Iwasawa's analytic method to
higher dimensions} (see, in particular, \cite{P2,P1}).

The higher adeles have been introduced exactly for this purpose.  So
we have the following underlying

 { \bf Problem 2.} {\em
Develop a measure theory and harmonic analysis on $n$-dimensional
local fields.}

Note that the $n$-dimensional local fields are not locally compact
topological spaces for $n>1$ and by Weil's theorem the existence of
the Haar measure (in the usual sense) on a topological group implies
its locally compactness.

In the paper, we show how to construct a harmonic analysis and a
measure theory in the first non-trivial case of algebraic surface
$X$ over a finite field $k$ (or two-dimensional schemes over $\Z$).
In this situation, the corresponding adelic objects will look as
follows. Let $P$ be a closed point of $X$, $C\subset X$ be an
irreducible curve such that $P\in C$.

If  $X$ and $C$ are smooth at $P$ then let $t\in {\cal O}_{X,P}$ be
a local equation of $C$ at $P$ and $u\in {\Cal O}_{X,P}$ be such
that $u|_C\in {\Cal O}_{C,P}$ is a local parameter at $P$. Denote by
$\wp$ the ideal defining the curve $C$ near $P$. Now we can
introduce the two-dimensional local field attached to the pair $P,
C$ by the procedure including completions and localizations: $
K_{P,C} = \mbox{Frac}(\widehat{(\hat{\Cal O}_{X,P})_\wp}) =
\mbox{Frac}(\widehat{\Cal O}_{x,C}) = k(P)((u))((t))$. This
definition can be easily extended to the case of a singular point
$P$ on the curve $C$. Then one can take the adelic product of the
local fields over all the pairs $P, C$. We get the adelic space
$\A_X$ on the surface $X$.

If $C$ is a curve then the space $\A_C$ contains important subspaces
$\A_0 = K = k(C)$ of principal adeles and $ \A_1 = \pr_{x \in C}
\cO_x$ of integral adeles, from which one constructs the adelic
complex: $ \A_0 \oplus \A_1 \rightarrow \A_C $. In dimension $2$,
there is a much more complicated structure of subspaces in $\A_X$.
Among the others, it includes subspaces $\A_{12} = {\prod}'_{x\in C}
\widehat{\Cal O}_{x,C}$ and, more generally, $\A_{12}(D)$, where the
divisor $D$ indicates on some regularity condition around the $D$
(see \cite{PF} and a brief exposition in \cite{P, P1}).

In dimension $1$, the group $\A_C$ is a locally compact group and
the spaces $\A_0$ and $\A_1$ are, correspondingly, discrete and
compact. The construction of analysis starts with a definition of
functional spaces. Let $V$ be a finite dimensional vector space
either over an adelic ring $\A_C$, or over an one-dimensional local
field $K$ with finite residue field $\Fq$, or a subspace of type
$\A_0$ or $\A_1$. We put
$$
\begin{array}{rcl}
{\cal D}(V) &=& \{\mbox{locally constant functions with compact support}\}\\
\tilde{{\cal E}}(V) &=& \{\mbox{uniformly locally constant functions}\}\\
{\cal E}(V) &=& \{\mbox{all locally constant functions}\}\\
{\cal D}'(V) &=& \{\mbox{dual to}~{\cal D} \mbox{, i.e. all distributions}\}\\
\tilde{{\cal E}}'(V) &=& \{\mbox{ "continuously" dual
to}~{\tilde{\cal E}}\}\\
{\cal E}'(V)&=& \{\mbox{"continuosulsy" dual to} ~{\cal E} \mbox{,
i.e. distributions with compact support}\} \mbox{.}
\end{array}
$$

These are the classical spaces introduced by F. Bruhat \cite{Br}.
The harmonic analysis includes definitions of direct and inverse
images in some category $C_1$ of the spaces like $V$, a definition
of the Fourier transform $\F$ as a map from ${\cal D}'(V) \otimes
\mu (V)^*$ to ${\cal D}'(\check V)$ as well for the other types of
spaces. Here $\mu (V)$ is a space of the Haar measures on $V$ and
$\check V$ is a dual object (continuous linear functionals from  the
full dual space). The main result is the following Poisson formula
$$
\F(\delta_{W,\mu_0} \otimes \mu^{-1}) = \delta_{{W^{\perp} \cap
\check{V}},\mu^{-1}/\mu^{-1}_0}
$$
for any closed subgroup $i \, : \, W \rightarrow V$. Here $\mu_0 \in
\mu(W)$, $\mu \in \mu(V)$, $\delta_{W,\mu_0} = i_*(\mu_0)$ (we use
that, by definition of measures,  $\mu(W) \subset \D'(W)$, and the
direct image $i_*$ maps  ${\cal D}'(W)$ to ${\cal D}'(V)$), and
${W^{\perp} \cap \check{V}}$ is the annihilator of $W$ in $\check
V$.

All these constructions are carefully done in  section~\ref{1-dim}
(in section~\ref{0-case}, we collected the corresponding facts for
the case of dimension $0$, section~\ref{not} contains the basic
notations used throughout the paper). The main issues here are the
theory of $C_0$ and $C_1$ spaces (section~\ref{c_1-spaces}),
definition of the spaces of functions and distributions ${\cal D}$,
  ${\cal D}'$, ${\cal E}$, ${\cal E}'$, $\tilde{{\cal E}}$,
 $ \tilde{{\cal E}}'$ (section~\ref{sfd}), invariant Haar
measures (section~\ref{imi}),  Fourier transform (section~\ref{FT}),
 Poisson formula (section~\ref{PF}), Fubini formula and table of all
 possible direct and inverse images for the
 morphisms in the category $C_1$ (section~\ref{DI}),  and the base change rules (section~\ref{bcr}).
 The latter one is very useful when we extend the
 theory to the higher dimensions.
 The exposition is completely independent  of the standard
ones.

This general formalism can be applied to the adelic space on a curve
$C$ and  we immediately get the following facts of analysis on the
self-dual group $\A_C$: $ \F(\delta_{\A_1(D)}) = \mbox{vol}(\A_1(D))
\delta_{\A_1((\omega)-D)}, $ and $ \F( \delta_K) =
\mbox{vol}(\A/K)^{-1}\delta_K $ for the standard subgroups in
$\A_C$(attached to a divisor $D$ and to the
principal adeles, respectively). This easily implies Riemann-Roch and Serre
duality for divisors on the curve $C$ (see\cite{P1}).

For dimension $2$, the space $\A_X$ has a filtration by the
subspaces $\A_{12}(D)$ where $D$ runs through the Cartier divisors
on $X$. The quotients $\A_{12}(D)/\A_{12}(D')$ will be  spaces of
the type considered above in the one-dimensional situation. This
allows us to introduce the following spaces of functions
(distributions):
$$
\begin{array}{rcccl}
{\cal D}_{P_0}(V) &=& \mbox{lim} & \mbox{lim} & {\cal D}(P/Q) \otimes \mu(P_0/Q) \mbox{,} \\
        &&  \stackrel{\longleftarrow}{j^*} & \stackrel{\longleftarrow}{i_*} & \\
{\cal D'}_{P_0}(V) &=& \mbox{lim} & \mbox{lim} & {\cal D}'(P/Q)  \otimes \mu(P_0/Q)^* \mbox{,} \\
        &&  \stackrel{\longrightarrow}{j_*} & \stackrel{\longrightarrow}{i^*} & \\
{\cal E}(V) &=& \mbox{lim} & \mbox{lim} & {\cal E}(P/Q) \mbox{,} \\
&& \stackrel{\longleftarrow}{j^*} &
\stackrel{\longrightarrow}{i^*} & \\
{\cal E'}(V) &=& \mbox{lim} & \mbox{lim} & {\cal E}'(P/Q) \mbox{,} \\
&&  \stackrel{\longrightarrow}{j_*} &
\stackrel{\longleftarrow}{i_*} & \\
\tilde{{\cal E}}(V) &=& \mbox{lim} & \mbox{lim} & \tilde{{\cal E}}(P/Q) \mbox{,} \\
&&  \stackrel{\longrightarrow}{i^*}    &
\stackrel{\longleftarrow}{j^*} & \\
\tilde{{\cal E}}'(V) &=& \mbox{lim} & \mbox{lim} & \tilde{{\cal E}}'(P/Q) \mbox{,} \\
& & \stackrel{\longleftarrow}{i_*}
 &  \stackrel{\longrightarrow}{j_*}\\
&&&& \\
\end{array}
$$
where $P \supset Q \supset R$ are some elements of the filtration in
$\A_X$(more generally, in a reasonable filtered space $V$ with
locally compact quotients), $P_0$ is a fixed subspace from the
filtration and $j: Q/R \rightarrow P/R$, $i: P/R \rightarrow P/Q$
are the canonical maps. Note that both spaces ${\cal D}_{P_0}(V)$,
${\cal D}_{P_0}'(V)$ are ${\cal E}(V)$-modules. In contrast to the one-dimensional
case, there are no simple inclusions between these spaces.

Just as in the case of dimension $1$, one can introduce  some
category $C_2$ (\cite{Osip}) of filtered spaces like $V$  that
includes {\em all} components of the adelic complex such as
$\A_{12}(D)$ and the space $\A_X$ itself.

We develop the theory of $C_2$ spaces in section~\ref{c_2-spaces}.
In this category, it is possible to define (section~\ref{vm}) the
spaces of virtual measures $\mu(F(i) | F(j))$ which are crucial for
the definition of the spaces of functions and distributions
(section~\ref{bs}). In section~\ref{ft} we introduce the Fourier
transform $\F$, which preserves the spaces ${\cal D}$ and ${\cal
D'}$, but interchanges the spaces ${\cal E}$ and $\tilde{{\cal
E'}}$. Next section~\ref{ce} contains a study of a central extension
of automorphism group of an object $V$ in the category $C_2$ and its
action on the functional spaces such as ${\cal D}_{P_0}(V)$ and
${\cal D}'_{P_0}(V)$. The definition of central extension goes up to
the paper of V. G. Kac and D. H. Peterson, \cite{KP}, but the
definition of an action of group  was introduced by M. M. Kapranov
for the case of local fields \cite{K}. We then  study direct and
inverse images (section~\ref{secdi}), base change rules
(section~\ref{bcr2}), the Fourier transform $\F$
(section~\ref{ftdi}), which interchanges direct and inverse images,
and, at last, we can prove a generalization of the Poisson formula
(section~\ref{2PF}).
 It is important that we have two basically
 different types of this formula (theorems~\ref{tpf1} and~\ref{th3}).
The formulation includes  characteristic functions $\delta_W$ of
subspaces $W\subset V$ and resembles the formula for dimension $1$.
It is important that for a class of spaces $V$ (but not for $\A_X$
itself) there exists an invariant measure ($\znakneponyaten_{\mu}$
or $\znakneponyaten$ ), defined up to a constant, as an element of
${\cal D'}(V)$ or ${\cal D}(V)$, an analogue of the classical Haar
measure (see section~\ref{2PF}).

There is also an analytical expression for the intersection number
of two divisors based on an adelic approach to the intersection
theory, \cite{P,PF}. As a corollary, we get an analytical proof of
the (easy part of) Riemann-Roch theorem for divisors on $X$. This
will be exposed in a separate paper.

The main principles of the construction of harmonic analysis were
outlined by the second author (see, in particular \cite{P1,P2,P11}).
The key point here is an observation concerning  structure of the
classical Bruhat functional spaces on local fields or adelic spaces.
Namely, they can be presented as  double (projective and/or
inductive) limit of functional spaces on finite groups
(in particular, on finite-dimensional vector spaces over a finite
field). Thus, the construction of harmonic analysis must  start with
the case of dimension 0 (finite-dimensional vector spaces over a
finite field representing a scheme of dimension 0 or finite abelian
groups) and then be developed by induction to the higher dimensions.
This plan is completely fulfilled in this paper in the cases of
dimensions 1 and 2.

  An important
contribution belongs to M. M. Kapranov \cite{K} who suggested to use
a trick from a construction of the Sato Grassmanian in the theory of
integrable systems. The trick is to use in the above definition of
the spaces ${\cal D}_{P_0}(V)$ and ${\cal D'}_{P_0}(V)$ the spaces
$\mu(P_0/Q)$ of measures instead of $\mu(P/Q)$. Without this trick
with the space of measures we cannot define the functional spaces
for {\em all} vector spaces from the category $C_2$ and,
in particular, for the whole adelic space $\A_X$.

The second author has tried to construct the harmonic analysis in
several ways. The most transparent and the easiest one to use is to
attach several types of functional spaces to the objects of a {\em
category} of vector spaces arising from some subspaces of the whole
adelic space. He worked with a concrete category arising from the
sub-factors of the adelic space. The key point here is a thorough
study of the conditions on  homomorphisms of vector spaces that
guarantee the existence of direct or inverse images.

Approximately at the same time,  the first author independently
introduced a notion of $C_n$ structure in the category of vector
spaces~\cite{Osip}. With this notion at hand, we can develop the
analysis in a very general setting, for any objects of the category
$C_2$ and hopefully for $C_n$. The crucial point is that the $C_n$
structure exists for the adelic spaces of {\em any} $n$-dimensional
Noetherian scheme \cite[Theorem 2.1]{Osip}. The principal advantage
of this approach is that one can make all the constructions {\em
simultaneously} for local and global cases. The category $C_1$
contains as a full subcategory the category of linearly locally
compact vector spaces (introduced and extensively used by S.
Lefshetz \cite{L}) and one can use the classical harmonic analysis
under these circumstances.

Note here that the previous publications \cite{P1,P2,P11} contained
a stupid mistake. Two different types  of locally constant
functions, ${\cal E}(V)$ and $\tilde{{\cal E}}(V)$ has not been
distinguished (the same was true for distributions).

The present work can certainly be extended, first to the case of
spaces arising from arithmetic surfaces. The corresponding category
$C_2$ must include more general abelian groups that appear from
local fields of unequal characteristics and the fields of
archimedean type. In the last case these are the vector spaces over
$\R$ or $\C$. There are no principal obstacles to do that and we
hope to consider this issue in another paper. It seems more difficult but still
possible to extend the theory to the case of higher dimensions, so
the category~$C_n$ \footnote{In the beginning of July 2007, our
paper was put on the net arXiv:0707.1766 [math.AG].  Then, there
appeared a paper by A. Deitmar arXiv:0708.0322 [math.NT], where he
claimed that he constructed harmonic analysis for all dimensions
$n$. Unfortunately, the first versions of his work contained serious
mistakes that were pointed out to the author in our letters. The
latest version of his work contains only a construction of the
spaces like ${\cal E}$, which is in fact  a rather trivial
generalization of our work.}.

In the present work, we develop the harmonic analysis on vector
spaces defined either over a local field  $K$ or over an adelic ring
$\A$. In classical case of dimension 1, the analysis can be
developed on arbitrary varieties (defined either over $K$, or over $\A$).
This was already done by Bruhat in local case \cite{Br}. For
arbitrary varieties defined over a two-dimensional local field $K$,
the analysis was constructed by D. Gaitsgory and D.~A. Kazhdan in
\cite{GK} for the purposes of representation theory of reductive
groups over the field $K$ (see also arXiv:math/0302174 [math.RT],
arXiv:math/0406282 [math.RT], arXiv:math/0409543 [math.RT]). It had
been preceded by \cite{K1}, where harmonic analysis was developed on
homogenous spaces such as $G(K)/G(O'_K)$ (introduced in \cite{P6}).
Note that the construction of harmonic analysis (over $\A$) is
 a non-trivial problem  even for the case $X = G_m$. This is the
subject for further discussion in a separate paper.

Different aspects of the harmonic analysis and  its possible
applications were exposed by the second author in several talks and
  lecture series: M\"unster Universit\"at 1999, 2002, Universit\'e
Paris VI 2002, 2005, Moscow (Steklov institute) 2002, 2006,
Saint-Petersburg (LOMI and Euler Institute) 2003, 2005, Salamanca
University 2004, Berlin (Humboldt Universit\"at) 2005, Oberwolfach
2005.

\section{Notations} \label{not}
For any vector space $V$ over a field $k$ we will denote by $V^*$
its dual space, i.e.,
$$V^* = \Hom\nolimits_k (V, k) \mbox{.} $$

For a vector $k$-subspace $H \subset V$ we denote by $H^{\bot}$ the
following $k$-subspace in the space $V^*$:
$$H^{\bot} \eqdef \{ u \in V^* : u(H) = 0    \} \mbox{.}$$

By $\ff(V)$ we denote the space of all functions on $V$ with values
in the field $\cc$.

For any $a \in V$, for any $f \in \ff(V)$ by $T_a(f)$ we denote the
following function from $\ff(V)$:
$$
T_a(f)(v) \eqdef f(v + a)  \mbox{, \quad} v \in V \mbox{.}
$$

For any $f \in \ff(V)$ by $\check{f}$ we denote the following
function from $\ff(V)$:
$$
\check{f}(v) \eqdef f(-v)  \mbox{, \quad} v \in V \mbox{.}
$$

For a finite set $A$ we denote by $\sharp A$ the number of elements
in the set $A$. By  $\Supp(f)$ we denote the support of a function
$f$.

We fix a finite field $\df_q$ and a nontrivial character $\psi :
\df_q \to {\bf C}^*$.

\section{Zero-dimensional case}
\label{0-case}
Let $V$ be a finite-dimensional vector space
over the field $\df_q$.

Then any element $f \in \ff(V)$ has the following unique
presentation:
$$
f =  \sum_{v \in V} a_v \cdot \delta_v  \mbox{,}
$$
where $a_v \in \cc$, and the functions $\delta_v  \in \ff(V)$ are
defined by the following rule
$$
\delta_v(w) = \left\{
\begin{array}{lrc}
1  \mbox{,} & \mbox{if} & w = v \\
0 \mbox{,} & \mbox{if} & w \ne v \mbox{.}
\end{array}
\right.
$$

On the space $\ff(V)$ we have the following pairing $< \, , \, >_V :
\ff(V) \times \ff(V)  \lto \cc$:
$$
<f,g>_V =  \sum_{v \in V} f(v) \cdot g(v) \mbox{.}
$$

We have the following obvious proposition.
\begin{prop} \label{p}
The pairing $<\, , \,>_V$ is a nondegenerate symmetric pairing  on
the finite-dimensional $\cc$-vector space $\ff(V)$, and therefore
this pairing determines a canonical isomorphism between the vector
space $\ff(V)$ and the dual vector space $\ff(V)^*$.
\end{prop}

Let $\pi : V \to W$ be a linear map between finite-dimensional
vector spaces over the field $\df_q$. We define the direct image
$\pi_* : \ff(V) \to \ff(W)$ as the following map:
$$
\pi_* (f)(w) = \left\{
\begin{array}{lrc}
0  & \mbox{, if} & w \notin \pi(V) \\
 \mathop{\mathop{\sum}\limits_{v \in V}}\limits_{\pi(v) = w} f(v)  & \mbox{, if} & w \in \pi(V) \mbox{,}
\end{array}
\right.
$$
 where $ f \in \ff(V)$, $w \in W $.

We define the inverse image  $\pi^* : \ff(W) \to \ff(V)$ as the
following map:
$$
\pi^* (g) (v) = f( \pi (v)) \quad \mbox{, where} \quad g \in \ff(W),
\quad v \in V \mbox{.}
$$

There is the following proposition, which follows easily from the
definitions.
\begin{prop} \label{pp}
The maps $\pi_*$ and $\pi^*$ are conjugate maps with respect  to the
pairings $<\,, \, >_V$ and $<\, , \,>_W$ on the spaces $\ff(V)$ and
$\ff(W)$, i.e.,
$$
<\pi^*(g), f>_V = <g, \pi_*(f) >_W
$$
for any $f \in \ff(V)$, $g \in \ff(W)$.
\end{prop}

There is also the following proposition, which also easily follows
from the definitions.
\begin{prop} \label{bch}
The following statements are satisfied.
\begin{enumerate}
\item Let $\pi_1 : V_1 \to V_2$ and $\pi_2 :
V_2 \to V_3$ are linear maps between the
vector spaces over the field $\df_q$, then
$$
(\pi_2)_* (\pi_1)_* (f) =  (\pi_2 \pi_1)_* (f) \quad \mbox{for any}
\quad f \in \ff(V_1) \mbox{,}
$$
$$
\mbox{and} \quad  \pi_1^* \pi_2^* (g) = ( \pi_2 \pi_1)^* (g) \quad
\mbox{for any} \quad g \in \ff(V_3) \mbox{.}
$$
\item Let $V$, $W$, $S$ are
finite-dimensional vector spaces over the
field $\df_q$. Let $\pi : V \to S$, $\alpha
: W \to S$ are linear maps between these
spaces. We consider the following cartesian
diagram of vector spaces:
$$
\begin{diagram}
\node{V \mathop{\times}_S W} \arrow{e,t}{\pi_W}
\arrow{s,l}{\alpha_V} \node{W} \arrow{s,r}{\alpha} \\
\node{V} \arrow{e,b}{\pi} \node{S}
\end{diagram} \mbox{.}
$$
Then we have the following base change rule
$$
\pi^* \alpha_* (h) = (\alpha_V)_* \pi_W^* (h) \quad \mbox{for any}
\quad h \in \ff(W) \mbox{.}
$$

\end{enumerate}
\end{prop}

For a finite-dimensional vector space $V$ over the field $\df_q$ we
define the Fourier transform $\F : \ff(V) \lto \ff(V^*)$ by the
following rule:
$$
\F (f) (u) = \sum_{v \in V } f(v) \cdot \overline{\psi(u(v))}
\mbox{,}
$$
where $u \in V^*$, $f \in \ff(V)$. We consider the following
example.
\begin{example} \label{ex1} \em
Let $H \subset V$ be a vector subspace. We consider the following
function $\delta_H =\sum\limits_{v \in H} \delta_v$ from $\ff(V)$.
Then it is easy to see by direct calculation that
$$
\F(\delta_H) = (\sharp H) \cdot \delta_{H^{\bot}} \mbox{.}
$$
 \em
\end{example}
\begin{nt}
{\em Example~\ref{ex1} is {\em the Poisson formula} for
$0$-dimensional case.}
\end{nt}

We have the following  proposition, which collects the properties of
the Fourier transform $\F$.

\begin{prop} \label{f0}
The following statements are satisfied.
\begin{enumerate}
\item \label{stf1} $\F$ is an isomorphism of
$\cc$-vector spaces. \item \label{stf2} $\F
\circ \F (f)  = \sharp V \cdot \check{f}  $
for any $f \in \ff(V)$. \item \label{stf3}
$<\F(f), g>_{V^*} = <f, \F(g)>_V $ for any
$f \in \ff(V)$, $g \in \ff(V^*)$. \item
\label{stf4} Let $\pi : V \to W$ be a linear
map between  the finite-dimensional vector
spaces over the field $\df_q$. Let $\pi' :
W^* \to V^*$ be the dual map. Then the
following diagrams are commutative:
\begin{equation} \label{d1}
\begin{diagram}
\node{\ff(V)} \arrow{e,t}{\pi_*}
\arrow{s,l}{\F} \node{\ff(W)} \arrow{s,r}{\F} \\
\node{\ff(V^*)} \arrow{e,b}{(\pi')^*} \node{\ff(W^*)}
\end{diagram}
\end{equation}
\qquad \qquad
\begin{equation}
\begin{diagram}  \label{d2}
\node{\ff(W)} \arrow{e,t}{\pi^*}
\arrow{s,l}{\frac{1}{\sharp W} \cdot \F} \node{\ff(V)} \arrow{s,r}{\F \cdot \frac{1}{\sharp V} } \\
\node{\ff(W^*)} \arrow{e,b}{(\pi')_*} \node{\ff(V^*)}
\end{diagram}
\end{equation}
\end{enumerate}
\end{prop}
\proof. Statements~\ref{stf2}, \ref{stf3} and diagram~(\ref{d1})
follow by direct calculations  on functions $\delta_v$. Now
statement~\ref{stf1} follows from statement~\ref{stf2}.
Diagram~(\ref{d2}) follows from diagram~(\ref{d1}) and
statement~\ref{stf2}.
\vspace{0.3cm}

\section{One-dimensional case} \label{1-dim}

\subsection{$C_0$ and $C_1$-spaces.} \label{c_1-spaces}
Here we recall the definitions from~\cite{Osip} about $C_0$ and
$C_1$-spaces.

By definition, the category $C_0$ is the category of
finite-dimensional vector spaces over a field $k$ with morphisms
coming from $k$-linear maps between vector spaces.

By definition, an admissible triple of finite-dimensional vector
spaces is an exact triple of  this vector spaces.

We say that $(I, F, V)$ is a filtered $k$-vector space, if
\begin{itemize}
\item[(i)] $V$ is a vector space over the field $k$,
\item[(ii)] $I$ is a
partially ordered set, such that for any $i,j \in  I$ there are $k,l
\in I$ with $k \le i \le l$ and $k \le j \le l$,
\item[(iii)] $F$ is a
function from $I$ to the set of $k$-vector subspaces of $V$ such
that if  $i \le j $ are any from $I$, then $F(i) \subset F(j)$,
\item[(iv)] $ \bigcap\limits_{i \in I}  F(i) = 0$  and $
\bigcup\limits_{i \in I} F(i) = V $.
\end{itemize}

We say that a filtered vector space $(I_1, F_1, V)$ {\em dominates}
another filtered vector space $(I_2, F_2, V)$ when there is an
order-preserving   function $\phi: I_2 \to I_1$ such that
\begin{itemize}
\item[(i)] for any $i \in I_2$ we have $F_1(\phi (i)) = F_2 (i)$,
\item[(ii)] for any $j \in I_1$ there are $i_1, i_2 \in I_2$ such that $ \phi(i_1) \le j \le \phi(i_2)$.
\end{itemize}

By definition, objects of the category $C_1$, i.e. $\Ob(C_1)$, are
filtered $k$-vector spaces $(I, F, V)$ such that  for any $i \le j
\in I$  the $k$-vector space $F(j) / F(i)$ is finite-dimensional
over the field $k$.

Let $E_1 = (I_1, F_1, V_1)$ and $E_2 = (I_2, F_2, V_2)$ be from
$\Ob(C_1)$. Then, by definition, a set $\Mor_{C_1}(E_1, E_2)$
consists of elements $A \in \Hom_k (V_1, V_2)$ under the following
condition: for any $j \in I_2$ there is an $i \in I_1$ such that $A
(F_1(i)) \subset F_2(j)$.

Let $E_1 = (I_1, F_1, V_1)$, $E_2 = (I_2,
F_2, V_2)$ and $E_3=(I_3, F_3, V_3)$ be from
$\Ob(C_1)$. Then we say that
\begin{equation} \label{adt}
0 \arrow{e} E_1  \arrow{e,t}{\alpha}  E_2 \arrow{e,t}{\beta}   E_3
\arrow{e} 0
\end{equation}
is {\em an admissible triple } from $C_1$ when
 the following conditions are satisfied:
\begin{itemize}
\item[(i)]
$$
0 \arrow{e} V_1  \arrow{e,t}{\alpha}  V_2 \arrow{e,t}{\beta}   V_3
\arrow{e} 0
$$
is an exact triple of $k$-vector spaces,
\item[(ii)] \label{itaa}
the filtration $(I_1, F_1, V_1)$ dominates the filtration $(I_2,
F'_1, V_1)$, where $F'_1 (i) = F_2(i) \cap V_1$  for any $i \in
I_2$,
\item[(iii)] \label{itbb}
the filtration $(I_3, F_3, V_3)$ dominates the filtration $(I_2,
F'_3, V_3)$, where $F'_3(i) =  F_2(i) / F_2(i) \cap V_1$ for any $i
\in I_2$,
\item[(iv)] for any $i \le j \in I_2$
$$
0 \lto \frac{F'_1(j)}{F'_1(i)} \lto \frac{F_2(j)}{F_2(i)}  \lto
\frac{F'_3(j)}{F'_3(i)} \lto 0
$$
is an exact triple of finite-dimensional $k$-vector spaces.
\end{itemize}

We say that $\alpha$ from an admissible triple~(\ref{adt}) is {\em
an admissible monomorphism}, and $\beta$ from an admissible
triple~(\ref{adt}) is { \em an admissible epimorphism}.

\begin{nt} \label{dop} \em
Let  $E_2= (I_2, F_2, V_2)$ be a $C_1$-space over the field $k$. Let
\begin{equation} \label{z}
0 \arrow{e} V_1  \arrow{e,t}{\alpha}  V_2 \arrow{e,t}{\beta}   V_3
\arrow{e} 0
\end{equation}
be an exact triple of $k$-vector spaces. We define on the spaces
$V_1$ and $V_3$ the structures of filtered spaces  $(I_i, F_i,
V_i)$, $i =1, 3$, where $I_1=I_3=I_2$, and $F_1(i) \eqdef F_2(i)
\cap V_1$, $F_3(i) \eqdef \beta(F_2(i))$ for $i \in I_2$.

Then the filtered space  $E_1= (I_1, F_1, V_1)$ will be always a
$C_1$-space. But filtered space  $E_3= (I_3, F_3, V_3)$ will be a
$C_1$-space, if and only if the following condition holds:
$$
\label{zamkn} \mbox{for any} \quad x \in V_2 \setminus V_1 \quad
\mbox{there exists} \quad i \in I_2 \mbox{ such that} \quad (x +
F_2(i)) \cap V_1 = \emptyset \mbox{.}
\leqno(*)
$$

Moreover, under  condition~$(*)$, exact triple~(\ref{z}) gives the
admissible triple
$$0 \arrow{e} E_1
\arrow{e,t}{\alpha}  E_2 \arrow{e,t}{\beta}   E_3 \arrow{e} 0
\mbox{,}$$ i.e. $\alpha$ will be an admissible monomorphism, $\beta$
will be an admissible epimorphism.

For example, we consider $k((t))$ as a  $C_1$-space with filtration,
given by all powers of maximal ideal of the ring $k[[t]]$. Then
$k$-subspace $k[t]$ with induced filtration will be a $C_1$-space,
but the  embedding $k[t] \subset k((t))$ will not be an admissible
monomorphism. And embedding $k[t^{-1}] \subset k((t))$ will be an
admissible monomorphism (with the induced filtration).
\end{nt}

Let $E = (I, F, V)$ be a $C_1$-space. We
define a $C_1$-space $\check{E} = (I^0, F^0,
\check{V})$ by the following way. The space
$\check{V} \subset V^*$ is defined as
\begin{equation}
\check {V} \eqdef \bigcup_{i \in I} F(i)^{\bot} \mbox{.}
\end{equation}
The set $I^0$ is a partially ordered set, which has the same set as
$I$, but with the inverse order then $I$. For $j \in I^0$ the
subspace $F^0(j) \eqdef F(j)^{\bot}  \subset \check{V}$.

We note also that
$$F(j)^{\bot} = (V/F(j))^* =\mathop{\Lim_{\gets}}_{i \ge j} (F(i)/F(j))^* \mbox{.}$$

If $E_1, E_2 \in \Ob(C_1)$ and $\theta \in \Mor_{C_1}(E_1,E_2)$.
Then there is canonically $\check{\theta} \in \Mor_{C_1}
(\check{E_2}, \check{E_1})$. If
$$
0 \arrow{e} E_1  \arrow{e,t}{\alpha}  E_2 \arrow{e,t}{\beta}   E_3
\arrow{e} 0
$$
is an admissible triple of $C_1$-spaces, where $E_i = (I_i, F_i,
V_i)$, $ 1 \le i \le 3$. Then there is canonically the following
admissible triple of $C_1$-spaces:
$$
0 \arrow{e} \check{E_3} \arrow{e,t}{\check{\beta}} \check{E_2}
\arrow{e,t}{\check{\alpha}}  \check{E_1} \arrow{e} 0 \mbox{.}
$$

\begin{lemma} \label{lbc}
Let
$$
0 \arrow{e} E_1  \arrow{e,t}{\alpha}  E_2 \arrow{e,t}{\beta}   E_3
\arrow{e} 0
$$
be an admissible triple of $C_1$-spaces.
\begin{enumerate}
\item Let $D \in \Ob(C_1)$, and $\gamma \in
\Mor_{C_1} (D, E_3)$. Then there is the
following admissible triple of $C_1$-spaces
\begin{equation} \label{adm1}
0 \arrow{e} E_1  \arrow{e,t}{\gamma_{\alpha}} E_2
\mathop{\times}_{E_3} D \arrow{e,t}{\gamma_{\beta}} D \arrow{e} 0
\end{equation}
and $\beta_{\gamma} \in \Mor_{C_1}(E_2 \mathop{\times}\limits_{E_3}
D, E_2)$ such that the following diagram is commutative:
\begin{equation} \label{adm2}
\begin{diagram}
\node{0} \arrow{e} \node{E_1}  \arrow{e,t}{\gamma_{\alpha}}
\arrow{s,=} \node{E_2 \mathop{\times}_{E_3} D}
\arrow{e,t}{\gamma_{\beta}}
\arrow{s,l}{\beta_{\gamma}} \node{D} \arrow{e} \arrow{s,r}{\gamma} \node{0} \\
\node{0} \arrow{e} \node{E_1}  \arrow{e,t}{\alpha}  \node{E_2}
\arrow{e,t}{\beta} \node{E_3} \arrow{e} \node{0}
\end{diagram}
\end{equation}
\item Let $A \in \Ob(C_1)$, and $\theta \in
\Mor_{C_1}(E_1, A)$.
 Then there is the following admissible triple of
$C_1$-spaces
\begin{equation}  \label{adm'1}
0 \arrow{e} A  \arrow{e,t}{\theta_{\alpha}} A  \mathop{\amalg}_{E_1}
E_2 \arrow{e,t}{\theta_{\beta}} E_3 \arrow{e} 0
\end{equation}
and $\alpha_{\theta} \in \Mor_{C_1}(E_2, A
\mathop{\amalg}\limits_{E_1} E_2)$ such that the following diagram
is commutative:
\begin{equation}  \label{adm'2}
\begin{diagram}
\node{0} \arrow{e} \node{E_1}  \arrow{e,t}{\alpha}
\arrow{s,l}{\theta} \node{E_2} \arrow{e,t}{\beta}
\arrow{s,r}{\alpha_{\theta}} \node{E_3} \arrow{e} \arrow{s,=}  \node{0} \\
\node{0} \arrow{e} \node{A}  \arrow{e,t}{\theta_{\alpha}} \node{A
\mathop{\amalg}\limits_{E_1} E_2} \arrow{e,t}{\theta_{\beta}}
\node{E_3} \arrow{e} \node{0}
\end{diagram}
\end{equation}
\end{enumerate}
\end{lemma}
\proof. We prove the first statement of the
lemma. Let $E_i = (I_i, F_i, V_i)$, $ 1 \le
i \le 3$, and
 $D = (K, H, T)$.

 We construct $E_2 \mathop{\times}\limits_{E_3}  D = (J, G, W) \in \Ob(C_1)$ in the
 following way. We define a $k$-vector space
 $$
 W = V_2 \mathop{\times}_{V_3} T \eqdef \{ (e,d) \in V_2 \times T \quad \mbox{such that} \quad  \beta(e) = \gamma(d)
 \} \mbox{.}
 $$
We define a partially ordered set
$$
J = \{ (i,j) \in I_2 \times K  \quad \mbox{such that} \quad
\gamma(H(j)) \subset \beta(F_2(i)) \}  \mbox{,}
$$
where $(i_1, j_1) \le (i_2, j_2)$ iff $i_1 \le i_2$ and $j_1 \le
j_2$. We define a function $G$ from $J$ to the set of $k$-vector
subspaces of $W$ as
$$
G((i,j)) = (F_2(i) \times H(j) ) \cap W  \mbox{,}
$$
where the intersection is inside $V_2 \times T$. Then $E_2
\mathop{\times}\limits_{E_3} D$ is well defined as a $C_1$-space,
and the maps $\beta_{\gamma}$ and $\gamma_{\beta}$ are projections.

Let $(i_1, j_1), (i_2, j_2) \in J$, and
$(i_1, j_1) \le (i_2, j_2)$. Then we have
the following commutative diagram of
finite-dimensional $k$-vector spaces:
$$
\begin{diagram}
\dgARROWLENGTH=1.3em
 \node{0} \arrow{e} \node{\frac{F_2(i_2) \cap V_1} {F_2(i_1) \cap
V_1} } \arrow{e}  \arrow{s,=} \node{\frac{G((i_2, j_2))}
{G((i_1,j_1))}} \arrow{e} \arrow{s} \node{\frac{H(j_2)}{H(j_1)}}
\arrow{e} \arrow{s} \node{0}
\\ \node{0} \arrow{e} \node{\frac{F_2(i_2) \cap V_1}{F_2(i_1) \cap
V_1}} \arrow{e}  \node{\frac{F_2(i_2)}{F_2(i_1)}} \arrow{e}
\node{\frac{F_2(i_2)/ (F_2(i_2) \cap V_1)}{F_2(i_1)/ (F_2(i_1)\cap
V_1)}} \arrow{e} \node{0}
\end{diagram}
$$
where
$$
\frac{G((i_2, j_2))} {G((i_1,j_1))} = \frac{F_2(i_2)}{F_2(i_1)}
 \mathop{\times}_{\frac{F_2(i_2)/ (F_2(i_2) \cap
V_1)}{F_2(i_1)/ (F_2(i_1)\cap V_1)}} \frac{H(j_2)}{H(j_1)}  \mbox{.}
$$

From this commutative diagram we obtain that triple~{(\ref{adm1})}
is an admissible triple of $C_1$-spaces, and diagram~{(\ref{adm2})}
is a commutative diagram. The first statement of the lemma is
proved.

We prove the second statement of the lemma,
which is the dual statement to the first
statement of the lemma. Let $A= (K', H',
T')$.

We construct $A  \mathop{\amalg}\limits_{E_1} E_2 = (J', G', W') \in
\Ob(C_1)$ in the following way. We remark that $T' \amalg V_2 = T'
\times V_2$.
 We define a
$k$-vector space
$$
W' = T' \mathop{\amalg}_{V_1}   V_2 \eqdef (T' \times V_2) / E \quad
\mbox{, where} \quad E = (\theta \times \alpha)  V_1  \mbox{.}
$$
We define a partially ordered set
$$
J' = \{(i,j) \in K' \times I_2   \; : \; \theta(F_2(j) \cap V_1)
\subset H'(i) \} \mbox{,}
$$
where $(i_1, j_1) \le (i_2, j_2)$ iff $i_1 \le i_2$ and $j_1 \le
j_2$. We define a function $G'$ from $J'$ to the set of $k$-vector
subspaces  of $W'$ as
$$
G'((i,j)) = ((H'(i) \times F_2(j)) + E / E )      \quad \subset
\quad W'  \mbox{.}
$$
Then $A  \mathop{\amalg}\limits_{E_1} E_2$ is well defined as a
$C_1$-space.

Let $(i_1, j_1), (i_2, j_2)  \in J'$, and
$(i_1, j_1) \le (i_2, j_2)$. Then we have
the following commutative diagram of
finite-dimensional k-vector spaces:
$$
\begin{diagram}
\dgARROWLENGTH=0.9em
 \node{0}  \arrow{e} \node{\frac{F_2(j_2) \cap V_1} {F_2(j_1) \cap
V_1} } \arrow{e}  \arrow{s} \node{\frac{F_2(j_2)}{F_2(j_1)}}
\arrow{e} \arrow{s} \node{\frac{F_2(j_2)/ (F_2(j_2) \cap
V_1)}{F_2(j_1)/ (F_2(j_1)\cap V_1)}} \arrow{e} \arrow{s,=} \node{0}
\\ \node{0}  \arrow{e} \node{\frac{H'(i_2)}{H'(i_1)}} \arrow{e}
\node{\frac{G'((i_2,j_2))}{G'((i_1,j_1))}} \arrow{e}
\node{\frac{F_2(j_2)/ (F_2(j_2) \cap V_1)}{F_2(j_1)/ (F_2(j_1)\cap
V_1)}}
 \arrow{e} \node{0}
\end{diagram}
$$
where
$$
\frac{G'((i_2, j_2))} {G'((i_1,j_1))} = \frac{H'(i_2)}{H'(i_1)}
 \mathop{\amalg}_{ \frac{F_2(j_2) \cap V_1} {F_2(j_1) \cap
V_1} } \frac{F_2(j_2)}{F_2(j_1)}  \mbox{.}
$$

From this commutative diagram we obtain that triple~{(\ref{adm'1})}
is an admissible triple of $C_1$-spaces, and diagram~{(\ref{adm'2})}
is a commutative diagram. The second statement of the lemma is
proved.
\vspace{0.3cm}

\begin{defin}
Let $(I, F, V), (J, G, V) \in \Ob(C_1)$. We say that there is an
{\it equivalence} $(I, F, V) \sim (J, G, V)$ iff there is a
collection of $C_1$-spaces $(I_l, F_l, V)$, $1 \le l \le n$
 such that
\begin{itemize}
\item[(i)] $I_1= I$, $F_1 = F$, $I_n = J$, $F_n =G$,
\item[(ii)] for any $1 \le l \le n-1$ either
$(I_l, F_l, V)$ dominates $(I_{l+1}, F_{l+1}, V)$, or $(I_{l+1},
F_{l+1}, V)$ dominates $(I_l, F_l, V)$.
\end{itemize}
\end{defin}

Let $\theta : E \lto D$ be an admissible
epimorphism of $C_1$-spaces. Let $\beta : B
\lto D$ be an admissible monomorphism of
$C_1$-spaces. Then from construction (see
the proof of lemma~\ref{lbc}) we have that
$\beta_{\theta} : E \mathop{\times}\limits_D
B \lto  B$ is an admissible epimorphism of
$C_1$-spaces, and $\theta_{\beta} : E
\mathop{\times}\limits_D B \lto E$ is an
admissible monomorphism.

Moreover, Let $N = E
\mathop{\times}\limits_D B$. From
construction we have that $D \sim B
\mathop{\amalg}\limits_ N E$ .We have the
admissible epimorphism
$\check{\theta_{\beta}} : \check{E} \lto
\check{N}$, and the admissible monomorphism
$\check{\beta_{\theta}} : \check{B} \lto
\check{N}$. From construction (see the proof
of lemma~\ref{lbc}) we have that $\check{D}
\sim \check{E}
\mathop{\times}\limits_{\check{N}}
\check{B}$.

\begin{defin}
We say that a $C_1$-space $(I_1, F_1, V_1)$ is a {\it compact}
$C_1$-space iff there is an $i \in I_1$ such that $ F_1(i) = V_1$.
We say that a $C_1$-space $(I_2, F_2, V_2)$ is a {\it discrete}
$C_1$-space iff there is an $j \in I_2$ such that $F_2(j) = \{0\}$.
\end{defin}

It follows from definition that if $E_1$ is a compact $C_1$-space,
then $\check{E_1}$ is a discrete $C_1$-space. If $E_2$ is a discrete
$C_1$-space, then $\check{E_2}$ is a compact $C_1$-space.

\begin{defin}
We say that a $C_1$-space $(I,F,V)$ is a complete $C_1$-space iff
$$
V = \mathop{\Lim_{\to}}_{i \in I}  \mathop{\Lim_{\gets}}_{j \le i}
F(i)/F(j)  \mbox{.}
$$
\end{defin}

It follows from definitions that  $E$ is a complete $C_1$-space iff
$\check{\check{E}} = E$.

\subsection{Spaces of functions and distributions.} \label{sfd}
We recall the basic definitions from~\cite{P1} and \cite{P11}.

 Let $E=(I,F,V)$ be a $C_1$-space
over the  field $\df_q$. For any $i,j,k \in I$ such that $i \le j
\le k$ we have the embedding map
$$
\alpha_{ijk} \quad : \quad F(j) / F(i)  \arrow{e,J} F(k) / F(i)
$$
and the surjection map
$$
\pi_{ijk}  \quad : \quad F(k) / F(i)  \arrow{e,A} F(k) /F(j)
$$
of finite-dimensional vector spaces over the field $\df_q$.

Then according to section~\ref{0-case} we have $4$  maps:
$$ \alpha_{ijk}^* \quad : \quad \ff(F(k)/  F(i)) \lto \ff(F(j)/F(i)) \mbox{,} $$
$$ (\alpha_{ijk})_* \quad : \quad  \ff(F(j)/ F(i))  \lto \ff(F(k)/
F(i)) \mbox{,} $$
$$ \pi_{ijk}^* \quad : \quad
\ff(F(k)/F(j))  \lto \ff(F(k) / F(i))  \mbox{,} $$
$$ (\pi_{ijk})_* \quad : \quad \ff(F(k) /F(i)) \lto \ff(F(k) /F(j)) \mbox{.}
$$

According to the maps described, we consider now all possible
combinations of projective and inductive limits of spaces $\ff(F(j)
/ F(i))$ with respect to the indices $i ,j \in I$, $i \le j$. We
define
\begin{equation} \label{d}
\D(E) \eqdef \mathop{\Lim_{\to}}_j \mathop{\Lim_{\to}}_i \ff(F(j) /
F(i)) = \mathop{\Lim_{\to}}_i \mathop{\Lim_{\to}}_j \ff(F(j) / F(i))
\end{equation}
with respect to the maps $\pi_{ijk}^*$ and $(\alpha_{ijk})_*$;
\begin{equation}  \label{d'}
\D'(E) \eqdef \mathop{\Lim_{\gets}}_j \mathop{\Lim_{\gets}}_i
\ff(F(j) / F(i)) = \mathop{\Lim_{\gets}}_i \mathop{\Lim_{\gets}}_j
\ff(F(j) / F(i))
\end{equation}
with respect to the maps $(\pi_{ijk})_*$ and $\alpha_{ijk}^*$;
\begin{equation} \label{e}
\E(E) \eqdef \mathop{\Lim_{\gets}}_j \mathop{\Lim_{\to}}_i \ff(F(j)
/ F(i))
\end{equation}
with respect to the maps  $\alpha_{ijk}^*$ and $\pi_{ijk}^*$;
\begin{equation} \label{te}
\tilde{\E}(E) \eqdef \mathop{\Lim_{\to}}_i \mathop{\Lim_{\gets}}_j
\ff(F(j) / F(i))
\end{equation}
with respect to the maps  $\pi_{ijk}^*$  and  $\alpha_{ijk}^*$;
\begin{equation} \label{e'}
\E' (E) \eqdef \mathop{\Lim_{\to}}_j \mathop{\Lim_{\gets}}_i
\ff(F(j) / F(i))
\end{equation}
with respect to the maps $(\alpha_{ijk})_*$ and $(\pi_{ijk})_*$;
\begin{equation} \label{te'}
\tilde{\E}' (E) \eqdef \mathop{\Lim_{\gets}}_i \mathop{\Lim_{\to}}_j
\ff(F(j) / F(i))
\end{equation}
 with respect to the maps $(\pi_{ijk})_*$ and
$(\alpha_{ijk})_*$.

Due to proposition~\ref{bch} we have that $\cc$-vector spaces
$\D(E)$, $\D'(E)$, $\E(E)$, $\tilde{\E(E)}$, $\E'(E)$,
$\tilde{\E'}(E)$ are well-defined.

Since  $\alpha_{ijk}^* \circ (\alpha_{ijk})_* $ is the identity map
for any $i \le j \le k \in I$, we have the following canonical
embeddings
\begin{equation} \label{em1}
\D(E) \arrow{e,J} \tilde{\E} (E)
\end{equation}
\begin{equation} \label{em2}
\tilde{\E'} (E) \arrow{e,J} \D'(E) \mbox{.}
\end{equation}

We have also the following canonical embeddings
\begin{equation} \label{em3}
\tilde{\E} (E) \arrow{e,J} \E(E)
\end{equation}
\begin{equation} \label{em4}
\E'(E) \arrow{e,J} \tilde{\E'} (E)  \mbox{.}
\end{equation}

Now we clarify the functional sense of spaces $\D(E)$, $\D'(E)$,
$\E(E)$, $\tilde{\E}(E)$, $\E'(E)$.

\begin{prop} \label{ch}
The following statements are satisfied.
\begin{enumerate}
\item The space $\D(E) \subset \ff(V)$ such that:
$$
f \in \D(E)   \quad   \mbox{iff there are} \quad i \le j \in I \quad
\mbox{such that}
$$
$$
f(v) = 0 \quad  \mbox{for any} \quad v \notin F(j) \quad \mbox{and}
\quad  T_w(f) = f  \quad  \mbox{for any} \quad w \in F(i) \mbox{,}
$$
i.e., the space $\D(E)$ consists of all locally constant functions
with compact support on the space $V$.
\item We suppose that
$E=(I,F,V)$ is a complete $C_1$-space. Then the space $\E(E) \subset
\ff(V)$ such that:
$$
f \in \E(E)   \quad   \mbox{iff for any} \quad v \in V \quad
\mbox{there is} \quad i \in I \quad \mbox{(wich depends on $v$)}
$$
$$
 \mbox{such that} \quad f (v + w) = f (v) \quad \mbox{for any} \quad w \in F(i)  \mbox{,}
$$
i.e., the space $\E(E)$ consists of all locally constant functions
on the space $V$.
\item The space $\tilde{\E}(E) \subset
\ff(V)$ such that:
$$
f \in \tilde{\E}(E)   \quad   \mbox{iff there is} \quad i \in I
\quad \mbox{such that} \quad T_w(f) = f \quad \mbox{for any} \quad w
\in F(i) \mbox{.}
$$
\item
\label{s4}
 The space $\D'(E) = \D(E)^*$, i.e., the space $\D'(E)$
consists of distributions on the space $V$.
\end{enumerate}
\end{prop}
\proof. The first, the second, and the third
statement of this proposition are just
reformulation of formulas~(\ref{d}),
(\ref{e}), (\ref{te}). In the second
statement we demand that $V$ is a complete
$C_1$-space, because then for any locally
constant function $f \in \ff(V)$ for any $j
\in I$ there is an $i \in I $ such that $f(v
+ w) = f(v)$ for any $v \in F(j)$, $w \in
F(i)$.

To prove the fourth statement we remark that
$$
\D(E)^* = \mathop{\Lim_{\gets}}_j \mathop{\Lim_{\gets}}_i \ff(F(j) /
F(i))^* \mbox{,}
$$
because the dual space to the space given by inductive limits of
vector spaces is given by projective limits of dual vector spaces.
Now we apply propositions~\ref{p} and \ref{pp} to obtain that
$\D'(E) = \D(E)^*$.
\vspace{0.3cm}

We can not speak about the value of distribution $G \in \D'(E)$ in
the point $v \in V$. But we can say, whether $G$ is equal to zero or
not in the point $v \in V$. We have the following definition.
\begin{defin} \label{ds}
Let $G \in \D'(E)$.
\begin{enumerate}
\item
 For $v \in V$ we say $G(v) = 0$ iff there is $i \in I$ such that
$G (f) = 0$ for any $f \in \D(V)$ with $\Supp f \subset \{ v + F(i)
\}$.
\item
$$
\Supp G \eqdef \{ v \in V \; : \; G(v) \ne 0 \} \mbox{.}
$$
\end{enumerate}
\end{defin}

Due to definition~\ref{ds} we have the following proposition.
\begin{prop} \label{propd}
Let $E=(I,F,V)$ be a complete $C_1$-space. For $G \in \D'(E)$ we
have that
$$
G \in \E'(E) \quad \mbox{iff there is} \quad j \in I \quad
\mbox{such that} \quad \Supp G \subset  F(j) \mbox{,}
$$
i.e., the space $\E'(E)$ consists of all distributions with compact
support on the space $V$.
\end{prop}
\proof follows from the forth statement of proposition~\ref{ch} and
formula~(\ref{e'}).
\vspace{0.3cm}

\begin{nt} \label{nt2} {\em
Using formulas~(\ref{d}) -- (\ref{te'}) we obtain that
$$ \E(E) = \tilde{\E}(E) = \D(E) \quad  \mbox{and} \quad \E'(E) = \tilde{\E}'(E) = \D'(E)$$
 for any compact $C_1$-space $E=(I, F,V)$.

Using formulas~(\ref{d}) -- (\ref{te'}) we  obtain that
$$\tilde{\E}'(E)=\E'(E) = \D(E) \quad \mbox{and} \quad \tilde{\E}(E)=\E(E)
= \D'(E)$$
 for any discrete $C_1$-space $E=(I, F, V)$. }
\end{nt}

\begin{nt}
\label{nt3} {\em Let $\{ B_k \}$ be a projective system of vector
spaces with surjective transition maps. Then we have a canonical
nondegenerate pairing between vector spaces
$$\mathop{\Lim_{\gets}}_k B_k  \qquad  \mbox{and} \qquad \mathop{\Lim_{\to}}_k
B_k^* \mbox{.} $$ Using it and formulas~(\ref{e}), (\ref{e'}),
(\ref{te'}), (\ref{te}) we have  canonical nondegenerate pairings
$$ < , >_{\E(E)}  \quad : \quad  \E (E)  \times \E'(E) \lto \cc \mbox{,} $$
$$
< , >_{\tilde{\E}(E)}  \quad : \quad  \tilde{\E}' (E)  \times
\tilde{\E}(E) \lto \cc \mbox{.}
$$
These pairings, being restricted to $\D(E) \times \E'(E)$ and to
$\tilde{\E}' (E) \times \D(E)$  correspondingly, coniside with the
pairings restricted from the pairing between $\D(E)$ and $\D'(E)$
(see statement~\ref{s4} of proposition~\ref{ch}). }
\end{nt}

\subsection{Invariant measures and integral} \label{imi}
Let $E=(I, F, V)$ be a $C_1$-space over the field $\df_q$. For any
$G \in \D'(E)$, for any $a \in V$ we denote by $T_a (G)$ the
following distribution from $\D'(E)$:
$$
T_a(G)(f) \eqdef G(T_{-a}(f)) \mbox{, \quad} f \in \D(E)  \mbox{.}
$$

We define
$$
\D'(E)^V \eqdef \{ G \in \D'(E) \mbox{\quad such that \quad} T_a(G)
=G \mbox{\quad for any \quad} a \in V \} \mbox{.}
$$

We have the following proposition.
\begin{prop} \label{pm}
$\D'(E)^V$ is a $1$-dimensional vector space over the field $\cc$.
\end{prop}
\proof. We note that the space $\D(E)$ is generated by functions
$T_a (\delta_{F(i)})$, $a \in V$, $i \in I$, where the function
$$
\delta_{F(i)}(v) \eqdef \left\{
\begin{array}{lrc}
1  \mbox{,} & \mbox{if} &  v \in F(i) \\
0 \mbox{,} & \mbox{if} & v \notin F(i) \mbox{.}
\end{array}
\right.
$$
Besides, we have for any $i \le j \in I$
$$
\delta_{F(j)} = \sum_{b \in F(j)/F(i)} T_{a_b}(\delta_{F(i)})
\mbox{,}
$$
where  $a_b \in F(j)$ is some lift of $b \in F(j) / F(i)$.

Therefore for fixed element $i \in I$ it is enough to define $G \in
\D'(E)^V$ only on $\delta_{F(i)} \in \D(E)$. The space of these
values is a $1$-dimensional vector space over the field $\cc$.
\vspace{0.3cm}

We denote the $1$-dimensional $\cc$-vector space $\D'(E)^V$ by
$\mu(E)$ and call it the space of {\em invariant measures} on $E$.
For any $\mu \in \mu(E)$ and for any $f \in \D(E)$ we denote
$$
\int\limits_V f \, d \mu \eqdef \mu (f)  \mbox{.}
$$

For any $i \in I$ we denote
$$
\mu (F(i)) \eqdef \mu (\delta_{F(i)}) = \int\limits_V \delta_{F(i)}
\, d \mu  \mbox{.}
$$

We note that the $\cc$-vector space $\ff(V)$ is a $\cc$-algebra,
i.e., for any \linebreak $f, g \in \ff(V)$ we have $f \cdot g (v)
\eqdef f(v) \cdot g(v)$, where $v \in V$. From definitions of
subspaces $ \E(E) \subset \ff(V)$ and $\D(E) \subset \ff(V)$ we have
that
\begin{equation}
\E(E) \cdot \D(E) \subset \D(E) \mbox{.}
\end{equation}
Therefore for fixed $\mu \in \mu(E)$ it is well-defined the
following map
$$
\begin{array}{c}
\I_{\mu} \quad : \quad \E(E) \lto \D'(E) \\
\I_{\mu} (f) (g)  \eqdef \int\limits_V f \cdot g \, d \mu \mbox{,
\quad} f \in \E(E) \mbox{, \quad} g \in \D(E) \mbox{.}
\end{array}
$$

We have the following proposition.
\begin{prop}
Let $\mu \in \mu(E)$, $\mu \ne 0$. There is the following $6$-edges
commutative diagram of embeddings:
\begin{equation} \label{de}
\begin{diagram}
\node[3]{\D'(E)} \\
\node{\E(E)} \arrow{ene,t,J}{\I_{\mu}} \node[4]{\tilde{\E}'(E)}
\arrow{wnw,J} \\
\node{\tilde{\E}(E)}  \arrow{n,J} \node[4]{\E'(E)} \arrow{n,J} \\
\node[3]{\D(E)} \arrow{wnw,J} \arrow{ene,b,J}{\I_{\mu}}
\end{diagram}
\end{equation}
\end{prop}
\proof. Most edges are embeddings which given by
formulas~(\ref{em1})-(\ref{em4}).

We proof that the map $\I_{\mu} : \E(E) \to \D'(E)$ is an embedding.
Indeed, let $f \in \E(E)$, $f \ne 0$. Then there is $v \in V$ such
that $f(v) \ne 0$. By definition of the space $\E(E)$ there is $i
\in I$ such that $f (v + w) = f(v)$ for any $w \in F(i)$. Then
$\I_{\mu}(f)(T_v(\delta_{F(i)})) = f(v) \cdot \mu(F(i)) \ne 0$.
Therefore $\I_{\mu}(f) \ne 0$.

Now we prove that $\I_{\mu} (\D(E)) \subset \E'(E)$. Indeed, let $f
\in \D(E)$, $\Supp f \subset F(j)$ for some $j \in I$. Then
$\I_{\mu} (f) \in \mathop{\Lim\limits_{\gets}}\limits_{i \le j}
\ff(F(j) / F(i))$, because $\I_{\mu}(f)(g) = \I_{\mu} (f) (g
\mid_{F(j)})$ for any $g \in \D(E)$. (Here $g \mid_{F(j)}$ is the
function restricted to $F(j)$ and then extended to $V$ by $0$
outside $F(j)$). Therefore $\I_{\mu}(f) \in \E'(E)$.
\vspace{0.3cm}

\begin{nt} \label{rm4}
{\em Let $f \in \E(E)$, $G \in \D'(E)$. Then we define $f \cdot G
\in D'(E)$ by the following rule:
$$
f \cdot G (g) \eqdef G(f \cdot g) \mbox{, \quad} g \in \D(E)
\mbox{.}
$$
From definitions we have for any $f_1, f_2 \in \E(E)$
$$
\I_{\mu}(f_1 \cdot f_2) = f_1 \cdot \I_{\mu}(f_2)= f_2 \cdot
\I_{\mu}(f_1) \mbox{.}
$$

From definitions we have also for any $a \in V$ and for any $f \in
\E(E)$
$$
\I_{\mu} (T_a(f)) = T_a (\I_{\mu}(f)) \mbox{.}
$$
 }
\end{nt}

\subsection{Fourier transform} \label{FT}
Let $E=(I, F, V)$ be a complete $C_1$-space over the field $\df_q$.
Then $\check{\check{E}} = E$.

Any $a \in V$ defines the function $\psi_a \in
\tilde{\E}(\check{E})$:
$$
\psi_a (w) \eqdef \psi (a(w)) \mbox{, \quad} w \in \check{V}
\mbox{.}
$$

\begin{defin} \label{def2}
We fix $\mu \in \mu(E)$. We define {\it the Fourier transform}
$\F_{\mu} : \D(E) \to \ff(\check{V})$:
$$
\F_{\mu}(f)(w) \eqdef \int\limits_V f \cdot \overline{\psi_w} \, d
\mu  \mbox{, \quad} f \in \D(E)  \mbox{, \quad} w \in \check{V}
\mbox{.}
$$
\end{defin}

\begin{example} \label{ex2}
{\em Let $i \in I$. Then by easy
calculations we obtain that
$$
\F_{\mu} (\delta_{F(i)}) = \mu(F(i)) \cdot \delta_{F(i)^{\bot}}
\mbox{\quad .}
$$
}
\end{example}

\begin{lemma} \label{l1}
Let $a \in V$, $b \in \check{V}$, $f \in
\D(V)$. Then
\begin{enumerate}
\item
$$
\F_{\mu}(T_a(f)) = \psi_a \cdot \F_{\mu}(f)  \mbox{,}
$$
\item
$$
\F_{\mu}(\psi_b \cdot f) = T_{-b} \cdot \F_{\mu}(f) \mbox{.}
$$
\end{enumerate}
\end{lemma}
\proof follows at once from definition of Fourier transform by an
easy calculation.
\vspace{0.3cm}

From example~\ref{ex2} and lemma~\ref{l1} we see that $\F_{\mu} (T_a
(\delta_{F(i)}))  \in \D(\check{E})$ for any $a \in V$, $i \in I$ .
Since the space $\D(E)$ is generated by $T_a (\delta_{F(i)})$, $a
\in V$, $i \in I$, we obtain that $\F_{\mu} (\D(E)) \subset
\D(\check{E})$.

The Fourier transform depends linearly on $\mu \in \mu(E)$.
Therefore the following map $\F$ is well-defined:
\begin{equation} \label{forf}
\begin{array}{c}
\F \quad : \quad \D(E) \otimes_{\cc} \mu(E) \lto \D(\check{E}) \\
\F(f \otimes \mu) \eqdef \F_{\mu}(f) \mbox{, \quad} f \in \D(E)
\mbox{, \quad} \mu \in \mu(E) \mbox{.}
\end{array}
\end{equation}

For any complete $C_1$-space $E=(I,F,V)$ we have a canonical
isomorphism
\begin{equation} \label{fm1}
\mu(E) \otimes_{\cc} \mu(\check{E}) \lto \cc \mbox{,}
\end{equation}
where if $\mu \in \mu(E)$, $\mu \ne 0$, then $\mu \otimes \mu^{-1}
\mapsto 1$, and $\mu^{-1} \in \mu (\check{E}) = \mu (E)^*$,
$\mu^{-1} (F(i)^{\bot}) \eqdef \mu (F(i))^{-1}$ for any $i \in I$.

\begin{prop} \label{propf}
 Let $\mu \in \mu(E)$, $\mu \ne 0$. Then
\begin{enumerate}
\item
$\F_{\mu^{-1}} \circ \F_{\mu}(f) = \check{f}$ for any $f \in \D(E)$;
\item
$\F_{\mu}$ is an isomorphism of $\cc$-vector spaces $\D(E)$ and
$\D(\check{E})$.
\end{enumerate}
\end{prop}
\proof. The first statement follows by direct calculations on the
functions $T_a (\delta_{F(i)})$, $a \in V$, $i \in I$ by using
example~\ref{ex2} and lemma~\ref{l1}, because the space $\D(E)$ is
generated by the functions $T_a (\delta_{F(i)})$, $a \in V$, $i \in
I$.

The second statement of this proposition follows from the first
statement.
\vspace{0.3cm}

\begin{defin} \label{def3}
We fix $\mu \in \mu(E)$ such that $\mu \ne 0$. We define the Fourier
transform $\F_{\mu^{-1}} : \D'(E) \to \D'(\check{E})$ in the
following way:
$$
\F_{\mu^{-1}} (G)(f) \eqdef G ( \F_{\mu^{-1}} (f) ) \mbox{, \quad} G
\in \D'(E) \mbox{, \quad} f \in \D(\check{E}) \mbox{.}
$$
\end{defin}

From definitions~\ref{def2} and~\ref{def3} we have that
$$
\F_{\mu^{-1}} \circ \I_{\mu} (f) = \I_{\mu^{-1}} \circ \F_{\mu}(f)
\mbox{, \quad} f \in \D(E) \mbox{.}
$$

From definition~\ref{def3} and formulas~(\ref{forf}) we have that
the following map $\F$ is well-defined:
$$
\begin{array}{c}
\F \quad : \quad \D'(E) \otimes_{\cc} \mu(\check{E}) \lto \D'(\check{E}) \\
\F(G \otimes \mu^{-1}) \eqdef \F_{\mu^{-1}}(G) \mbox{, \quad} G \in
\D(E) \mbox{, \quad} \mu \in \mu(E) \mbox{, \quad} \mu \ne 0
\mbox{.}
\end{array}
$$

\begin{lemma} \label{l2}
Let $G \in \D'(E)$, $a \in V$, $b \in \check{V}$. Then
\begin{enumerate}
\item
$$
\F_{\mu^{-1}} (T_a (G)) = \psi_a \cdot \F_{\mu^{-1}} (G) \mbox{,}
$$
\item
$$
\F_{\mu^{-1}} (\psi_b \cdot G) = T_{-b} (\F_{\mu^{-1}} (G)) \mbox{.}
$$
\end{enumerate}
\end{lemma}
\proof. Using lemma~\ref{l1}, we have for any $f \in \D(\check{E})$:
 $$
 \begin{array}{c}
  \F_{\mu^{-1}} (T_a (G))(f) = (T_a (G))
(\F_{\mu^{-1}} (f) ) = G (T_{-a}(\F_{\mu^{-1}}(f))) = \\
= G (\F_{\mu^{-1}} (\psi_a \cdot f)) = \F_{\mu^{-1}} (G) (\psi_a
\cdot f) = (\psi_a \cdot \F_{\mu^{-1}} (G)) (f) \mbox{.}
\end{array}
$$
$$
\begin{array}{c}
\F_{\mu^{-1}}(\psi_b \cdot G)(f) = (\psi_b \cdot G)
(\F_{\mu^{-1}}(f)) = G (\psi_b \cdot \F_{\mu^{-1}}(f)) = \\
= G (\F_{\mu^{-1}} (T_b (f))) = \F_{\mu^{-1}} (G) (T_b (f)) = T_{-b}
(\F_{\mu^{-1}}(G(f))) \mbox{.}
\end{array}
$$
\vspace{0.3cm}

 For any $G \in \D'(E)$ we define $\check{G} \in \D'(E)$:
$$
\check{G} (f) \eqdef G (\check{f}) \mbox{, \quad} f \in \D(E)
\mbox{.}
$$
We note that for any $g \in \E(E)$ we have
$$
\I_{\mu} (\check{g}) = \I_{\mu}(g)\check{\vphantom{\I}} \mbox{\quad
.}
$$

\begin{prop} \label{prop9}
Let $\mu \in \mu(E)$, $\mu \ne 0$. Then
\begin{enumerate}
\item
$\F_{\mu} \circ \F_{\mu^{-1}}(G) = \check{G}$ for any $G \in
\D'(E)$;
\item
$\F_{\mu^{-1}}$ is an isomorphism of $\cc$-vector spaces $\D'(E)$
and $\D'(\check{E})$.
\end{enumerate}
\end{prop}
\proof follows at once from proposition~\ref{propf}.
\vspace{0.3cm}

Now we define the Fourier transform $\F \; : \; \E(E) \lto
\D'(\check{E})$ by the following rule:
\begin{equation}  \label{four'e}
\F(g)  \eqdef \F_{\mu^{-1}} \circ \I_{\mu} (g) \mbox{, \quad} g \in
\E(E) \mbox{ .}
\end{equation}
We note that this map $\F$ does not depend on the choice of $\mu \in
\mu(E)$, although we use $\mu \in \mu(E)$, $\mu \ne 0$ for the
definition of this map $\F$.

\begin{nt} \label{nt5}
{\em Let $i \le j \in I$, then in
section~\ref{0-case} we defined the Foutier
transform
$$
\F \quad : \quad \ff(F(j)/F(i)) \lto \ff((F(j)/F(i))^*)  \mbox{.}
$$

Now using successively  projective and inductive limits from
formulas~(\ref{d})-(\ref{te'}) and explicit description of
$\check{V}$ from subsection~\ref{c_1-spaces}, by means of
statement~\ref{stf4} of proposition~\ref{f0} we redefine the Fourier
transform  $\F_{\mu^{-1}} : \D'(E) \to \D'(\check{E})$. This new
definition of Fourier transform coincides with the definition of
Fourier transform which was given in this section above. }
\end{nt}

Using this definition of Fourier transform, we have the following
proposition.
\begin{prop} \label{four'e'}
The following statements are satisfied.
\begin{enumerate}
\item \label{stnf1} $\F$ gives an isomorphism of $\cc$-vector spaces $\E(E)$ and
$\tilde{\E}'(\check{E})$.
\item  \label{stnf2}  $\F \mid_{\tilde{\E}(E)}$ gives an isomorphism of $\cc$-vector
spaces $\tilde{\E}(E)$ and $\E'(\check{E})$.
\item \label{stnf3}  For any $f \in \E(E)$, for any $g \in \tilde{\E}(\check{E})$ we
have
\begin{equation} \label{ff}
<\F(f), g>_{\tilde{\E}(\check{E})} = <f, \F(g)>_{\E(E)}  \mbox{.}
\end{equation}
\item  \label{stnf4}  For any $G \in \tilde{\E}'(\check{E})$, for any $v \in V$ we
have
\begin{equation} \label{fi}
\F^{-1}(G)(v) = <G, \psi_v>_{\tilde{\E}(\check{E})} \mbox{.}
\end{equation}
\end{enumerate}
\end{prop}
\proof. Statements~\ref{stnf1}, \ref{stnf2} follow from
remark~\ref{nt5} and formulas~(\ref{e})-(\ref{te'}).
Statement~\ref{stnf3} follows from statement~\ref{stf3} of
proposition~\ref{f0}, statements~\ref{stnf1}-\ref{stnf2} of this
proposition and remark~\ref{nt3}.

Statement~\ref{stnf4} follows from statement~\ref{stnf3} in the
following way. For any $v \in V$ we have
$$ \delta_v \in \E'(E) \quad : \quad \delta_v(g) \eqdef g(v) \mbox{,
\quad} g \in \D(E) \mbox{.}
$$
We have $<f, \delta_v>_{\E(E)} = f(v)$ for any $f \in \E(E)$. We
have also  $\F(\psi_v) = T_{-v}(\F(1)) = T_{-v}(\delta_0)=
\delta_v$. Now formula~(\ref{fi}) follows from formula~(\ref{ff}).
\vspace{0.3cm}

\subsection{Poisson formula for 1-dimensional case.} \label{PF}
Let
$$
0 \arrow{e} E_1  \arrow{e,t}{\alpha}  E_2 \arrow{e,t}{\beta}   E_3
\arrow{e} 0
$$
be an admissible triple of complete $C_1$-spaces over the field
$\df_q$, where \linebreak $E_i = (I_i, F_i, V_i)$, $ 1 \le i \le 3$.
Then we have the following admissible triple of complete
$C_1$-spaces:
$$
0 \arrow{e} \check{E_3} \arrow{e,t}{\check{\beta}} \check{E_2}
\arrow{e,t}{\check{\alpha}}  \check{E_1} \arrow{e} 0 \mbox{.}
$$

Besides, we have the following canonical isomorphism:
\begin{equation}  \label{fm2}
\mu(E_1) \otimes_{\cc} \mu(E_3) \lto \mu(E_2) \mbox{,}
\end{equation}
where $ (\mu_1 \otimes \mu_3) (F_2(j)) \eqdef \mu_1 (F_2(j) \cap V_1
) \cdot \mu_3 (\beta (F_2(j))) $, $j \in I_2$, $\mu_1 \in \mu(V_1)$,
$\mu_3 \in \mu(V_3)$.

We define the map $\alpha^*$:
\begin{equation} \label{ff1}
\alpha^* \quad : \quad \D(E_2) \lto \D(E_1)  \mbox{,}
\end{equation}
where $\alpha^* (f)(v) \eqdef f (\alpha(v)) $, $f \in \D(E_2)$, $v
\in V_1$. This map is well-defined.

As conjugate map to the map $\alpha^*$ we define the map $\alpha_*$:
\begin{equation}  \label{ff2}
\alpha_* \quad : \quad \D'(E_1) \lto \D'(E_2) \mbox{.}
\end{equation}

We note that if we fixed $\mu_1 \in \mu(E_1)$, $\mu_1 \ne 0$, and
fixed $\mu_2 \in \mu(E_2)$, $\mu_2 \ne 0$, then from
formulas~(\ref{fm1}) and (\ref{fm2}) we have well-defined $\mu_1
\otimes \mu_2^{-1} \in \mu(\check{E}_3)$.

We define {\em the characteristic function } of subspace  $E_1$ as
distribution   $\delta_{E_1, \mu_1} \in \D'(E_2)$, where
$$
  \delta_{E_1,
\mu_1} (f) \eqdef \int\limits_{V_1} \alpha^*(f) \, d \mu_1 \quad
\mbox{for any} \quad f \in \D(E_2)\mbox{.}
$$
Analogously we can define  $\delta_{\check{E}_3, \mu_1 \otimes
\mu_2^{-1}} \in \D'(\check{E}_2)$.

\begin{Th}[\it Poisson formula].
For any $f \in \D(E_2)$, for any  $\mu_1 \in \mu(E_1)$, $\mu_1 \ne
0$, and $\mu_2 \in \mu(E_2)$, $\mu_2 \ne 0$  the following
equivalent statements are satisfied.
\begin{enumerate}
\item
\begin{equation} \label{fp}
\int\limits_{V_1} \alpha^*(f) \, d \mu_1 =\int\limits_{\check{V_3}}
\check{\betta}^* (\F_{\mu_2} (f)) \, d(\mu_1 \otimes \mu_2^{-1})
\quad \mbox{.}
\end{equation}

\item
\begin{equation} \label{fp'}
\F_{\mu_2^{-1}} (\delta_{E_1, \mu_1}) =\delta_{\check{E}_3, \mu_1
\otimes \mu_2^{-1}} \quad \mbox{.}
\end{equation}
\end{enumerate}
\end{Th}
\proof.  We have the following obvious properties:
\begin{equation} \label{eqn1}
T_a (\delta_{E_1, \mu_1}) = \delta_{E_1, \mu_1} \quad \mbox{for any}
\quad a \in  V_1   \quad \mbox{;}
\end{equation}
\begin{equation} \label{eqn2}
\psi_{b} \cdot \delta_{E_1, \mu_1} = \delta_{E_1, \mu_1} \quad
\mbox{for any} \quad b \in  \check{V}_3 \quad \mbox{.}
\end{equation}
There are analogous properties for $\delta_{\check{E}_3, \mu_1
\otimes \mu_2^{-1}}$.

By lemma~\ref{l2} we have that properties given by
formulas~(\ref{eqn1})-(\ref{eqn2}) are changed with each other under
the Fourier transform.  We show that $\delta_{E_1, \mu_1}$ is
determined by these properties uniquely up to some constant. (And
the analogy  is true for $\delta_{\check{E}_3, \mu_1 \otimes
\mu_2^{-1}}$.) Indeed, from formula~(\ref{eqn2}) it follows that
$\delta_{E_1, \mu_1} = \alpha_* (G) $ for some $G \in \D'(E_1)$. Now
from formula~(\ref{eqn1}) and proposition~\ref{pm} we obtain that
$G(g) = c \cdot \int\limits_{V_1} g \, d \mu_1 $ for some $c \in
\cc$ and any $g \in \D(E_1)$.

Therefore we proved formula~(\ref{fp'}) up to some constant $d \in
\cc$. To obtain that $d =1$ we calculate formula~(\ref{fp'}) on
$\delta_{F_2(j)} \in \D(E_2)$ for some $j \in I_2$ using explicit
example~\ref{ex2}.
\vspace{0.3cm}

\subsection{Direct and inverse images.} \label{DI}
Let
$$
0 \arrow{e} E_1  \arrow{e,t}{\alpha}  E_2 \arrow{e,t}{\beta}   E_3
\arrow{e} 0
$$
be an admissible triple of $C_1$-spaces over the field $\df_q$,
where $E_i = (I_i, F_i, V_i)$, $ 1 \le i \le 3$.

Now we define the map $\beta_*$:
\begin{equation} \label{ff3}
\betta_*  \quad :  \quad \D(E_2) \otimes_{\cc} \mu(E_1) \lto \D(E_3)
\mbox{,}
\end{equation}
where
\begin{equation} \label{ff3'}
\beta_* (f \otimes \mu_1)(w) \eqdef \int\limits_{\betta^{-1}(w)} f
|_{ \betta^{-1}(w)} d \mu_1  \eqdef \int\limits_{V_1} T_{w'}(f) d
\mu_1  =  \mu_1(T_{w'}(f)|_{V_1}).
\end{equation}
Here $f \in \D(E_2)$, $\mu_1 \in \mu(E_1)$,  $w' \in \betta^{-1}(w)$
and then $\betta^{-1}(w) = w' + V_1$. The result doesn't depend on a
choice of the $w'$.

As conjugate map to the map $\beta_*$ we define the map $\beta^*$:
\begin{equation}  \label{ff4}
\beta^* \quad : \quad \D'(E_3) \otimes_{\cc} \mu(E_1)  \lto \D'(E_2)
\mbox{.}
\end{equation}

We have the following property of the map $\beta_*$.
\begin{prop}[\it Fubini formula] \label{fubf}
Let $f \in \D(E_2)$, $\mu_i \in \mu(E_i)$ for $1 \le i \le 3$. Let
$\mu_1 \otimes \mu_3 = \mu_2$ under the map~(\ref{fm2}). Then
$$
\int\limits_{E_3} \beta_* (f \otimes \mu_1) d \mu_3 =
\int\limits_{E_2} f d \mu_2 \mbox{.}
$$
\end{prop}
\proof. We know that the space $\D(E_2)$ is generated by functions
$T_a (\delta_{F_2(i)})$, where $a \in V_2$, $i \in I_2$. Therefore
we check the Fubini formula on these functions.
\vspace{0.3cm}

Now we extend the direct and inverse images  of admissible
monomorphisms and epimorphisms (see formulas~(\ref{ff1}),
(\ref{ff2}), (\ref{ff3}), (\ref{ff4})) to other types of functions
and distributions from diagram~(\ref{de}). The problem is that these
maps are not defined for arbitrary admissible monomorphisms and
epimorphisms.

\begin{prop}. \label{if}
The following maps of spaces of functions on  $C_1$-spaces are
well-defined.

$$
\begin{array}{lllll}
1. & \alpha^* \quad :& \quad \D(E_2)  \lto \D(E_1) \mbox{.}
&&\\[0.3cm]

2. &  \beta^* \quad :& \quad \D(E_3) \lto \D(E_2) \mbox{,} \quad
\mbox{where} & E_1 ~\mbox{is  compact} \mbox{,} & V_1 =
\mbox{Ker}(\beta)\mbox{.}\\[0.3cm]

3. &   \alpha_* \quad :& \quad \D(E_1) \lto \D(E_2)  \mbox{,} \quad
\mbox{where}&  E_3 ~\mbox{is discrete}\mbox{,} & V_3 =
\mbox{Coker}(\alpha)\mbox{.}\\[0.3cm]

4. & \beta_* \quad :& \quad \D(E_2) \otimes_{\cc} \mu(E_1) \lto
\D(E_3)
\mbox{.}&&\\[0.3cm]

5. & \alpha^* \quad : & \quad {\E}(E_2) \lto
{\E}(E_1) \mbox{,}& &\\[0.2cm]

&& \mbox{and } \; \alpha^* (\tilde{\E}(E_2))  \subset
\tilde{\E}(E_1)
\mbox{.}&&\\[0.3cm]

6. &  \beta^* \quad : & \quad {\E}(E_3)
\lto {\E}(E_2) ~\mbox{,}&&\\[0.2cm]

&& \mbox{and } \; \beta^* (\tilde{\E}(E_3))  \subset \tilde{\E}(E_2)
\mbox{.}&&\\[0.3cm]

7. & \alpha_* \quad : & \quad {\E}(E_1) \lto {\E}(E_2)\mbox{,} \quad
\mbox{where}  & E_3
~\mbox{is discrete}\mbox{,} & V_3 = \mbox{Coker}(\alpha) ~\mbox{,}\\[0.2cm]

&& \mbox{and } \; \alpha_* (\tilde{\E}(E_1))  \subset
\tilde{\E}(E_2)
\mbox{.}\\[0.3cm]

8. &  \beta_* \quad : & \quad {\E}(E_2) \otimes_{\cc} \mu(E_1) \lto
{\E}(E_3)\mbox{,} \quad \mbox{where} & E_1
~\mbox{is compact}\mbox{,} & V_1 = \mbox{Ker}(\beta)~\mbox{,}\\[0.2cm]

 && \mbox{and } \; \beta_* (\tilde{\E}(E_2) \otimes_{\cc} \mu(E_1))
\subset \tilde{\E}(E_3) \mbox{.}&&
\end{array}
$$
\end{prop}
\proof.  For any map $\gamma: V \to W$ we have the map $\gamma^*:
\ff(W) \to \ff(V)$ by the rule: $\gamma^* (f)(v) \eqdef f
(\gamma(v))$, where $f \in \ff(W)$, $v \in V$.

For any embedding $\lambda : Y \lto Z$ we have the map $\lambda_*:
\ff(Y) \lto \ff(Z)$ by the following rule:
$$
\lambda_* (f)(z) = \left\{
\begin{array}{lrc}
0  & \mbox{, if} & z \notin \lambda(Y) \\
  f(y)  & \mbox{, if} & z = \lambda(y) \mbox{,}
\end{array}
\right.
$$
 where $ f \in \ff(Y)$, $z \in Z $, $y \in Y$.

 According to remark~\ref{nt2} we have that $\E(E) = \D(E)$ for any
 compact $C_1$-space $E=(I,F,V)$. Therefore for the compact $C_1$-space $E_1$
 it is well defined the map
 $$
\beta_* \quad : \quad \E(E_2) \otimes_{\cc} \mu(E_1) \lto \E(E_3)
 $$
by formula~(\ref{ff3'}). Now using proposition~\ref{ch} we verify
the statements of this proposition. The proposition is proved
\vspace{0.3cm}

Now as conjugate maps to the maps described in proposition~\ref{if}
we define the direct and inverse images of admissible monomorphisms
and epimorphisms for the spaces of all distributions $\D'(E_i)$, $1
\le i \le 3$ when these maps are defined. To calculate the
corresponding direct and inverse images for the other types of
distributions from diagram~(\ref{de}) we have to introduce the
topology on the spaces of functions $\E(E)$ and $\tilde{\E}(E)$ for
a $C_1$-space $E=(I, F, V)$.

Let $\{B_l\}$ be a projective system of topological $k$-vector
spaces with surjective transition maps (over a discrete field $k$).
Then $\mathop{\mathop{\Lim\limits_{\gets}}}\limits_l B_l $ is a
topological $k$-vector space with a base of neighbourhoods of $0$
given by $\phi_l^{-1}(V)$ for all $l$ and all open $k$-subspaces $V$
of $B_l$, where $\phi_l :
\mathop{\mathop{\Lim\limits_{\gets}}}\limits_l B_l \to B_l$ is the
natural map. By $\Hom_k^{cont}(\cdot, \cdot )$ we denote a $k$-space
of all $k$-linear continuous maps between topological $k$-vector
spaces. Then we have a canonical isomorphism:
\begin{equation} \label{ci1}
\Hom\nolimits_k^{cont}(
\mathop{\mathop{\Lim\limits_{\gets}}}\limits_l B_l,k) =
\mathop{\mathop{\Lim\limits_{\to}}}\limits_l \Hom\nolimits_k^{cont}
(B_l, k) \mbox{.}
\end{equation}

Let $\{A_m\}$ be an inductive system of topological $k$-vector
spaces with injective transition maps  (over a discrete field $k$).
Then $\mathop{\mathop{\Lim\limits_{\to}}}\limits_m A_m $ is a
topological $k$-vector space with a base of neighbourhoods of zero
given by all $k$-subspaces $U \subset
\mathop{\mathop{\Lim\limits_{\to}}}\limits_m A_m$ such that $U \cap
A_m$ is an open $k$-subspace in $A_m$ for any $m$. Then we have a
canonical isomorphism:
\begin{equation} \label{ci2}
\Hom\nolimits_k^{cont}( \mathop{\mathop{\Lim\limits_{\to}}}\limits_m
A_m,k) = \mathop{\mathop{\Lim\limits_{\gets}}}\limits_m
\Hom\nolimits_k^{cont} (A_m, k) \mbox{.}
\end{equation}

Now we apply these reasonings to the field $k = \cc$ (with a
discrete topology) and the $\cc$-vector spaces $\tilde{\E}(E)$ and
$\E(E)$. We use formulas~(\ref{te}) and~(\ref{e}). We consider the
discrete topology on every finite-dimensional $\cc$-vector space
$\ff(F(i)/F(j))$ for any $j \le i \in I$.

Then $\E(E)$ is a topological $\cc$-vector space with the base of
neighbourhoods of $0$ given by the following $\cc$-subspaces
$\E(E)_i$ (for all $i \in I$):
$$
\E(E)_i  \eqdef \{f \in \E(E) \quad \mbox{such that} \quad f
\mid_{F(i)} = 0  \} \mbox{.}
$$

Let for any $j \in I$
$$\tilde{\E}(E)_j \eqdef \{f \in  \tilde{\E}(E) \quad \mbox{such that} \quad T_{w}f =f
\quad \mbox{for any} \quad w \in F(j) \} \mbox{.} $$ Then
$\tilde{\E}(E) = \bigcup\limits_{j \in I}\tilde{\E}(E)_j$, and
$\tilde{\E}(E)$ is a topological $\cc$-vector space with the base of
neighbourhoods of zero given by all $\cc$-subspaces $U \subset
\tilde{\E}(E)$ such that for any $j \in I$ there is $i \in I$, $i
\ge j$ such that
$$
U \cap \tilde{\E}(E)_j \quad \supset \quad  \E(E)_i \cap
\tilde{\E}(E)_j  \mbox{.}
$$

We have that the canonical embedding $\tilde{\E}(E) \subset \E(E)$
is a continuous map with respect to these topologies, but the
restriction of the topology of $\E(E)$ to $\tilde{\E}(E)$ doesn't
coincide with the topology of $\tilde{\E}(E)$.

Using formulas~(\ref{e'}) and~(\ref{te'}) for the $\cc$-spaces
$\E'(E)$ and $\tilde{\E}'(E)$, and also isomorphisms~(\ref{ci1})
and~(\ref{ci2}) we obtain the following canonical isomorphisms
\begin{equation} \label{for}
\qquad \E'(E) = \Hom\nolimits_{\cc}^{cont}(\E(E), \cc) \qquad
\mbox{and} \qquad \tilde{\E}'(E) =
\Hom\nolimits_{\cc}^{cont}(\tilde{\E}(E), \cc) \mbox{.}
\end{equation}

We have the following lemma.
\begin{lemma} \label{lll}
Suppose that conditions of proposition~\ref{if} are satisfied. Then
\begin{enumerate}
\item
The  maps  $\alpha_*$, $\alpha^*$, $\beta_*$, $\beta^*$ (when they
are defined) are continuous on the topological $\cc$-vector spaces
$\E(E_i)$, $1 \le i \le 3$.
\item
The  maps $\alpha_*$, $\alpha^*$, $\beta_*$, $\beta^*$ (when they
are defined) are continuous on the topological $\cc$-vector spaces
$\tilde{\E}(E_i)$, $1 \le i \le 3$.
\end{enumerate}
\end{lemma}
\proof. We use  explicit descriptions of open subsets on $\E(E)$ and
$\tilde{\E}(E)$ given above. For example, to prove that $\beta_*$ is
a continuous map from $\cc$-vector space $\E(E_2)$  to $\cc$-vector
space $\E(E_3)$ (or from $\cc$-vector space $\tilde{\E}(E_2)$ to
$\cc$-vector space $\tilde{\E}(E_3)$) when $E_1$ is a compact
$C_1$-space, we use that
$$
\beta_*( \E(E_2)_i \otimes_{\cc} \mu(E_1)) \quad \subset \quad
\E(E_3)_{i'}
$$
when $i \in I_2$, $F_2(i) \supset V_1$, and $i' \in I_3$, $F_3(i') =\beta(F_2(i))$. And also
$$
\beta_*(\tilde{\E}(E_2)_j \otimes_{\cc} \mu(E_1)) \quad \subset
\quad \tilde{\E}(E_3)_{j'}
$$
when $j \in I_2$, $j' \in I_3$, and $F_3(j')=\beta(F_2(j))$.
\vspace{0.3cm}

Using this lemma and formulas~(\ref{for}) we consider the direct and
inverse images as conjugate maps to the maps considered in
proposition~\ref{if}.

\begin{prop}. \label{if'}
The following maps of spaces of distributions on $C_1$-spaces are
well-defined.
$$
\begin{array}{lllll}
1. & \alpha_* \quad :& \quad \D'(E_1)  \lto \D'(E_2) \mbox{.} &&\\[0.3cm]

2. &  \beta_* \quad :& \quad \D'(E_2) \lto \D'(E_3) \mbox{,} \quad
\mbox{where} & E_1 ~\mbox{is compact} \mbox{,} & V_1 =
\mbox{Ker}(\beta)\mbox{.}\\[0.3cm]

3. &   \alpha^* \quad :& \quad \D'(E_2) \lto \D'(E_1) \mbox{,} \quad
\mbox{where}& E_3 ~\mbox{is discrete}\mbox{,} & V_3 =
\mbox{Coker}(\alpha)\mbox{.}\\[0.3cm]

4. & \beta^* \quad :& \quad \D'(E_3) \otimes_{\cc} \mu(E_1) \lto
\D'(E_2)
\mbox{.}&&\\[0.3cm]

5. & \alpha_* \quad : & \quad \tilde{\E}'(E_1) \lto
\tilde{\E}'(E_2) ~\mbox{,}& &\\[0.2cm]

&& \mbox{and } \; \alpha_* (\E'(E_1))  \subset \E'(E_2)
\mbox{.}&&\\[0.3cm]

6. &  \beta_* \quad : & \quad \tilde{\E}'(E_2)
\lto \tilde{\E}'(E_3) ~\mbox{,}&&\\[0.2cm]

&& \mbox{and } \; \beta_* (\E'(E_2))  \subset \E'(E_3)
\mbox{.}&&\\[0.3cm]

7. & \alpha^* \quad : & \quad \tilde{\E}'(E_2) \lto
\tilde{\E}'(E_1)\mbox{,} \quad  \mbox{where}  & E_3
~\mbox{is discrete}\mbox{,} & V_3 = \mbox{Coker}(\alpha) ~\mbox{,}\\[0.3cm]

&& \mbox{and } \; \alpha^* (\E'(E_2))  \subset \E'(E_1)
\mbox{.}\\[0.3cm]

8. &  \beta^* \quad : & \quad \tilde{\E}'(E_3) \otimes_{\cc}
\mu(E_1) \lto \tilde{\E}'(E_2)\mbox{,} \quad  \mbox{where} & E_1
~\mbox{is compact}\mbox{,} & V_1 = \mbox{Ker}(\beta)~\mbox{,}\\[0.2cm]

 && \mbox{and } \; \beta^* (\E'(E_3) \otimes_{\cc} \mu(E_1))
\subset \E'(E_2) \mbox{.}&&
\end{array}
$$
\end{prop}
\proof. We use proposition~\ref{if}, lemma~\ref{lll}, and
formulas~(\ref{for}).
\vspace{0.3cm}

\begin{nt} { \em
Using remark~\ref{nt3} we obtain that conjugate maps defined in
proposition~\ref{if'} are compatible with respect to the canonical
embeddings:
$$
\E'(E_i)  \subset \tilde{\E}'(E_i) \subset \D'(E_i) ,  \qquad 1 \le
i \le 3 \mbox{.}
$$
}
\end{nt}

From definitions of the maps $\alpha^*$, $\alpha_*$, $\beta^*$,
$\beta_*$ we have the following {\em projections formulas}. We fix
$\mu_1 \in \mu(E_1)$.

For any $f \in \D(E_2)$ and  $g \in \E(E_3)$
\begin{equation} \label{pf1}
\beta_*(f \cdot \beta^*(g) \otimes \mu_1) = \beta_*(f \otimes \mu_1)
\cdot g \mbox{.}
\end{equation}

For any $f \in \E(E_2)$, $g \in \E(E_3)$, and a compact $C_1$-space
$E_1$ we have
\begin{equation} \label{pf1'}
\beta_*(f \cdot \beta^*(g) \otimes \mu_1) = \beta_*(f \otimes \mu_1)
\cdot g \mbox{.}
\end{equation}

For any $ f \in \E(E_1)$, $g \in \E(E_2)$, and a discrete
$C_1$-space $E_3$
\begin{equation} \label{pf2}
\alpha_* (f \cdot \alpha^* (g)) = \alpha_*(f) \cdot g  \mbox{.}
\end{equation}
\begin{prop} \label{prop15}
 Let $\mu_i \in \mu(E_i)$, $1 \le i \le 3$ such that
$\mu_1 \otimes \mu_3 = \mu_2$ according to map~(\ref{fm2}).
 We have the following formulas.
\begin{enumerate}
\item $\I_{\mu_2}(\beta^*(g)) = \beta^*(\I_{\mu_3}(g))$ for any $g \in
\E(E_3)$.
\item Let $E_1$ be a compact $C_1$-space. Then
$\I_{\mu_3}(\beta_*(f \otimes \mu_1)) = \beta_*(\I_{\mu_2}(f))$ for
any $f \in \E(E_2)$.
\end{enumerate}
\end{prop}
\proof follows from Fubini formula (proposition~\ref{fubf}) and
projection formulas~(\ref{pf1}) and~(\ref{pf1'}).
\vspace{0.3cm}

\subsection{Composition of maps and base change rules.} \label{bcr}
Let
$$
0 \arrow{e} E_1  \arrow{e,t}{\alpha}  E_2 \arrow{e,t}{\beta}   E_3
\arrow{e} 0
$$
be an admissible triple of $C_1$-spaces over the field $\df_q$,
where $E_i = (I_i, F_i, V_i)$, $ 1 \le i \le 3$.

 Let
$$
0 \arrow{e} L  \arrow{e,t}{\alpha'}  H \arrow{e,t}{\beta'}   E_2
\arrow{e} 0
$$
be another admissible triple of $C_1$-spaces over the field $\df_q$,
where $L = (J_1, T_1, U)$,  $H= (J_2, T_2, W)$.

We have then the following admissible triple of $C_1$-spaces:
$$
0 \arrow{e} H \mathop{\times}_{E_2} E_1 \arrow{e,t}{\beta'_{\alpha}}
H \arrow{e,t}{\beta \beta'} E_3 \arrow{e} 0  \mbox{.}
$$
 Let $H' = H \mathop{\times}\limits_{E_2} E_1 =(K, P, Z)$ as
 $C_1$-space.

\begin{prop}. \label{pr16} We have the following formulas.
\begin{enumerate}
\item For any $f \in \D(H)$, $\nu \in \mu(L)$, $\mu \in
\mu(E_1)$
\begin{equation}
\label{e1} (\beta \beta')_* (f \otimes (\nu \otimes \mu)) =  \beta_*
(\beta'_* (f \otimes \nu) \otimes \mu) \mbox{.}
\end{equation}
\item For any $G \in \D'(E_3)$, $\nu \in \mu(L)$, $\mu \in
\mu(E_1)$
\begin{equation}
\label{e2} (\beta \beta')^* (G \otimes (\nu \otimes \mu)) =
(\beta')^*  (\beta^*  (G \otimes \mu) \otimes \nu)  \mbox{.}
\end{equation}
\item For any $f \in \E(E_3)$
\begin{equation}
\label{e3} (\beta \beta')^* (f) =  (\beta')^*  \beta^* (f) \mbox{.}
\end{equation}
\item For any $G \in \tilde{\E}'(H)$
\begin{equation}
\label{e4}
 (\beta \beta')_* (G) =  \beta_* \beta'_* (G)
\mbox{.}
\end{equation}

\vspace{0.3cm}

\noindent We suppose further that  $E_1$ and $L$ are compact
$C_1$-spaces. Then $H \mathop{\times}\limits_{E_2} E_1$ is a compact
 $C_1$-space, and the following formulas are satisfied.

\item For any $f \in \D(E_3)$
\begin{equation} \label{pr16f5}
(\beta \beta')^*(f) = (\beta')^* \beta^* (f) \mbox{.}
\end{equation}
\item For any $G \in \D'(H)$
\begin{equation} \label{pr16f6}
(\beta \beta')_*(G)= \beta_* (\beta')_* (G) \mbox{.}
\end{equation}
\item For any $f \in \E(H)$, $\nu \in \mu(L)$, $\mu \in
\mu(E_1)$
\begin{equation} \label{pr16f7}
 (\beta \beta')_* (f \otimes (\nu \otimes \mu)) =  \beta_*
(\beta'_* (f \otimes \nu) \otimes \mu) \mbox{.}
\end{equation}
\item For any $G \in \tilde{\E}'(E_3)$, $\nu \in \mu(L)$, $\mu \in
\mu(E_1)$
\begin{equation} \label{pr16f8}
 (\beta \beta')^* (G \otimes (\nu \otimes \mu)) = (\beta')^*
(\beta^*  (G \otimes \mu) \otimes \nu)  \mbox{.}
\end{equation}
\end{enumerate}
\end{prop}
\proof. We have by formula~(\ref{fm2}) that the following equalities
are canonically satisfied:
$$
\mu(E_1) \otimes_{\cc} \mu(E_3) = \mu(E_2 ) \qquad \mbox{,} \qquad
\mu(L) \otimes_{\cc} \mu(E_2) = \mu (H) \mbox{,}
$$
$$
\mu(H') \otimes_{\cc} \mu(E_3) = \mu(H)  \mbox{.}
$$
Therefore we have that
$$
\mu(L) \otimes_{\cc} \mu(E_1) = \mu(H')  \mbox{.}
$$
The last formula follows also from the following admissible triple
of $C_1$-spaces:
\begin{equation} \label{rastr}
0 \arrow{e} L \arrow{e,t}{\alpha_{\alpha'}}
 H \mathop{\times}_{E_2} E_1
 \arrow{e,t}{\alpha_{\beta'}} E_1 \arrow{e} 0  \mbox{.}
\end{equation}
Therefore for  $\nu \in \mu(L)$, $\mu \in \mu(E_1)$ we have
canonically  $\nu \otimes \mu \in \mu \left(H
\mathop{\times}\limits_{E_2} E_1\right)$.

Now formula~(\ref{e1}) follows from Fubini formula
(proposition~\ref{fubf}). Formula~(\ref{e2}) is the conjugate
formula to formula~(\ref{e1}). Formula~(\ref{e3}) is evident. We
restrict formula~(\ref{e3}) to the subspace~$\tilde{\E}(E_3)$. Then
formula~(\ref{e4}) is the conjugate formula to the last formula.

From triple~(\ref{rastr}) we obtain that if  $E_1$ and $L$ are
compact  $C_1$-spaces, then $H \mathop{\times}\limits_{E_2} E_1$ is
also a compact  $C_1$-space.

Now formula~(\ref{pr16f5}) follows from the definitions of inverse
images. Formula~(\ref{pr16f6}) is the conjugate formula to
formula~(\ref{pr16f5}). Formula~(\ref{pr16f7}) follows from Fubini
formula (proposition~\ref{fubf}).  Formula~(\ref{pr16f8}) follows
from formula~(\ref{e2}). The proof is finished.
\vspace{0.3cm}

Let
$$
0 \arrow{e} E_1  \arrow{e,t}{\alpha}  E_2 \arrow{e,t}{\beta}   E_3
\arrow{e} 0
$$
be an admissible triple of $C_1$-spaces over the field $\df_q$,
where $E_i = (I_i, F_i, V_i)$, $ 1 \le i \le 3$. Let $\alpha' : E_2
\lto H'$ be an admissible monomorphism, where $H'= (J', T', W')$,
that is, there is an admissible triple of $C_1$-spaces
$$
0 \arrow{e} E_2  \arrow{e,t}{\alpha'}  H' \arrow{e,t}{\beta'}   L'
\arrow{e} 0  \mbox{.}
$$

We have the following admissible triple of $C_1$-spaces
$$
0 \arrow{e} E_1 \arrow{e,t}{\alpha' \alpha} H'
\arrow{e,t}{\alpha'_{\beta}} E_3 \mathop{\amalg}\limits_{E_2} H'
\arrow{e} 0  \mbox{.}
$$

\begin{prop}. \label{pr17}
We have the following formulas.
\begin{enumerate}
\item  For any $f \in \D(H')$
\begin{equation} \label{e'1}
 (\alpha' \alpha)^* (f) = \alpha^* (\alpha')^* (f) \mbox{.}
\end{equation}
\item For any $G \in \D'(E_1)$
\begin{equation}
\label{e'2} (\alpha' \alpha)_* (G) = (\alpha')_* \alpha_* (G)
\mbox{.}
\end{equation}
\item For any $f \in \E(H')$
\begin{equation}
\label{e'3} (\alpha' \alpha)^* (f) = \alpha^* (\alpha')^* (f)
\mbox{.}
\end{equation}
\item For any $G \in \tilde{\E}'(E_1)$
\begin{equation}
\label{e'4}  (\alpha' \alpha)_* (G) = (\alpha')_* \alpha_* (G)
\mbox{.}
\end{equation}

\vspace{0.3cm}

\noindent We suppose further that  $E_3$ and $L'$ are discrete
$C_1$-spaces. Then $E_3 \mathop{\amalg}\limits_{E_2} H'$ is a
discrete  $C_1$-space, and the following formulas are satisfied:

\item For any
$f \in \D(E_1)$
\begin{equation} \label{pr17f5}
(\alpha' \alpha)_* (f) = (\alpha')_* \alpha_* (f) \mbox{.}
\end{equation}
\item For any $G \in \D'(H')$
\begin{equation} \label{pr17f6}
(\alpha' \alpha)^* (G) = \alpha^* (\alpha')^*(G) \mbox{.}
\end{equation}
\item For any $f \in \E(E_1)$
\begin{equation} \label{pr17f7}
(\alpha' \alpha)_* (f) = (\alpha')_* \alpha_* (f) \mbox{.}
\end{equation}
\item
For any $G \in \tilde{\E}'(H')$
\begin{equation} \label{pr17f8}
(\alpha' \alpha)^* (G) = \alpha^* (\alpha')^*(G) \mbox{.}
\end{equation}
\end{enumerate}
\end{prop}
\proof. Formula~(\ref{e'1}) follows from definitions of the
corresponding maps. Formula~(\ref{e'2}) is the conjugate formula to
 formula~(\ref{e'1}). Formula~(\ref{e'3}) follows also from definitions of the
corresponding maps. We restrict formula~(\ref{e'3}) to the subspace
$\tilde{\E}(W')$. Then formula~(\ref{e'4}) is the conjugate formula
to this last formula.

From admissible triple of $C_1$-spaces
$$
0 \arrow{e} E_3 \arrow{e,t}{\beta_{\alpha'}} E_3
\mathop{\amalg}\limits_{E_2} H' \arrow{e,t}{\beta_{\beta'}} L'
 \arrow{e} 0
$$
we obtain that if $E_3$ and $L'$ are discrete $C_1$-spaces, then
$E_3 \mathop{\amalg}\limits_{E_2} H'$ is also a discrete
$C_1$-space.

Now formulas~(\ref{pr17f5}) and (\ref{pr17f7}) follow at once from
definitions. Formula~(\ref{pr17f6}) is the conjugate formula to
formula~(\ref{pr17f5}). Formula~(\ref{pr17f8}) follows from
formula~(\ref{pr17f6}). The proposition is proved.
\vspace{0.3cm}

Let
$$
0 \arrow{e} E_1  \arrow{e,t}{\alpha}  E_2 \arrow{e,t}{\beta}   E_3
\arrow{e} 0
$$
be an admissible triple of $C_1$-spaces over the field $\df_q$,
where $E_i = (I_i, F_i, V_i)$, $ 1 \le i \le 3$. Let $\gamma : D
\lto E_3$ be an admissible monomorphism, where $D = (R, S, Y)$. That
is, there is the following admissible triple of $C_1$-spaces
$$
0 \arrow{e} D \arrow{e,t}{\gamma} E_3 \arrow{e,t}{\delta} B
 \arrow{e} 0 \mbox{.}
$$

Then we have the following commutative diagram (see
section~\ref{c_1-spaces}):
$$
\begin{diagram}
\node{E_2 \mathop{\times}\limits_{E_3} D}
\arrow{s,l}{\beta_{\gamma}} \arrow{e,t}{\gamma_{\beta}} \node{D} \arrow{s,r}{\gamma} \\
\node{E_2} \arrow{e,b}{\beta} \node{E_3}
\end{diagram}
$$
where $\gamma_{\beta}$ is an admissible epimorphism, and
$\beta_{\gamma}$ is an admissible monomorphism. Let $X'=E_2
\mathop{\times}\limits_{E_3} D = (N, Q, X)$ as a $C_1$-space.

\begin{prop}[\it Base change rules]. \label{pr18}
We have the following formulas.
\begin{enumerate}
\item For any $f \in \D(E_2)$, $\mu \in
\mu(E_1)$
\begin{equation}  \label{e''1}
\gamma^* \beta_{*} (f \otimes \mu) = (\gamma_{\beta})_*
(\beta_{\gamma}^* (f) \otimes \mu)  \mbox{.}
\end{equation}
\item For any $G \in \D'(D)$, $\mu \in
\mu(E_1)$
\begin{equation} \label{e''2}
\beta^* (\gamma_*(G) \otimes \mu)=(\beta_{\gamma})_*
\gamma_{\beta}^* (G \otimes \mu)  \mbox{.}
\end{equation}
\item For any $f' \in \E(E_3)$
\begin{equation} \label{e''3}
\gamma_{\beta}^* \gamma^* (f') = \beta_{\gamma}^* \beta^* (f')
\mbox{.}
\end{equation}
\item For any $G' \in \tilde{\E}'(X')$
\begin{equation} \label{e''4}
\beta_* (\beta_{\gamma})_* (G') = \gamma_* (\gamma_{\beta})_* (G')
\mbox{.}
\end{equation}
\item
For any $f' \in \D(E_3)$ and compact  $C_1$-space $E_1$
\begin{equation} \label{pr18f4a}
\gamma_{\beta}^* \gamma^* (f') = \beta_{\gamma}^* \beta^* (f')
\mbox{.}
\end{equation}
\item
For any $G' \in \D'(X')$ and compact  $C_1$-space $E_1$
\begin{equation} \label{pr18f4b}
\beta_* (\beta_{\gamma})_* (G') = \gamma_* (\gamma_{\beta})_* (G')
\mbox{.}
\end{equation}
\item For any $f \in \D(X')$, $\mu \in \mu(E_1)$ and discrete  $C_1$-space $B$
\begin{equation} \label{pr18f4c}
\beta_* ((\beta_{\gamma})_* (f) \otimes \mu) = \gamma_*
(\gamma_{\beta})_* (f \otimes \mu) \mbox{.}
\end{equation}
\item  For any   $G \in \D'(E_3)$, $\mu \in \mu(E_1)$ and discrete  $C_1$-space $B$
\begin{equation} \label{pr18f4d}
\gamma_{\beta}^* (\gamma^* (G) \otimes \mu ) = \beta_{\gamma}^*
\beta^* (G \otimes \mu) \mbox{.}
\end{equation}

\vspace{0.3cm}

\noindent We suppose further that
  $E_1$ is a compact
$C_1$-space,  and $B$ is a discrete  $C_1$-space. Then the following
formulas are satisfied.

\item For any  $f \in \D(D)$
\begin{equation} \label{pr18f5}
\beta^* \gamma_*(f) =(\beta_{\gamma})_* \gamma_{\beta}^* (f)
\mbox{.}
\end{equation}
\item For any  $G \in \D'(E_2)$
\begin{equation}  \label{pr18f6}
\gamma^* \beta_{*} (G) = (\gamma_{\beta})_* \beta_{\gamma}^* (G)
\mbox{.}
\end{equation}
\item For any  $f \in \E(X')$, $\mu \in \mu(E_1)$
\begin{equation} \label{pr18f7}
\beta_* ((\beta_{\gamma})_* (f) \otimes \mu ) = \gamma_*
(\gamma_{\beta})_* (f \otimes \mu) \mbox{.}
\end{equation}
\item For any  $G \in \tilde{\E}'(E_3)$, $\mu \in \mu(E_1)$
\begin{equation} \label{pr18f8}
\gamma_{\beta}^* (\gamma^* (G) \otimes \mu) = \beta_{\gamma}^*
\beta^* (G \otimes \mu) \mbox{.}
\end{equation}
\end{enumerate}
\end{prop}
\proof. Formula~(\ref{e''1}) follows from the corresponding
definitions. Formula~(\ref{e''2}) is the conjugate formula to
formula~(\ref{e''1}). Formula~(\ref{e''3}) follows also from the
corresponding definitions. To obtain formula~(\ref{e''4}), we
restrict formula~(\ref{e''3}) to the subspace $\tilde{\E}(E_3)$ and
then take the conjugate formula to this formula.
Formulas~(\ref{pr18f4a}) and~(\ref{pr18f4c}) follow at once from the
corresponding definitions. Formula~(\ref{pr18f4b}) is the conjugate
formula to formula~(\ref{pr18f4a}), formula~(\ref{pr18f4d}) is the
conjugate formula  to formula~(\ref{pr18f4c}).
Formulas~(\ref{pr18f5}) and~(\ref{pr18f7}) follow at once from the
corresponding definitions. Formula~(\ref{pr18f6}) is the conjugate
formula to formula~(\ref{pr18f5}). Formula~(\ref{pr18f8}) follows
from formula~(\ref{pr18f4d}).
\vspace{0.3cm}

\subsection{Direct and inverse images and Fourier transform.}
\label{dif} Let
$$
0 \arrow{e} E_1  \arrow{e,t}{\alpha}  E_2 \arrow{e,t}{\beta}   E_3
\arrow{e} 0
$$
be an admissible triple of complete $C_1$-spaces over the field
$\df_q$, where \linebreak $E_i = (I_i, F_i, V_i)$, $ 1 \le i \le 3$.
Then we have the following admissible triple of complete
$C_1$-spaces:
$$
0 \arrow{e} \check{E_3} \arrow{e,t}{\check{\beta}} \check{E_2}
\arrow{e,t}{\check{\alpha}}  \check{E_1} \arrow{e} 0 \mbox{.}
$$

\begin{prop}
We have the following commutative diagrams.\footnote{Here and in the
sequel, the expressions of type
  $\beta_* \otimes \mu(E_3)$ are shortenings of more correct
  $\beta_* \otimes {\rm Id}_{\mu(E_3)}$}:

\begin{equation} \label{eq1}
\begin{diagram}
\node{\D(E_2) \otimes_{\cc} \mu(E_2)} \arrow{e,t}{\beta_* \otimes
\mu(E_3)}
\arrow{s,l}{\F} \node{\D(E_3) \otimes_{\cc} \mu(E_3)} \arrow{s,r}{\F} \\
\node{\D(\check{E_2})} \arrow{e,b}{\check{\beta}^*}
\node{\D(\check{E_3})}
\end{diagram}
\end{equation}

\begin{equation}  \label{eq2}
\begin{diagram}
\node{\D(E_2)}  \arrow{e,t}{\alpha^*}
\arrow{s,l}{\mu(\check{E_2})\otimes \F} \node{\D(E_1)}
\arrow{s,r}{\F \otimes \mu(\check{E_1})}
\\
\node{\D(\check{E_2})\otimes_{\cc} \mu(\check{E_2})}
\arrow{e,b}{\check{\alpha}_* \otimes \mu(\check{E_1})}
\node{\D(\check{E_1}) \otimes_{\cc} \mu(\check{E_1})}
\end{diagram}
\end{equation}

\begin{equation}  \label{eq3}
\begin{diagram}
\node{\D'(E_3) \otimes_{\cc} \mu(E_3)^*}  \arrow{e,t}{\beta^*
\otimes \mu(E_2)^*} \arrow{s,l}{\F} \node{\D'(E_2) \otimes_{\cc}
\mu(E_2)^*} \arrow{s,r}{\F}
\\
\node{\D'(\check{E_3})} \arrow{e,b}{\check{\beta}_* }
\node{\D'(\check{E_2})}
\end{diagram}
\end{equation}

\begin{equation}  \label{eq4}
\begin{diagram}
\node{\D'(E_1)}  \arrow{e,t}{\alpha_*} \arrow{s,l}{\mu(E_1) \otimes
\F } \node{\D'(E_2)}  \arrow{s,r}{\F \otimes \mu(E_2)}
\\
\node{\D'(\check{E_1}) \otimes_{\cc} \mu(E_1)}
\arrow{e,b}{\check{\alpha}^* \otimes \mu(E_1)}
\node{\D'(\check{E_2}) \otimes_{\cc} \mu(E_2)}
\end{diagram}
\end{equation}
\end{prop}
\proof. We use formulas which are analogous to formulas~(\ref{fm1})
and~(\ref{fm2}):
$$
\mu(E_1) \otimes_{\cc} \mu(E_3) = \mu(E_2)  \mbox{,}
$$
$$
\mu(\check{E_3}) \otimes_{\cc} \mu(\check{E_1}) = \mu(\check{E_2})
\mbox{,}
$$
$$
\mu(E_i) \otimes_{\cc} \mu(\check{E_i}) = \cc,  \quad  1 \le i \le 3
 \mbox{.}
$$

Now diagram~(\ref{eq1}) is equivalent to the following formula:
\begin{equation} \label{for1}
\int\limits_{V_3} \beta_*(f) \cdot \overline{\psi_w} \, d \mu_3 =
\int\limits_{V_2} f \cdot \overline{\psi_{\check{\beta}(w)}} \, d
\mu_2 \mbox{,}
\end{equation}
which has to be satisfied for any $f \in \D(E_2)$, any $w \in
\check{V_3} $, and $\mu_3 = \mu_1 \otimes \mu_2$.

Now using  $\overline{\psi_{\check{\beta}(w)}}
=\beta^*(\overline{\psi_w}) $, projection formula~(\ref{pf1}) and
Fubini formula (proposition~\ref{fubf}), we obtain
formula~(\ref{for1}).

Diagram~(\ref{eq2}) follows from diagram~(\ref{eq1}) and
proposition~\ref{propf}. Diagram~(\ref{eq3}) is the conjugate
diagram to diagram~(\ref{eq2}). Diagram~(\ref{eq4}) is the conjugate
diagram to diagram~(\ref{eq1}).
\vspace{0.3cm}

\begin{prop}
We have the following commutative diagrams.

\begin{equation}  \label{eqq1}
\begin{diagram}
\node{\E(E_3)}   \arrow{e,t}{\beta^* } \arrow{s,l}{\F}
\node{\E(E_2)}
 \arrow{s,r}{\F}
\\
\node{\tilde{\E}'(\check{E_3})} \arrow{e,b}{\check{\beta}_* }
\node{\tilde{\E}'(\check{E_2})}
\end{diagram}
\end{equation}

\begin{equation}  \label{eqq2}
\begin{diagram}
\node{\tilde{\E}'(E_1)}  \arrow{e,t}{\alpha_*} \arrow{s,l}{ \F }
\node{\tilde{\E}'(E_2)}  \arrow{s,r}{\F}
\\
\node{\E(\check{E_1}) } \arrow{e,b}{\check{\alpha}^*}
\node{\E(\check{E_2})}
\end{diagram}
\end{equation}

\begin{equation}  \label{eqq3}
\begin{diagram}
\node{\E(E_2)}  \arrow{e,t}{\alpha^*} \arrow{s,l}{\F} \node{\E(E_1)}
\arrow{s,r}{\F}
\\
\node{\tilde{\E}'(\check{E_2})} \arrow{e,b}{\check{\alpha}_* }
\node{\tilde{\E}'(\check{E_1}) }
\end{diagram}
\end{equation}

\begin{equation} \label{eqq4}
\begin{diagram}
\node{\tilde{\E}'(E_2)} \arrow{e,t}{\beta_* }
\arrow{s,l}{\F} \node{\tilde{\E}'(E_3)} \arrow{s,r}{\F} \\
\node{\E(\check{E_2})} \arrow{e,b}{\check{\beta}^*}
\node{\E(\check{E_3})}
\end{diagram}
\end{equation}
\end{prop}
\proof. Diagram~(\ref{eqq1}) follows from diagram~(\ref{eq3}),
proposition~\ref{prop15}, formula~(\ref{four'e}), and
proposition~\ref{four'e'}.

Diagram~(\ref{eqq2}) follows from diagram~(\ref{eqq1}),
proposition~\ref{four'e'}, and the following formula: $ \F \circ \F
(g) = \check{g}$ for any $g \in \E(E_i)$ (which follows from
proposition~\ref{prop9}).

We change in diagram~(\ref{eqq1}) the spaces $\E(E_i)$ on their
subspaces $\tilde{\E}(E_i)$ ($i=2, 3$) correspondingly, and the
spaces $\tilde{\E}'(\check{E_i})$ also on their subspaces
$\E'(\check{E_i})$ ($i=2, 3$) correspondingly. Then by
proposition~\ref{four'e'}, this new diagram is commutative. And the
conjugation to this diagram gives diagram~\ref{eqq3}. (We used
remark~\ref{nt3}.)

Diagram~(\ref{eqq4}) follows from diagram~(\ref{eqq3}) and the
following formula: $ \F \circ \F (g) = \check{g}$ for any $g \in
\E(E_i)$.
\vspace{0.3cm}

\subsection{The invariance property.} \label{rem}
For any $C_1$-space $E = (I, F, V)$ we defined the spaces $\mu(E)$,
$\D(E)$, $\D'(E)$, $\E(E)$, $\E'(E)$, $\tilde{\E}(E)$,
$\tilde{\E}'(E)$.

 \begin{prop}. \label{prop21}
Let  $E'$ be another  $C_1$-space such that  $E \sim E'$. Then there
are the following canonical isomorphisms: $\mu(E) = \mu(E')$,
$$
\D(E) = \D(E') \quad , \quad \D'(E) = \D'(E') \quad , \quad \E(E) =
\E(E') \mbox{,}
$$
$$
\E'(E) = \E'(E')\quad , \quad \tilde{\E}(E) = \tilde{\E}(E') \quad ,
\quad \tilde{\E}'(E) = \tilde{\E}'(E') \mbox{.}
$$
\end{prop}
\proof \ follows from the corresponding definitions.
\vspace{0.3cm}
\begin{nt}.
{\em The isomorphisms from proposition~\ref{prop21} are well-defined
with respect to isomorphisms~(\ref{fm1}) and~(\ref{fm2}), direct and
inverse images, and Fourier transform.}
\end{nt}

\section{Two-dimensional case}

\subsection{$C_2$-spaces}  \label{c_2-spaces}
We recall the definition of $C_2$-space from~\cite{Osip}.

We fix a field $k$. By definition, objects of the category $C_2$,
i.e. $\Ob(C_2)$, are filtered $k$-vector spaces $(I, F, V)$ with the
following additional structures
\begin{itemize}
\item[(i)] for any $i \le j \in I$  on the $k$-vector space $F(j) / F(i)$
it is given a structure $E_{i,j} \in \Ob(C_1)$,
\item[(ii)] for any $i \le j \le k \in I$
$$
0 \lto E_{i,j}  \lto E_{i,k}  \lto  E_{j,k} \lto 0
$$
is an admissible triple from $C_1$.
\end{itemize}

Let $E_1 = (I_1, F_1, V_1)$ and $E_2 = (I_2, F_2, V_2)$ be from
$\Ob(C_2)$. Then, by definition, $\Mor_{C_2}(E_1, E_2)$ consists of
elements $A \in \Hom_k (V_1, V_2)$ such that the following
conditions hold:
\begin{itemize}
\item[(i)] \label{i1} for any $i \in I_1$ there is an $j \in I_2$ such that  $A (F_1(i)) \subset F_2(j)$,
\item[(ii)] \label{i2}  for any $j \in I_2$ there is an $i \in I_1$ such that  $A (F_1(i)) \subset F_2(j)$,
\item[(iii)] \label{i3}  for any $i_1 \le i_2 \in I_1$ and $j_1 \le j_2 \in I_2$ such that $A (F_1(i_1)) \subset F_2(j_1)$
and $A (F_1(i_2)) \subset F_2(j_2)$ we have that the induced
$k$-linear map
$$   \bar{A} : \frac{F_1(i_2)}{F_1(i_1)} \lto \frac{F_2(j_2)}{F_2(j_1)}
$$
is an element from
$$\Mor\nolimits_{C_1}(\frac{F_1(i_2)}{F_1(i_1)}, \frac{F_2(j_2)}{F_2(j_1)})  \mbox{.}$$
\end{itemize}

Let $E_1 = (I_1, F_1, V_1)$, $E_2 = (I_2,
F_2, V_2)$ and $E_3=(I_3, F_3, V_3)$ be from
$\Ob(C_2)$. Then we say that
\begin{equation} \label{adt'}
0 \arrow{e} E_1 \arrow{e,t}{\alpha} E_2 \arrow{e,t}{\beta} E_3 \arrow{e}
0
\end{equation}
is {\em an admissible triple} from $C_2$ when
 the following conditions are satisfied:
\begin{itemize}
\item[(i)]
$$
0 \arrow{e} V_1 \arrow{e,t}{\alpha} V_2 \arrow{e,t}{\beta} V_3
\arrow{e} 0
$$
is an exact triple of $k$-vector spaces,
\item[(ii)] \label{itaa'}
the filtration $(I_1, F_1, V_1)$ dominates the filtration $(I_2,
F'_1, V_1)$, where \linebreak $F'_1 (i) = F_2(i) \cap V_1$  for any
$i \in I_2$,
\item[(iii)] \label{itbb'}
the filtration $(I_3, F_3, V_3)$ dominates the filtration $(I_2,
F'_3, V_3)$, where \linebreak $F'_3(i) =  F_2(i) / F_2(i) \cap V_1$
for any $i \in I_2$,
\item[(iv)] for any $i \le j \in I_2$
\begin{equation} \label{trojkaa1}
0 \lto \frac{F'_1(j)}{F'_1(i)} \lto \frac{F_2(j)}{F_2(i)}  \lto
\frac{F'_3(j)}{F'_3(i)} \lto 0
\end{equation}
is an admissible triple from $C_1$. (By definition of $\Ob(C_2)$, on
every vector space from triple~(\ref{trojkaa1}) it is given the
structure of $\Ob(C_1)$).
\end{itemize}

We say that $\alpha$ from an admissible triple~(\ref{adt'}) is {\em
an admissible monomorphism}, and $\beta$ from an admissible
triple~(\ref{adt'}) is {\em an admissible epimorphism.}

\begin{nt} \label{dop2} \em
Let $E_2= (I_2, F_2, V_2)$ be a $C_2$-space over the field $k$. Let
\begin{equation} \label{z2}
0 \arrow{e} V_1  \arrow{e,t}{\alpha}  V_2 \arrow{e,t}{\beta}   V_3
\arrow{e} 0
\end{equation}
be an exact triple of  $k$-vector spaces. We define on spaces $V_1$
and $V_3$ the structures of filtered spaces  $(I_i, F_i, V_i)$, $i
=1, 3$, where $I_1=I_3=I_2$, and $F_1(i) \eqdef F_2(i) \cap V_1$,
$F_3(i) \eqdef \beta(F_2(i))$ for $i \in I_2$.

By definition, for any $i \le j \in I_2$ on the space
$F_2(j)/F_2(i)$ it is given the structure of $C_1$-space, and we
have the following exact triple of $k$-vector spaces:
\begin{equation} \label{asd}
0 \lto F_1(j)/F_1(i) \lto F_2(j)/F_2(i) \lto F_3(j)/F_3(i) \lto 0
\mbox{.}
\end{equation}
Therefore, as in the case of remark~\ref{dop}, we consider the
induced filtration on the   $k$-space $F_1(j)/ F_1(i)$ and the
factor-filtration on the  $k$-space $F_3(j)/F_3(i)$.

The filtrations, so defined, gives always the structure of
$C_2$-space $E_1= (I_1, F_1,V_1)$. At the same time, the
factor-filtrations, which are defined above, will give the structure
of
 $C_2$-space $E_3= (I_3, F_3, V_3)$
iff the following conditions hold:
$$ \label{zamkn2}
 \quad \mbox{for any} \quad x \in V_2 \setminus V_1 \quad
\mbox{there exist} \quad l \in I_2 \mbox{ such that} \quad (x +
F_2(l)) \cap V_1 = \emptyset;
\leqno(**)
$$
$$
 \quad \mbox{for any} \quad i \le j \in I_2 \quad \mbox{for the exact triple~(\ref{asd}) } \leqno(***)
$$
$$  \mbox{condition~$(*)$ from
remark~\ref{dop} is satisfied} \mbox{.}
$$

Moreover, under  conditions~$(**)$, $(***)$,  exact
triple~(\ref{z2}) gives the admissible triple
$$0 \arrow{e} E_1
\arrow{e,t}{\alpha}  E_2 \arrow{e,t}{\beta}   E_3 \arrow{e} 0
\mbox{,}$$ i.e. $\alpha$ will be an admissible monomorphism, $\beta$
will be an admissible epimorphism.
\end{nt}

We say that a $C_2$-space $E_1 = (I_1, F_1, V)$ {\em dominates}
another $C_2$-space $E_2 = (I_2, F_2, V)$ if the following
conditions are satisfied:
\begin{itemize}
\item[(i)] the filtration $(I_1, F_1, V)$ dominates the filtration
$(I_2, F_2, V)$,
\item[(ii)] for any $i \le j \in I_2$  the $C_1$-space ${E_1}_{i,j}=F_2(j)/F_2(i)$
with filtration induced by $C_2$-structure on $E_1$ dominates the
$C_1$-space ${E_2}_{i,j}= F_2(j)/F_2(i)$ with filtration induced by
$C_2$-structure on $E_2$.
\end{itemize}

We consider the group $\Aut_{C_2}(E) \eqdef \Mor_{C_2}(E, E)^*$
(i.e., invertible elements in $k$-algebra $\Mor_{C_2}(E, E)$) for
any $C_2$-space $E$. By~\cite[prop.2.1]{Osip}, if a $C_2$-space
$E_1$ dominates a $C_1$-space $E_2$, then canonically
$\Aut_{C_2}(E_1) =\Aut_{C_2}(E_2)$.

Let $E = (I, F, V)$ be a $C_2$-space. We
define a $C_2$-space $\check{E} = (I^0, F^0,
\check{V})$ by the following way. The space
$\check{V} \subset V^*$ is defined as
\begin{equation}
\check {V} \eqdef \mathop{\Lim_{\rightarrow}}_{j \in I}
\mathop{\Lim_{\leftarrow}}_{i \ge j} \check{E}_{j,i} \mbox{,}
\end{equation}
where the $C_1$-space $\check{E}_{j,i}$ is constructed from the
$C_1$-space $E_{j,i}= F(i)/F(j)$ (see section~\ref{c_1-spaces}). The
set $I^0$ is a partially ordered set, which has the same set as $I$,
but with the inverse order than $I$. For $j \in I^0$ the subspace
\begin{equation}
F^0(j) \eqdef \mathop{\Lim_{\leftarrow}}_{i \le j \in I^0}
\check{E}_{j,i} \quad \subset \quad \check{V} \mbox{.}
\end{equation}

If $E_1, E_2 \in \Ob(C_2)$ and $\theta \in \Mor_{C_2}(E_1,E_2)$,
then there is canonically $\check{\theta} \in \Mor_{C_2}
(\check{E_2}, \check{E_1})$. If
$$
0 \arrow{e} E_1  \arrow{e,t}{\alpha}  E_2 \arrow{e,t}{\beta}   E_3
\arrow{e} 0
$$
is an admissible triple of $C_2$-spaces, then there is canonically
the following admissible triple of $C_2$-spaces:
$$
0 \arrow{e} \check{E_3} \arrow{e,t}{\check{\beta}} \check{E_2}
\arrow{e,t}{\check{\alpha}}  \check{E_1} \arrow{e} 0 \mbox{.}
$$

\begin{lemma} \label{lbc2}
Let
$$
0 \arrow{e} E_1  \arrow{e,t}{\alpha}  E_2 \arrow{e,t}{\beta}   E_3
\arrow{e} 0
$$
be an admissible triple of $C_2$-spaces.
\begin{enumerate}
\item Let $D \in \Ob(C_2)$, and $\gamma \in
\Mor_{C_2} (D, E_3)$. Then there is the following admissible triple
of $C_2$-spaces
\begin{equation} \label{adm12}
0 \arrow{e} E_1  \arrow{e,t}{\gamma_{\alpha}} E_2
\mathop{\times}_{E_3} D \arrow{e,t}{\gamma_{\beta}} D \arrow{e} 0
\end{equation}
and $\beta_{\gamma} \in \Mor_{C_2}(E_2 \mathop{\times}\limits_{E_3}
D, E_2)$ such that the following diagram is commutative:
\begin{equation} \label{adm22}
\begin{diagram}
\node{0} \arrow{e} \node{E_1}  \arrow{e,t}{\gamma_{\alpha}}
\arrow{s,=} \node{E_2 \mathop{\times}_{E_3} D}
\arrow{e,t}{\gamma_{\beta}}
\arrow{s,l}{\beta_{\gamma}} \node{D} \arrow{e} \arrow{s,r}{\gamma} \node{0} \\
\node{0} \arrow{e} \node{E_1}  \arrow{e,t}{\alpha}  \node{E_2}
\arrow{e,t}{\beta} \node{E_3} \arrow{e} \node{0}
\end{diagram}
\end{equation}
\item Let $A \in \Ob(C_2)$, and $\theta \in
\Mor_{C_2}(E_1, A)$.
 Then there is the following admissible triple of
$C_2$-spaces
\begin{equation}  \label{adm'12}
0 \arrow{e} A  \arrow{e,t}{\theta_{\alpha}} A  \mathop{\amalg}_{E_1}
E_2 \arrow{e,t}{\theta_{\beta}} E_3 \arrow{e} 0
\end{equation}
and $\alpha_{\theta} \in \Mor_{C_2}(E_2, A
\mathop{\amalg}\limits_{E_1} E_2)$ such that the following diagram
is commutative:
\begin{equation}  \label{adm'22}
\begin{diagram}
\node{0} \arrow{e} \node{E_1}  \arrow{e,t}{\alpha}
\arrow{s,l}{\theta} \node{E_2} \arrow{e,t}{\beta}
\arrow{s,r}{\alpha_{\theta}} \node{E_3} \arrow{e} \arrow{s,=}  \node{0} \\
\node{0} \arrow{e} \node{A}  \arrow{e,t}{\theta_{\alpha}} \node{A
\mathop{\amalg}\limits_{E_1} E_2} \arrow{e,t}{\theta_{\beta}}
\node{E_3} \arrow{e} \node{0}
\end{diagram}
\end{equation}
\end{enumerate}
\end{lemma}
\proof~follows by repeating of reasonings of the proof of
lemma~\ref{lbc} when in this proof the fibered product (or the sum)
of factors of  filtration is understood as the fibered product (or
the sum) of $C_1$-spaces, constructed in lemma~\ref{lbc}.
\vspace{0.3cm}

\begin{defin}
We say that a $C_2$-space $(I_1, F_1, V_1)$ is a  $cC_2$-space iff
there is an $i \in I_1$ such that $ F_1(i) = V_1$.

We say that a $C_2$-space $(I_2, F_2, V_2)$ is an $dC_2$-space iff
there is an $j \in I_2$ such that $F_2(j) = \{0\}$.
\end{defin}

It follows from definition that if $E_1$ is a  $cC_2$-space, then
$\check{E_1}$ is an $dC_2$-space. If $E_2$ is an $dC_2$-space, then
$\check{E_2}$ is a  $cC_2$-space.

\begin{defin}
We say that a $C_2$-space $(I_1, F_1, V_1)$ is a $cfC_2$-space iff
for any $i_1 \ge j_1 \in I_1$ the $C_1$-space $F(i_1)/F(j_1)$ is a
compact $C_1$-space.

We say that a $C_2$-space $(I_2, F_2, V_2)$ is a  $dfC_2$-space iff
for any $i_2 \ge j_2 \in I_2$ the $C_1$-space $F(i_2)/F(j_2)$ is a
discrete $C_1$-space.
\end{defin}

It follows from definition that if $E_1$ is a $cfC_2$-space, then
$\check{E_1}$ is a  $dfC_2$-space. If $E_2$ is a  $dfC_2$-space,
then $\check{E_2}$ is a  $cfC_2$-space.

\begin{defin}
We say that a $C_2$-space $(I,F,V)$ is a complete $C_2$-space if the
following conditions are satisfied:
\begin{enumerate}
\item
For any $i \ge j \in I$ the $C_1$-space $E_{j,i}= F(i)/F(j)$ is a
complete $C_1$-space.
\item
$$
V = \mathop{\Lim_{\to}}_{i \in I}  \mathop{\Lim_{\gets}}_{j \le i}
F(i)/F(j)  \mbox{.}
$$
\end{enumerate}
\end{defin}

It follows from definitions that $E$ is a complete $C_2$-space iff
$\check{\check{E}} = E$.

\subsection{Virtual measures.} \label{vm}
Let $E = (I, F, V)$ be a $C_2$-space over the field $\df_q$.

For any $i,j \in I$ we define a $1$-dimensional $\cc$-vector space
of virtual measures
\begin{equation} \label{virm}
\mu(F(i) \mid F(j)) \eqdef \mathop{\mathop{\Lim_{\lto}}_{l \in
I}}_{l \le i, l \le j} \Hom\nolimits_{\cc} (\mu(F(i)/ F(l)) \, , \,
\mu (F(j)/ F(l))) \mbox{,}
\end{equation}
where  to take  the inductive limit we need the identities which
follow from formula~(\ref{fm2}): for $l' \le l \in I$
$$
\mu(F(i)/ F(l')) = \mu(F(i)/ F(l)) \otimes_{\cc} \mu(F(l)/ F(l'))
$$
$$
\mu(F(j)/ F(l')) = \mu(F(j)/ F(l)) \otimes_{\cc} \mu(F(l)/ F(l'))
\mbox{.}
$$
And
$$
\textstyle f \in \Hom\nolimits_{\cc} (\mu(F(i)/ F(l)), \mu (F(j)/
F(l))) \mapsto f' \in \Hom\nolimits_{\cc} (\mu(F(i)/ F(l')), \mu
(F(j)/ F(l')))  \mbox{,}$$
 where $f'(a \otimes c) \eqdef f(a) \otimes c$, $a$ is any
from $\mu(F(i)/ F(l))$, $c$ is any from $\mu(F(l)/ F(l'))$. This map
is an isomorphism.

\begin{prop} \label{pci}
For any $i,j,l \in I$ there is a canonical isomorphism
$$
\gamma \quad : \quad \mu(F(i) \mid F(j)) \otimes_{\cc} \mu(F(j) \mid
F(l))  \arrow{e} \mu(F(i) \mid F(l))
$$
such that the following diagram of associativity is commutative for
any $i,j,l,n \in I$:
$$
\begin{diagram}
 \node{\scriptstyle \mu(F(i) \mid F(j)) \otimes_{\cc}
\mu(F(j) \mid F(l)) \otimes_{\cc} \mu(F(l) \mid F(n)) }  \arrow{e}
\arrow{s}
\node{\scriptstyle \mu(F(i) \mid F(l)) \otimes_{\cc} \mu(F(l) \mid F(n))} \arrow{s} \\
\node{\scriptstyle \mu(F(i) \mid F(j)) \otimes_{\cc} \mu(F(j) \mid
F(n))} \arrow{e} \node{\scriptstyle \mu(F(i) \mid F(n))}
\end{diagram}
$$
\end{prop}
\proof. We have a canonical map:
\begin{equation} \label{eqt}
\begin{diagram}
\node{ \Hom\nolimits_{\cc} (\mu(F(i)/ F(l')), \mu (F(j)/ F(l')))
\otimes_{\cc} \Hom\nolimits_{\cc} (\mu(F(j)/ F(l')), \mu
(F(l)/F(l')))} \arrow{s} \\
\node{ \Hom\nolimits_{\cc} (\mu(F(i)/ F(l')), \mu (F(l)/ F(l')))}
\end{diagram}
\end{equation}
which satisfies the associativity diagram. And this map commutes
with the inductive limit from the definition of $\mu (\cdot \mid
\cdot)$. We obtain the map $\gamma$ after the taking the inductive
limit in~(\ref{eqt}).
\vspace{0.3cm}

\begin{nt} \label{nci} { \em
We have the following canonical isomorphisms. For any $i, j \in I$
$$
\mu(F(i) \mid F(i)) =\cc \quad , \quad \mu(F(i) \mid F(j)) =
\mu(F(j) \mid F(i))^*  \mbox{.}
$$
Let $i \le j \in I$. Then
$$
\mu(F(i) \mid F(j)) = \mu (F(j)/F(i)) \quad , \quad \mu(F(j) \mid
F(i))= \mu (F(j)/F(i))^*  \mbox{.}
$$ }
\end{nt}

\subsection{Basic spaces} \label{bs}
Let $E = (I, F, V)$ be a $C_2$-space over
the field $\df_q$.  Let $i \ge j \ge l \ge n
\in I$. Then $F(i) \supset F(j) \supset F(l)
\supset F(n)$, and from the definition of a
$C_2$-space we have
$$
F(j)/ F(n) \sim (F(i)/ F(n)) \mathop{\times}\limits_{F(i)/F(l)}
(F(j)/F(l)) $$
$$F(i)/ F(l)  \sim (F(i)/F(n))
\mathop{\amalg}\limits_{F(j)/F(n)} (F(j)/F(l))
$$
as $C_1$-spaces.

For any $i \ge j \ge l \in I$ we have the following admissible
triple of $C_1$-spaces:
$$
0 \arrow{e} F(j)/ F(l)  \arrow{e,t}{\alpha_{lji}} F(i)/F(l)
\arrow{e,t}{\beta_{lji}} F(i)/F(j) \arrow{e} 0 \mbox{.}
$$

We fix some $o \in I$. Then from section~\ref{bcr} we have that the
following spaces over the field $\cc$ are well defined.

\begin{equation} \label{2d}
\begin{array}{c}
\D_{F(o)}(E) \eqdef \mathop{\Lim\limits_{\leftarrow}}\limits_i
\mathop{\Lim\limits_{\leftarrow}}\limits_l \D(F(i) / F(l))
\otimes_{\cc} \mu(F(l) \mid F(o)) =  \qquad \qquad \qquad \qquad \qquad \\
\qquad \qquad \qquad \qquad \qquad
=
\mathop{\Lim\limits_{\leftarrow}}\limits_l
\mathop{\Lim\limits_{\leftarrow}}\limits_i \D(F(i) / F(l))
\otimes_{\cc} \mu(F(l) \mid F(o))
\end{array}
\end{equation}
with respect to the maps $(\beta_{lji})_* \otimes \mu(F(j) \mid
F(o))$:
$$
 \D(F(i)/F(l)) \otimes_{\cc} \mu(F(l) \mid F(0)) \arrow{e} \D(F(i)/
F(j)) \otimes_{\cc} \mu(F(j) \mid F(o)) \mbox{,}
$$
and the maps $\alpha^*_{lji} \otimes \mu(F(l) \mid F(o))$:
$$ \D(F(i)/ F(l)) \otimes_{\cc} \mu(F(l) \mid F(o)) \arrow{e}
\D(F(j)/ F(l)) \otimes_{\cc} \mu(F(l) \mid F(o))
  \mbox{.}$$

\begin{equation}  \label{2d'}
\begin{array}{c}
\D'_{F(o)}(E) \eqdef \mathop{\Lim\limits_{\to}}\limits_i
\mathop{\Lim\limits_{\to}}\limits_l \D'(F(i) / F(l))
\otimes_{\cc} \mu(F(o) \mid F(l)) =  \qquad \qquad \qquad \qquad \qquad \\
\qquad \qquad \qquad \qquad \qquad
=
\mathop{\Lim\limits_{\to}}\limits_l
\mathop{\Lim\limits_{\to}}\limits_i \D'(F(i) / F(l)) \otimes_{\cc}
\mu(F(o) \mid F(l))
\end{array}
\end{equation}
with respect to the maps $\beta^*_{lji} \otimes \mu(F(o) \mid
F(l))$:
$$
\D'(F(i)/F(j)) \otimes_{\cc} \mu(F(o) \mid F(j)) \arrow{e}
\D'(F(i)/F(l)) \otimes_{\cc} \mu(F(o) \mid F(l)) \mbox{,}
$$
and the maps $(\alpha_{lji})_* \otimes \mu(F(o) \mid F(l))$:
$$ \D'(F(j)/ F(l)) \otimes_{\cc} \mu(F(o) \mid F(l)) \arrow{e}
\D'(F(i)/ F(l)) \otimes_{\cc} \mu(F(o) \mid F(l))
  \mbox{.}
$$

\begin{equation} \label{2e}
\E(E) \eqdef \mathop{\Lim_{\gets}}_i \mathop{\Lim_{\to}}_l \E(F(i) /
F(l))
\end{equation}
with respect to the maps $\alpha_{lji}^*$  and $\beta_{lji}^*$.

\begin{equation} \label{2te}
\tilde{\E}(E) \eqdef \mathop{\Lim_{\to}}_l \mathop{\Lim_{\gets}}_i
\tilde{\E}(F(i) / F(l))
\end{equation}
with respect to the maps  $\beta_{lji}^*$  and  $\alpha_{lji}^*$.

\begin{equation} \label{2e'}
\E' (E) \eqdef \mathop{\Lim_{\to}}_i \mathop{\Lim_{\gets}}_l
\E'(F(i) / F(l))
\end{equation}
with respect to the maps $(\alpha_{lji})_*$ and $(\beta_{lji})_*$.

\begin{equation} \label{2te'}
\tilde{\E}' (E) \eqdef \mathop{\Lim_{\gets}}_l \mathop{\Lim_{\to}}_i
\tilde{\E}'(F(i) / F(l))
\end{equation}
 with respect to the maps $(\beta_{lji})_*$ and
$(\alpha_{lji})_*$.

From definitions of these spaces we have the following canonical
isomorphism for any $o_1 \in I$:
$$
\D_{F(o)}(E)  \otimes_{\cc} \mu(F(o) \mid F(o_1))  \arrow{e}
\D_{F(o_1)}(E)
$$
such that the following diagram is commutative for any $o_2 \in I$:
$$
\begin{diagram}
 \node{\scriptstyle   \D_{F(o)}(E)  \otimes_{\cc}
\mu(F(o) \mid F(o_1)) \otimes_{\cc} \mu(F(o_1) \mid F(o_2))}
\arrow{e} \arrow{s} \node{\scriptstyle \D_{F(o_1)}(E) \otimes_{\cc}
\mu(F(o_1) \mid F(o_2))}
\arrow{s} \\
\node{\scriptstyle   \D_{F(o)}(E) \otimes_{\cc} \mu(F(o) \mid
F(o_2))} \arrow{e} \node{\scriptstyle   \D_{F(o_2)} (E)}
\end{diagram}
$$

Dually, we have the following canonical isomorphism for any $o_1 \in
I$:
$$
\mu(F(o_1) \mid F(o)) \otimes_{\cc}  \D'_{F(o)}(E) \arrow{e}
\D'_{F(o_1)}(E)
$$
such that the following diagram is commutative for any $o_2 \in I$:
$$
\begin{diagram}
 \node{\scriptstyle  \mu(F(o_2) \mid F(o_1)) \otimes_{\cc} \mu(F(o_1) \mid F(o)) \otimes_{\cc}
 \D'_{F(o)}(E) } \arrow{e} \arrow{s} \node{\scriptstyle \mu(F(o_2) \mid F(o_1)) \otimes_{\cc} \D'_{F(o_1)}(E)}
\arrow{s} \\
\node{\scriptstyle  \mu(F(o_2) \mid F(o)) \otimes_{\cc}
\D'_{F(o)}(E)} \arrow{e} \node{\scriptstyle   \D'_{F(o_2)} (E)}
\end{diagram}
$$

Since for any $i \ge l \in I$ there are canonical embeddings
$$
\tilde{\E}(F(i)/F(l))  \hookrightarrow \E(F(i)/F(l))
$$
$$
\E'(F(i)/F(l)) \hookrightarrow \tilde{\E}'(F(i)/F(l))  \mbox{,}
$$
there are the following canonical embeddings:
$$
\tilde{\E}(E)  \hookrightarrow \E(E)
$$
$$
\E'(E) \hookrightarrow \tilde{\E}'(E)  \mbox{.}
$$

Moreover, we have an embedding $\E(E) \lto \ff(V)$ such that
$\tilde{\E}(E)$ and $\E(E)$ are $\cc$-subalgebras of $\cc$-algebra
$\ff(V)$.

\begin{prop} \label{prm}
The following statements are satisfied.
\begin{enumerate}
\item
The $\cc$-vector space $\D_{F(o)}(E)$ is a module over the
$\cc$-algebra $\E(E)$.
\item
The $\cc$-vector space $\D'_{F(o)}(E)$ is a module over the
$\cc$-algebra $\E(E)$.
\end{enumerate}
\end{prop}
\proof. We construct the map:
\begin{equation} \label{exm}
\E(E)  \otimes_{\cc} \D_{F(o)}(E)  \arrow{e} \D_{F(o)}(E) \mbox{.}
\end{equation} For any $i_1 \ge j_1 \in I$ we have the following
multiplication map:
\begin{equation} \label{exm1}
 \scriptstyle \E(F(i_1)/F(j_1)) \otimes_{\cc} \D(F(i_1)/
F(j_1)) \otimes_{\cc} \mu(F(j_1) \mid F(o)) \arrow{e,t}{\lambda}
\D(F(i_1)/ F(j_1)) \otimes_{\cc} \mu(F(j_1) \mid F(o)) \mbox{.}
\end{equation}
Now from projection formula~(\ref{pf1}) we have the following
formula for any $i \ge j \ge l \in I$, for any $g \in \E(F(i)/
F(j))$, $f \in \D(F(i)/F(l))$,  $ \mu \in \mu(F(l) \mid F(o))$:
\begin{equation} \label{pfor}
((\beta_{lji})_* \otimes \mu (F(j)/F(o))) (\lambda(\beta^*_{lji}(g)
\otimes f \otimes \mu)) = \lambda( g \otimes ((\beta_{lji})_*
\otimes \mu (F(j)/F(o))) (f \otimes \mu) )  \mbox{.}
\end{equation}

Using formulas~(\ref{2e}) and~(\ref{2d}) we construct
map~(\ref{exm}) from  maps~(\ref{exm1}). Due to
formulas~(\ref{pfor}) it is well defined. We proved the first
statement of the proposition.

For any $i_1 \ge j_1 \in I$, any $f \in \E(F(i_1)/F(j_1))$ we
defined in remark~\ref{rm4} the action of $f$ on
$\D'(F(i_1)/F(j_1))$ as conjugate action to the action of $f$ on
$\D(F(i_1)/F(j_1))$. Therefore applying the conjugate formulas to
the formulas above, from formulas~(\ref{2e}) and~(\ref{2d'}) we
construct the conjugate action of $\E(E)$ on $\D'_{F(o)}(E)$.
\vspace{0.3cm}

\begin{nt} \label{nntt} { \em
Let $\{B_l\}$, $l \in J$ be a projective system of vector spaces
over a field $k$ with the surjective transition maps $\phi_{l_1,l_2}
: B_{l_1} \to B_{l_2}$ for any $l_1 \ge l_2 \in J$. Let $J^0$ be a
partially ordered set with the same set as $J$ but with inverse
order then $J$. Let $\{A_l\}$, $l \in J^0$ be an inductive system of
$k$-vector spaces with the injective transition maps $\psi_{l_2,
l_1} : A_{l_2} \to A_{l_1}$. Suppose for any $l \in J$ we have a
nondegenerate $k$-linear pairing
$$
<\cdot , \cdot >_l \quad B_l \times A_l  \arrow{e}  k
$$
such that for any $l_1 \ge l_2 \in I$, $x \in B_{l_1}$, $y \in
A_{l_2}$ we have
$$
<\phi_{l_1, l_2}(x), y>_{l_2} = <x, \psi_{l_2, l_1}(y)>_{l_1}
\mbox{.}
$$
Then we have canonically nondegenerate $k$-linear pairing between
$k$-vector spaces $\mathop{\Lim\limits_{\gets}}\limits_{l \in J}
B_l$ and $\mathop{\Lim\limits_{\to}}\limits_{l \in J^0} A_l$ which
is induced by pairings $<\cdot , \cdot >_l$, $l \in J$. }
\end{nt}

For  a $C_2$-space $E = (I, F, V)$ over the field $\df_q$, for any
$i \ge j \in I$ there are nondegenerate $\cc$-linear pairings
between $\D(F(i)/F(j))$ and $\D'(F(i)/F(j))$, between
$\E(F(i)/F(j))$ and $\E'(F(i)/F(j))$, between
$\tilde{\E}(F(i)/F(j))$ and $\tilde{\E}'(F(i)/F(j))$ (see
remark~\ref{nt3}).

Therefore applying remark~\ref{nntt} twice to
formulas~(\ref{2d})-(\ref{2te'}) we obtain that there are the
following nondegenerate pairings:
$$
<\cdot, \cdot >_{\D_{F(o)}(E)} \quad : \quad \D_{F(o)}(E) \times
\D'_{F(o)}(E) \arrow{e}  \cc
$$
$$
<\cdot, \cdot >_{\E(E)} \quad : \quad \E(E) \times \E'(E) \arrow{e}
\cc
$$
$$
<\cdot, \cdot >_{\tilde{\E}(E)} \quad : \quad \tilde{\E}'(E) \times
\tilde{\E}(E) \arrow{e} \cc \mbox{.}
$$

Now from the  construction of multiplications from
proposition~\ref{prm} we have for any $g \in \E(E)$, $f \in
\D_{F(o)}(E)$, $G \in \D'_{F(o)}(E)$:
\begin{equation}  \label{mf1}
<f \: , \: g \cdot G >_{\D_{F(o)}(E)} = < g \cdot f \: , \:
G>_{\D_{F(o)}(E)} \mbox{.}
\end{equation}

\begin{nt}
\label{remark9}
 { \em
Let $E = (I,F,V)$ be a $cfC_2$-space over the field $\df_q$. Then
for any $i \ge j \in I$ we have a canonical element
$$1_{ij} \in \mu
(F(i)/F(j)) \quad  : \quad 1_{ij} (\delta_{(F(i)/F(j))}) \eqdef 1
\mbox{.}$$ Therefore for any $k, l \in I$ we have a canonical
element $1_{kl} \in \mu (F(k) \mid F(l))$ such that for any $k,l,n
\in I$ according to proposition~\ref{pci}:
$$
\gamma( 1_{kl} \otimes 1_{ln})  = 1_{kn} \mbox{.}
$$
Taking into account these elements $1_{lo}$ and $1_{ol}$ for $l, o
\in I$ we can do not write $1$-dimensional $\cc$-spaces $\mu(F(l)
\mid F(o) )$ and $\mu(F(o) \mid F(l) )$ in formulas~(\ref{2d})
and~(\ref{2d'}) which define $\cc$-spaces $\D_{F(o)(E)}$ and
$\D'_{F(o)}(E)$. These $\cc$-spaces {\it do not depend} on the
choice of $o \in I$.

Similarly, Let $E = (I,F, V)$ be a $dfC_2$-space over the field
$\df_q$. Then for any $i \ge j \in I$ we have a canonical element
$$\delta_{ij} \in \mu
(F(i)/F(j)) \quad  : \quad \delta_{ij} (\{0\}) \eqdef 1 \mbox{.}$$
Therefore for any $k, l \in I$ we have a canonical element
$\delta_{kl} \in \mu (F(k) \mid F(l))$ such that for any $k,l,n \in
I$ according to proposition~\ref{pci}:
$$
\gamma( \delta_{kl} \otimes \delta_{ln})  = \delta_{kn} \mbox{.}
$$
Taking into account these elements $\delta_{lo}$ and $\delta_{ol}$
for $l, o \in I$ we can do not write $1$-dimensional $\cc$-spaces
$\mu(F(l) \mid F(o) )$ and $\mu(F(o) \mid F(l) )$ in
formulas~(\ref{2d}) and~(\ref{2d'}) which define $\cc$-spaces
$\D_{F(o)(E)}$ and $\D'_{F(o)}(E)$. These $\cc$-spaces {\it do not
depend} on the choice of $o \in I$.
 }
\end{nt}

\begin{nt} { \em
From a $C_2$-space $E= (I,F,V)$ over the field $\df_q$,  we
constructed $6$ $\cc$-spaces: $\D_{F(o)}(E)$, $\D'_{F(o)}(E)$ (for
any $o \in I$), $\E(E)$, $\tilde{\E}(E)$, $\E'(E)$, $\tilde{\E}'(E)$
by formulas~(\ref{2d})-(\ref{2te'}).

We can construct 10 $\cc$-spaces, where additional $\cc$-spaces
appear if for any $i \ge l \in I$ we put $\tilde{\E}(F(i)/F(l))$
instead $\E(F(i)/F(l))$ in formula~(\ref{2e}), $\E(F(i)/F(l))$
instead $\tilde{\E}(F(i)/F(l))$ in formula~(\ref{2te}),
$\E'(F(i)/F(l))$ instead $\tilde{\E}'(F(i)/F(l))$ in
formula~(\ref{2te'}), $\tilde{\E}'(F(i)/F(l))$ instead
$\E'(F(i)/F(l))$ in formula~(\ref{2e'}).

In the sequel we restrict ourselves to consider $6$ $\cc$-spaces
which are given by formulas~(\ref{2d})-(\ref{2te'}). }
\end{nt}

\subsection{Fourier transform} \label{ft}
Let $E = (I, F, V)$ be a complete $C_2$-space over the field
$\df_q$. Then we have $\check{\check{E}} =E$. We fix some $o \in I$.

\subsubsection{}
Let $i \ge j \in I$ be any. By
section~\ref{FT} there are the following
$\cc$-linear maps:
$$
f \in \E(F(i)/F(j))  \longmapsto \check{f} \in \E(F(i)/F(j))
$$
$$
G \in \tilde{\E}'(F(i)/F(j)) \longmapsto \check{G} \in
\tilde{\E}'(F(i)/F(j))
$$
$$
f \otimes \mu \in \D(F(i)/F(j)) \otimes_{\cc} \mu(F(j) \mid F(o))
\mapsto \check{f} \otimes \mu \in \D(F(i)/F(j)) \otimes_{\cc}
\mu(F(j) \mid F(o))
$$
$$
G \otimes \mu' \in \D'(F(i)/F(j)) \otimes_{\cc} \mu(F(o) \mid F(j))
\mapsto \check{G} \otimes \mu' \in \D(F(i)/F(j)) \otimes_{\cc}
\mu(F(o) \mid F(j))  \mbox{.}
$$
They map the $\cc$-subspaces $\tilde{\E}(F(i)/F(j)) \subset
\E(F(i)/F(j))$ and $\E'(F(i)/F(j)) \subset \tilde{\E}'(F(i)/F(j))$
to itself correspondingly.

These maps commute with direct and inverse images when we change $i
,j \in I$ to $i', j' \in I $, $i' \ge i$, $j' \le j$. Therefore
using formulas~(\ref{2d})-(\ref{2te'}) we have that the following
$\cc$-linear maps are well defined:
$$
f \in \D_{F(o)}(E)   \longmapsto \check{f} \in \D_{F(o)}(E)
$$
$$
G \in \D'_{F(o)}(E)   \longmapsto \check{G} \in \D'_{F(o)}(E)
$$
$$
f \in \E(E)   \longmapsto \check{f} \in \E(E)
$$
$$
G \in \tilde{\E}'(E)   \longmapsto \check{G} \in \tilde{\E}'(E)
\mbox{.}
$$
They map the $\cc$ subspaces  $\tilde{\E}(E) \subset \E(E)$ and
$\E'(E) \subset \tilde{\E}'(E)$ to itself correspondingly.

The square of any from these maps is the indentity map. From the
corresponding formulas in 1-dimensional case we have the following
formulas:
$$
<\check{f}, G>_{\D_{F(o)}(E)} = <f, \check{G}>_{\D_{F(o)}(E)} \quad
\mbox{for any} \quad f \in \D_{F(o)}(E), \quad G \in \D'_{F(o)}(E)
\mbox{;}
$$
$$
<\check{f}, G>_{\E(E)} = <f, \check{G}>_{\E(E)} \quad \mbox{for any}
\quad f \in \E(E), \quad G \in \E'(E) \mbox{;}
$$
$$
<\check{G}, f >_{\tilde{\E}(E)} = <G, \check{f}>_{\tilde{\E}(E)}
\quad \mbox{for any} \quad G \in\tilde{\E}'(E), \quad f \in
\tilde{\E}(E)  \mbox{.}
$$

\subsubsection{}
\label{2Ft}
According to section~\ref{c_2-spaces} we have the $C_2$-space
$\check{E} = (I^0, F^0, \check{V})$. For any $i \ge j \in I$, for
the $C_1$-space $E_{j,i} = F(i)/F(j)$ we have the $C_1$-space
$\check{E}_{j,i} = F^0(j) / F^0(i)$.

For any $l,n \in I$ we have
\begin{equation} \label{fff}
\mu(F(l) \mid F(n)) = \mu (F^0(l) \mid F^0(n))  \mbox{.}
\end{equation}
Indeed, let  $k \in I$ such that $k \le l$, $k \le n$. Then
according to proposition~\ref{pci}, remark~\ref{nci} and
formula~(\ref{fm1}) we have
$$
\mu (F(l) \mid F(n)) = \mu(F(l) \mid F(k)) \otimes_{\cc} \mu(F(k)
\mid F(n)) = $$
$$ = \mu(F(l)/F(k))^* \otimes_{\cc} \mu(F(n)/F(k)) =
 \mu(F^0(k) / F^0(l)) \otimes_{\cc} \mu(F^0(k) / F^0(n))^* =$$
$$
= \mu(F^0(l) \mid F^0(k)) \otimes_{\cc} \mu(F^0(k) \mid F^0(n))
=
\mu(F^0(l) \mid F^0(n)) \mbox{.}
$$

Therefore we have the following maps for any $i \ge j \in I$:
\begin{equation}
\begin{array}{cc}
\F  \otimes  \mu(F(i) \mid F(o)) \quad :
\\
  \D(F(i)/ F(j)) \otimes_{\cc} \mu (F(j) \mid F(o)) \arrow{e}
\D(F^0(j)/ F^0(i)) \otimes_{\cc} \mu(F^0(i) \mid F^0(o))
\end{array}
\end{equation}
\begin{equation}
\begin{array}{c}
\F \: \otimes \: \mu(F(o) \mid F(i)) \quad : \\  \D'(F(i)/ F(j)) \oc
\mu (F(o) \mid F(j)) \arrow{e} \D'(F^0(j)/ F^0(i)) \otimes_{\cc}
\mu(F^0(o) \mid F^0(i))
\end{array}
\end{equation}
\begin{equation} \label{f74}
\F \quad : \quad \E(F(i)/ F(j))  \arrow{e}  \tilde{\E}'(F^0(j)/
F^0(i))
\end{equation}
\begin{equation} \label{f75}
\F \quad : \quad \tilde{\E}'(F(i)/ F(j))  \arrow{e} \E(F^0(j)/
F^0(i))
\end{equation}
\begin{equation} \label{f76}
\F \quad : \quad \tilde{\E}(F(i)/ F(j))  \arrow{e} \E'(F^0(j)/
F^0(i))
\end{equation}
\begin{equation} \label{f77}
\F \quad : \quad \E'(F(i)/ F(j))  \arrow{e} \tilde{\E}(F^0(j)/
F^0(i)) \mbox{.}
\end{equation}

We remark that the map $\F$ from formula~(\ref{f74}) being
restricted to the $\cc$-subspace $\tilde{\E}(F(i)/ F(j)) \subset
\E(F(i)/ F(j))$ coincides with the map $\F$ from
formula~(\ref{f76}). The map $\F$ from formula~(\ref{f75}) being
restricted to the $\cc$-subspace $\E'(F(i)/ F(j)) \subset
\tilde{\E}'(F(i)/ F(j)) $ coincides with the map $\F$ from
formula~(\ref{f77}).

Now we use diagrams~(\ref{eq1})-(\ref{eq4}),
(\ref{eqq1})-(\ref{eqq4}), which connect Fourier transform with
direct and inverse images. We use it when we change $i ,j \in I$ to
$i', j' \in I $, $i' \ge i$, $j' \le j$, and take the limits
according to formulas~(\ref{2d})-(\ref{2te'}). We obtain that the
following $\cc$-linear maps ({\it two-dimensional Fourier
transforms}) are well defined:
\begin{equation} \label{ae1}
\F \quad : \quad \D_{F(o)}(E) \arrow{e} \D_{F^0(o)}(\check{E})
\end{equation}
\begin{equation} \label{ae2}
\F \quad : \quad \D'_{F(o)}(E) \arrow{e} \D'_{F^0(o)}(\check{E})
\end{equation}
\begin{equation} \label{ae3}
\F \quad : \quad \E(E) \arrow{e} \tilde{\E}'(\check{E})
\end{equation}
\begin{equation} \label{ae4}
\F \quad : \quad \tilde{\E}'(E) \arrow{e} \E(\check{E})
\end{equation}
\begin{equation} \label{ae5}
\F \quad : \quad \tilde{\E}(E) \arrow{e} \E'(\check{E})
\end{equation}
\begin{equation} \label{ae6}
\F \quad : \quad \E'(E) \arrow{e} \tilde{\E}(\check{E}) \mbox{.}
\end{equation}

We remark that from the construction and the corresponding
one-dimensional property it follows that the map $\F$ from
formula~(\ref{ae3}) being restricted to the $\cc$-subspace
$\tilde{\E}(E) \subset \E(E)$ coincides with the map $\F$ from
formula~(\ref{ae5}). The map $\F$ from formula~(\ref{ae4}) being
restricted to the $\cc$-subspace $\E'(E) \subset \tilde{\E}'(E) $
coincides with the map $\F$ from formula~(\ref{ae6}).

\begin{prop}
The following statements are satisfied.
\begin{enumerate}
\item The maps $\F$ are isomorphisms of $\cc$-vector spaces from
formulas~(\ref{ae1})-(\ref{ae6}).
\item For any $f \in \D_{F(o)}(E)$ and $G \in \D'_{F(o)}(E)$
$$
\F \circ \F (f) = \check{f}  \mbox{,} \quad  \F \circ \F (G) =
\check{G} \mbox{,}
$$
for any $f \in \E(E)$ and $G \in \tilde{\E}'(E)$
$$
\F \circ \F (f) = \check{f}  \mbox{,} \quad \F \circ \F (G) =
\check{G} \mbox{;}
$$
\item For any $f \in \D_{F(o)}(E)$ and $G \in \D'_{F(o)}(E)$
$$
<\F(f), G>_{\D_{F(o)}(E)} = <f, \F(G)>_{\D_{F(o)}(E)}  \mbox{,}
$$
for any $f \in \E(E)$ and $ g \in \tilde{\E}(E)$
$$
<\F(f), g>_{\tilde{\E}(E)} = <f, \F(g)>_{\E(E)}  \mbox{,}
$$
for any $ H \in \tilde{\E}'(E)$, $ G \in \E'(E)$
$$
<\F(H), G >_{\E(E)} = <H, \F(G)>_{\tilde{\E}(E)}  \mbox{.}
$$
\end{enumerate}
\end{prop}
\proof follows from the construction of two-dimensional Fourier
transform given in this section and the properties of 1-dimensional
Fourier transform (see section~\ref{FT}) applied to functions and
distributions on $C_1$-spaces $F(i)/F(j)$ for any $i \ge j \in I$.
\vspace{0.3cm}

\subsection{Central extension and its representations} \label{ce}
In this section we construct a central extension of group of
automorphisms of an object of the category $C_2$ and consider
representations of this group on the spaces of functions adn
distributions of the given object.
\subsubsection{}
\label{ss}
 Let a $C_2$-space   $E_1= (I_1, F_1, V)$ dominates a
$C_2$-space $E_2= (I_2, F_2, V)$ (over the field $\df_q$). We fix
some $o \in I_2$. Then we obtain the following canonical
isomorphisms:
$$ \D_{F(o)}(E_1) = \D_{F(o)}(E_2) \mbox{,} \qquad \qquad
  \D'_{F(o)}(E_1) = \D'_{F(o)}(E_2)  \mbox{,}
$$
$$
\E(E_1) = \E(E_2)\mbox{,} \qquad \qquad \E'(E_1) = \E'(E_2)\mbox{,}
$$
$$
\tilde{\E}(E_1) = \tilde{\E}(E_2) \mbox{,} \qquad \qquad
\tilde{\E}'(E_1) = \tilde{\E}'(E_2)   \mbox{.}
$$
Indeed, from formulas~(\ref{2d})-(\ref{2te'}) we obtain that
inductive and projective limits used in the definitions of these
spaces will be the same on  sets of indices depending on $I_1$ or
$I_2$. Moreover, by the same reasons, the two-dimensional Fourier
transform coincide for the spaces depending on $E_1$ or $E_2$.
(Compare it with section~\ref{rem}.)

\subsubsection{} \label{ssec}
Let $E= (I, F, V)$ be a $C_2$-space over a
field $k$. We construct the maximal
$C_2$-space $\tilde{E}$ over the field $k$,
which dominates the $C_2$-space $E$. By
$G(E)$ we denote the set of all
$k$-subspaces $W \subset V$ such that there
are some $i \ge j \in I$ with the
properties:
\begin{itemize}
\item
$ F(i)  \supset W \supset F(j)$ ,
\item
for any $x \in F(i)/F(j)$, $x \notin W/F(j)$ there is $l \in I_{ij}$
such that \linebreak
 $(x+ F_{ij}(l)) \cap (W/F(j)) = \emptyset $ (this intersection is inside $F(i)/F(j)$),
where $(I_{ij}, F_{ij}, F(i)/F(j))$ is the structure of $C_1$-space
on the $k$-space $F(i)/F(j)$.
\end{itemize}

Then $G(E)$ is a partially ordered set, which is ordered by
inclusions of subspaces. We remark that for any $k , l \in G(E)$
there are $p_1 \ge p_2 \in G(E)$  such that
$p_1 \ge k \ge p_2$, $p_1 \ge l \ge p_2$. Therefore $V$ is a
filtered $k$-vector space with filtration given by $G(E)$ (see section~\ref{c_1-spaces}).

For any $W_1 , W_2 \in G(E)$ such that $W_1 \supset W_2$ we
introduce the structure of $C_1$-space on the $k$-vector space $W_1
/W_2$ in the following way. There are $i \ge j \in I$ such that
$$
F(i ) \supset W_1  \supset W_2 \supset F(j) \mbox{.}
$$
 The $k$-space $F(i)/F(j)$ is a $C_1$-space, therefore it is filtered. We
restrict the filtration from the $k$-space $F(i)/F(j)$ to the
$k$-space $W_1 / F(j)$. After that we take the factor filtration on
the $k$-space $W_1/W_2$. By $G_{i,j}(W_1/W_2)$ we denote the set of
$k$-subspaces of $W_1/W_2$ which are given by this filtration. Now
by $Gr(W_1 /W_2)$ we denote the set of all $k$-subspaces $U \subset
W_1/W_2$ such that there are some $P , Q \in G_{i,j}(W_1/W_2)$ with
the property
$$
P  \supset U \supset Q \mbox{.}
$$
Then $Gr(W_1/W_2)$ is a partially ordered set, which is ordered by
inclusions of subspaces. The filtration on $W_1/W_2$ which is given
by elements of $Gr(W_1/W_2)$ defines a well-defined $C_1$-space on
$W_1/W_2$.

Thus we constructed the well-defined $C_2$-space $\tilde{E}$, which
is filtered by the set $G(E)$. We denote this $C_2$-space by
$(\tilde{I}, \tilde{F}, V)$.

\subsubsection{} Let $E=(I, F, V)$ be a $C_2$-space over the field
$\df_q$. Then for any $i_1 , i_2 \in \tilde{I}$ we have the virtual
measures space $\mu(\tilde{F}(i_1) \mid \tilde{F}(i_2))$. We have
$\Aut_{C_2}(E) = \Aut_{C_2}(\tilde{E})$. For any $g \in
\Aut_{C_2}(E)$, for any $i \in I$ we have from the definition of the
group $\Aut_{C_2}(E)$ that $g F(i) = \tilde{F}(p)$ for some $p \in
\tilde{I}$. Therefore for any $g , h \in \Aut_{C_2}(E)$, for any
$i_1 , i_2 \in I$ we have the well-defined $1$-dimensional
$\cc$-vector space \linebreak $ \quad \mu(g F(i_1) \mid h F(i_2))$.

For any $g \in  \Aut_{C_2}(E)$, for any $p \ge q \in \tilde{I}$ we
have the following $\cc$-isomorphism
$$
n_g \quad : \quad \mu (\tilde{F}(p) /\tilde{F}(q)) \arrow{e} \mu
(g\tilde{F}(p) /g\tilde{F}(q)) \mbox{,}
$$
where for any $U \in Gr  (g \tilde{F}(p) /g \tilde{F}(q))$, for any
$\mu \in \mu (\tilde{F}(p) /\tilde{F}(q))$:
$$
n_g (\mu) (\delta_{U}) \eqdef \mu (\delta_{g^{-1} U}) \mbox{.}
$$

For any $g \in \Aut_{C_2}(E)$, for any $p, q, s \in \tilde{I}$,  $s
\le p$, $s \le q$ we have the following $\cc$-isomorphism $m_g$:
$$
 \Hom\nolimits_{\cc} (\mu(\tilde{F}(p)/ \tilde{F}(s))  ,  \mu(\tilde{F}(q)/ \tilde{F}(s)))
\arrow{e} \Hom\nolimits_{\cc} (\mu(g\tilde{F}(p)/ g\tilde{F}(s))  ,
\mu(g\tilde{F}(q)/ g\tilde{F}(s))) \mbox{,}
$$
where for any $f \in \Hom\nolimits_{\cc} (\mu(\tilde{F}(p)/
\tilde{F}(s)) , \mu(\tilde{F}(q)/ \tilde{F}(s)))$:
$$
m_g (f) \eqdef n_g \circ f \circ n_{g^{-1}}  \mbox{.}
 $$

No we apply isomorphisms $m_g$ to the inductive limit  in
formula~(\ref{virm}), then we obtain the following $\cc$-linear
isomorphism for any $p, q \in \tilde{I}$:
$$
l_g \quad : \quad \mu (\tilde{F}(p) \mid \tilde{F}(q)) \arrow{e}
\mu(g\tilde{F}(p) \mid g\tilde{F}(q))  \mbox{.}
$$

We have for any $g_1, g_2 \in G$ that $l_{g_1 g_2} = l_{g_1}
l_{g_2}$. Besides, $\gamma (l_g (a) \otimes l_g(b) )= l_g \gamma (a
\otimes b)$ (see proposition~\ref{pci}) for any $a \in
\mu(\tilde{F}(p) \mid \tilde{F}(q))$, $b \in \mu(\tilde{F}(q) \mid
\tilde{F}(s))$, for any $p,q,s \in \tilde{I}$.

Let $p,q \in \tilde{I}$. Then for any $\mu
\in \mu (\tilde{F}(p) \mid \tilde{F}(q) )$,
$\mu \ne 0$ we define canonically $\mu^{-1}
\in \mu (\tilde{F}(q) \mid \tilde{F}(p) )$
 such that $\mu \otimes \mu^{-1} =1$ with respect to the following
 canonical isomorphism:
 $$
\mu (\tilde{F}(p) \mid \tilde{F}(q) )\otimes_{\cc} \mu (\tilde{F}(q)
\mid \tilde{F}(p) ) = \cc \mbox{.}
 $$

 Let $E= (I,F,V)$ be a $C_2$-space over the field $\df_q$. We fix some $o \in I$.
Then there is the following central extension of groups:
\begin{equation} \label{cextm}
1 \arrow{e} {\cc}^* \arrow{e}
\widehat{\Aut\nolimits_{C_2}(E)}_{F(o)} \arrow{e,t}{\Lambda}
\Aut\nolimits_{C_2}(E) \arrow{e} 1 \mbox{,}
\end{equation}
where
$$
\widehat{\Aut\nolimits_{C_2}(E)}_{F(o)} \eqdef \{ (g, \mu) \: \; :
\; \:
 g \in  \Aut\nolimits_{C_2}(E), \; \mu \in \mu(F(o) \mid g F(o)),
\; \mu \ne 0 \} \mbox{.}
$$
Here $\Lambda((g, \mu)) = g$. The operations $(g_1, \mu_1) \cdot
(g_2, \mu_2) = (g_1 g_2, \gamma(\mu_1 \otimes l_{g_1}(\mu_2)))$ and
$ (g, \mu)^{-1} = (g^{-1}, l_{g^{-1}}(\mu^{-1}))$ define the
structure of a group on the set
$\widehat{\Aut\nolimits_{C_2}(E)}_{F(o)}$. (The unit element of this
group is $(e,1)$, where $e$ is the unit element of the group
$\Aut\nolimits_{C_2}(E)$).

\begin{nt} {\em
For any $o_1 \in I$ there is a canonical isomorphism
$$
\alpha_{o, o_1} \; : \; \widehat{\Aut\nolimits_{C_2}(E)}_{F(o)} \lto
\widehat{\Aut\nolimits_{C_2}(E)}_{F(o_1)} \mbox{.}
$$
Indeed, we fix any $\nu \in \mu(F(o_1) \mid  F(o))$, $\nu \ne 0$.
Then
$$
\alpha_{o, o_1} ((g, \mu)) \eqdef (g, \gamma( \gamma(\nu \otimes
\mu) \otimes l_g(\nu^{-1} ))) \mbox{.}
$$
The map $\alpha_{o, o_1}$ does not depend on the choice of $\nu \in
\mu(F(o_1) \mid  F(o))$, $\nu \ne 0$.}
\end{nt}

\begin{nt} {\em For any  $C_2$-space $E= (I, F, V)$ over a field $k$,
any $o \in I$ it is possible to construct the central extension:
\begin{equation} \label{cext'}
0 \arrow{e} \dz \arrow{e} \widetilde{\Aut\nolimits_{C_2}(E)}_{F(o)}
\arrow{e} \Aut\nolimits_{C_2}(E) \arrow{e} 1 \mbox{.}
\end{equation}
If the field $k = \df_q$, then under the map
$$ a \in \dz  \longmapsto  q^{a} \in \cc^* \mbox{,}$$
 which we apply to the kernel of central extension~(\ref{cext'}),
 central extension~(\ref{cext'}) transfers to central
extension~(\ref{cextm}).

In the case if a $C_2$-space $E$ is constructed from the
two-dimensional local field $k((t_1))((t_2))$, then central
extension~(\ref{cext'}) was constructed in~\cite{O}.}
\end{nt}

\subsubsection{} Let $E=(I,F,V)$ be a $C_2$-space over the finite
field $\df_q$. In section~\ref{ssec} we constructed the maximal
$C_2$-space $\tilde{E} = (\tilde{I}, \tilde{F}, V)$ which dominates
the $C_2$-space $E$. We recall that $G(E)$ is the set of all
$k$-subspaces $W \subset V$ such that $W = \tilde{F}(p)$ for some $p
\in \tilde{I}$. Then any $g \in \Aut_{C_2}(E)$  induces an
isomorphism of the set $G(E)$.

For any $g \in \Aut_{C_2}(E)$, for any  $p \ge q \in \tilde{I}$ we
have the following maps:
\begin{equation} \label{form1}
r_g \quad : \quad    \D(\tilde{F}(p) / \tilde{F}(q)) \arrow{e}
\D(g\tilde{F}(p) / g\tilde{F}(q))
\end{equation}
\begin{equation} \label{form2}
r_g \quad : \quad    \tilde{\E}(\tilde{F}(p) / \tilde{F}(q))
\arrow{e} \tilde{\E}(g\tilde{F}(p) / g\tilde{F}(q))
\end{equation}
\begin{equation}  \label{form3}
r_g \quad : \quad    \E(\tilde{F}(p) / \tilde{F}(q)) \arrow{e}
\E(g\tilde{F}(p) / g\tilde{F}(q)) \mbox{,}
\end{equation}
where $r_g(f) (v) \eqdef f (g^{-1} v)$ for
any $v \in g\tilde{F}(p) / g\tilde{F}(q)$
and any $f$ from $\D(\tilde{F}(p) /
\tilde{F}(q))$, or any $f$ from
$\tilde{\E}(g\tilde{F}(p) / g\tilde{F}(q))$,
or any $f$ from $\E(\tilde{F}(p) /
\tilde{F}(q))$.

Applying formulas~(\ref{form2})-(\ref{form3}) to
formulas~(\ref{2te}) and (\ref{2e})  we obtain the maps for any $g
\in \Aut\nolimits_{C_2}(E)$
$$
r_g \quad : \quad \tilde{\E}(\tilde{E}) \arrow{e}
\tilde{\E}(\tilde{E})
$$
$$
r_g \quad : \quad \E(\tilde{E}) \arrow{e} \E(\tilde{E})
$$
such that $r_g r_h = r_{gh}$ for any $g,h \in \Aut_{C_2}(E)$.

Therefore using section~\ref{ss} we obtain the representations of
the group $\Aut_{C_2}(E)$ by  maps $r_g$ on the $\cc$-spaces
$\tilde{\E}(E)$ and $\E(E)$ such that the first one is a
subrepresentation of the second one.

We fix some $o \in I$. Using formula~(\ref{form1}) we construct the
map $R_{\tilde{g}}$ for any $\tilde{g}= (g, \mu) \in
\widehat{\Aut\nolimits_{C_2}(E)}_{F(o)}$ for any $p \ge q \in
\tilde{I}$:
\begin{equation} \label{rf}
R_{\tilde{g}} \quad : \quad  \D(\tilde{F}(p) / \tilde{F}(q))
\otimes_{\cc} \mu(\tilde{F}(q) \mid F(o)) \arrow{e} \D(g
\tilde{F}(p) / g \tilde{F}(q)) \otimes_{\cc} \mu( g \tilde{F}(q)
\mid F(o))
\end{equation}
as composition of the map $r_g \otimes l_g$ with multiplication by
$\mu^{-1} \in \mu(g F(o) \mid F(o))$.

Applying this formula to formula~(\ref{2d}) we obtain the map for
any $\tilde{g} \in \widehat{\Aut\nolimits_{C_2}(E)}_{F(o)}$:
$$
R_{\tilde{g}} \quad : \quad \D_{F(o)}(\tilde{E}) \arrow{e}
\D_{F(o)}(\tilde{E})
$$
such that $R_{\tilde{g}} R_{\tilde{h}} = R_{\tilde{g} \tilde{h}}$.

Therefore using section~\ref{ss} we obtain the representation of the
group $ \widehat{\Aut\nolimits_{C_2}(E)}_{F(o)}$ by  maps
$R_{\tilde{g}}$ on the $\cc$-space $\D_{F(o)}(E)$.

From constructions we have the following formula:
\begin{equation}  \label{mf2}
R_{\tilde{g}}(f \cdot H) = r_{\Lambda(\tilde{g})}(f) \cdot
R_{\tilde{g}}(H)
\end{equation}
for any $\tilde{g} \in \widehat{\Aut\nolimits_{C_2}(E)}_{F(o)}$, $f
\in  \E(E)$, $H \in \D(E)$.

\vspace{0.5cm}

Dually to formulas~(\ref{form2})-(\ref{form3}), for any $g \in
\Aut_{C_2}(E)$, for any $p \ge q \in \tilde{I}$ we have the
following maps:
\begin{equation}  \label{form3'}
r'_g \quad : \quad    \E'(\tilde{F}(p) / \tilde{F}(q)) \arrow{e}
\E'(g\tilde{F}(p) / g\tilde{F}(q)) \mbox{,}
\end{equation}
\begin{equation} \label{form2'}
r'_g \quad : \quad    \tilde{\E}'(\tilde{F}(p) / \tilde{F}(q))
\arrow{e} \tilde{\E}'(g\tilde{F}(p) / g\tilde{F}(q)) \mbox{.}
\end{equation}

Applying formulas~(\ref{form3'})-(\ref{form2'}) to
formulas~(\ref{2e'})-(\ref{2te'}), and using section~\ref{ss}, we
obtain the representations of the group $\Aut_{C_2}(E)$  on the
$\cc$-spaces $\tilde{\E}'(E)$   and $\E'(E)$ by  maps $r_g'$   such
that the first representation is a subrepresentation of the second
one. By construction, we have for any $ f_1 \in \E(E)$, any $H_1 \in
\E'(E)$, and any $ g \in \Aut_{C_2}(E)$
$$
<r_g f_1 \:,\: r_g' H_1>_{\E(E)} = <f_1 \:,\: H_1>_{\E(E)}  \mbox{.}
$$
We have for any $ f_2 \in \tilde{\E}(E)$, any $H_2 \in
\tilde{\E}'(E)$, and any $ g \in \Aut_{C_2}(E)$
$$
<r_g' H_2 \:,\: r_g f_2>_{\tilde{\E}(E)} = <H_2 \:, \:
f_2>_{\tilde{\E}(E)} \mbox{.}
$$

Let $o \in I$. Dually to formula~(\ref{rf})
 we construct the map $R'_{\tilde{g}}$ for any $\tilde{g}= (g, \mu)
\in \widehat{\Aut\nolimits_{C_2}(E)}_{F(o)}$ for any $p \ge q \in
\tilde{I}$:
\begin{equation} \label{rf'}
R'_{\tilde{g}} \quad : \quad  \D'(\tilde{F}(p) / \tilde{F}(q))
\otimes_{\cc} \mu(F(o) \mid \tilde{F}(q) ) \arrow{e} \D'(g
\tilde{F}(p) / g \tilde{F}(q)) \otimes_{\cc} \mu( F(o) \mid  g
\tilde{F}(q))  \mbox{.}
\end{equation}

Applying formula~(\ref{rf'}) to formula~(\ref{2d'}), and using
section~\ref{ss}, we obtain the representation of the group $
\widehat{\Aut\nolimits_{C_2}(E)}_{F(o)}$  on the $\cc$-space
$\D'_{F(o)}(E)$ by maps $R'_{\tilde{g}}$. By construction, we have
the following formula for any $H \in \D'_{F(o)}(E)$, any $f \in
\D_{F(o)}(E)$, any $\tilde{g} \in
\widehat{\Aut\nolimits_{C_2}(E)}_{F(o)} $:
\begin{equation} \label{mf3}
< R'_{\tilde{g}}(H)  \:,\:  R_{\tilde{g}}(f)  >_{\D_{F(o)}(E)} = < H
 \:,\: f>_{\D_{F(o)}(E)} \mbox{.}
\end{equation}

From formulas~(\ref{mf3}), (\ref{mf2}), and~(\ref{mf1}) we obtain
the following formula for any $f \in \E(E)$, any $H \in
\D'_{F(o)}(E)$, and any $\tilde{g} \in
\widehat{\Aut\nolimits_{C_2}(E)}_{F(o)} $:
$$
R'_{\tilde{g}} (f \cdot H) = R_{\Lambda(\tilde{g})} (f) \cdot
R'_{\tilde{g}} (H)  \mbox{.}
$$

\subsubsection{}
\label{sec5.5.5}
Let $E = (I,F,V)$ be a complete $C_2$-space over the field $\df_q$.
Then $\check{E} = (I^0, F^0, \check{V})$, and $\check{\check{E}}=E$.
For any $g \in \Aut_{C_2}(E)$ we have canonically $\check{g} \in
\Aut_{C_2}(\check{E})$. And we have the following isomorphism of
groups:
\begin{equation} \label{eee1}
g \in \Aut\nolimits_{C_2}(E) \longmapsto \check{g}^{-1} \in
\Aut\nolimits_{C_2}(\check{E}) \mbox{.}
\end{equation}

Also for any $o \in I$ we have the following isomorphism of groups:
$$
(g, \mu) \in \widehat{\Aut\nolimits_{C_2}(E)}_{F(o)} \longmapsto
(\check{g}^{-1}, \mu)  \in
\widehat{\Aut\nolimits_{C_2}(\check{E})}_{F^{0}(o)} \mbox{,}
$$
where we use isomorphism~(\ref{fff}) to obtain that for any $g \in
\Aut_{C_2}(E)$ canonically
\begin{equation}  \label{eee2}
\mu(F(o) \mid g F(o)) = \mu(F^0(o) \mid \check{g}^{-1} F^0(o) )
\mbox{.}
\end{equation}

Using isomorphism~(\ref{eee1}) we obtain representations of the
group~$\Aut_{C_2}(E)$ on $\cc$-spaces $\E(\check{E})$,
$\tilde{\E}(\check{E})$, $\E'(\check{E})$, $\tilde{\E}'(\check{E})$.

Using isomorphism~(\ref{eee2}) we obtain representations of the
group $\widehat{\Aut\nolimits_{C_2}(E)}_{F(o)}$ on $\cc$-spaces
$\D_{F^0(o)}(\check{E})$ and $\D'_{F^0(o)}(\check{E})$.

We have the following proposition.
\begin{prop}
Let $E = (I,F,V)$ be a complete $C_2$-space over the field $\df_q$.
We fix some $o \in I$. Then
\begin{enumerate}
\item The Fourier transform $\F$ gives an isomorphism between
representations of the group $\Aut_{C_2}(E)$ on $\cc$-spaces $\E(E)$
and $\tilde{\E}' (\check{E})$ and on $\cc$-spaces $\tilde{\E}(E)$
and $\E' (\check{E})$.
\item The Fourier transform $\F$ gives an isomorphism between
representations of the group
$\widehat{\Aut\nolimits_{C_2}(E)}_{F(o)}$ on $\cc$-spaces
$\D_{F(o)}(E)$ and $\D_{F^0(o)} (\check{E})$ and on $\cc$-spaces
$\D'_{F(o)}(E)$ and $\D'_{F^0(o)} (\check{E})$.
\end{enumerate}
\end{prop}
\proof. Using definition of $2$-dimensional
Fourier transform $\F$ from
section~\ref{2Ft}, and using
formulas~(\ref{2d})-(\ref{2te'}), we reduce
the statements of this proposition to
corresponding statements about isomorphisms
of $C_1$-spaces and $1$-dimensional Fourier
transforms between them. The last statements follow from
the definition and properties of
$1$-dimensional Fourier transform, or can
be again reduced to $0$-dimensional case by
formulas~(\ref{d})-(\ref{te'}) using
remark~\ref{nt5}. The corresponding
statement in the $0$-dimensional case looks
as following. For any finite-dimensional
vector space $W$ over the field $\df_q$, for
any $f \in \ff(W)$, for any $u \in W^*$, for
any $g \in \Aut_{\sdf_q}(W)$ we have:
$$
\sum\limits_{w \in W} f(g^{-1}w) \overline{\psi (u(w))} =
\sum\limits_{w \in W} f(w) \overline{\psi (u(gw))} \mbox{.}
 $$
\vspace{0.3cm}

\subsection{Direct and inverse images}
\label{secdi} In this section we construct direct and inverse images
of spaces of functions and distributions on $C_2$-spaces.
\subsubsection{}
Let
$$
0\arrow{e} E_1\arrow{e,t}{\alpha} E_2\arrow{e,t}{\beta}E_3\arrow{e}0
$$
be an admissible triple of $C_2$-spaces over a finite field $\df_q$,
where $E_i=(I_i,F_i,V_i)$, \linebreak $1 \le  i \le 3$. By
definition, there are  order-preserving  functions:
\begin{itemize}
\item[(i)] $ \gamma \; \colon \; I_2\arrow{e} I_3\quad\mbox{such
that}\quad \beta(F_2(i))=F_3(\gamma(i))\quad\mbox{for any} ~i\in I_2
\mbox{,}
$
\item[(ii)] $ \veps \; \colon \; I_2\arrow{e} I_1\quad\mbox{such that}\quad
F_2(i)\cap V_1=F_1(\veps(i))\quad\mbox{for any} ~i\in I_2.
$
\end{itemize}

\begin{prop} \label{pr26}
We suppose that $E_1$ is a $cC_2$-space. Let $o\in I_2$. Then there
is the direct image
$$
\beta_* \quad \colon \quad
\D_{F_{2}(o)}(E_2)\otimes_{\cc}\mu(F_1(\veps(o))\mid V_1)\arrow{e}
\D_{F_{3}(\gamma(o))}(E_3).
$$
\end{prop}
\proof. For any $i,j\in I_2$ we have canonically
\begin{equation} \label{*}
\mu(F_2(i)\mid F_2(j))=\mu(F_1(\veps(i))\mid
F_1(\veps(j)))\otimes_{\cc} \mu(F_3(\gamma(i))\mid F_3(\gamma(j)))
\mbox{.}
\end{equation}

Let $k\in I_2$ such that
$F_1(\veps(k))=V_1$. Let $i\ge j\in I_2$ be
any such that $i\ge k\ge j$. Then we have an
admissible triple of $C_1$-spaces:
$$
0\arrow{e} V_1/F_1(\veps(j))\arrow{e}
F_2(i)/F_2(j)\arrow{e,t}{\beta_{ij}}
F_3(\gamma(i))/F_3(\gamma(j))\arrow{e} 0.
$$

We have a well-defined map:
$$
(\beta_{ij})_*   \quad \colon  \quad
\D(F_2(i)/F_2(j))\otimes_{\cc}\mu(F_1(\veps(j))\mid V_1)\arrow{e}
\D(F_3(\gamma(i))/F_3(\gamma(j))),
$$
where we use that $\mu(F_1(\veps(j))\mid V_1)=\mu(V_1 /
F_1(\veps(j)))$.

Therefore we have a well-defined map
$$
(\beta_{ij})_* \otimes \mu(F_3(\gamma(j))\mid F_3(\gamma(o))) \quad
\colon
$$
$$
 \D(F_2(i)/F_2(j)) \otimes_{\cc} \mu(F_1(\veps(j)) \mid
V_1)\otimes_{\cc} \mu(F_3(\gamma(j))\mid F_3(\gamma(o)))\arrow{e}
\qquad \qquad \qquad \qquad
$$
$$
\qquad \qquad \qquad \qquad \qquad \qquad \arrow{e}
\D(F_3(\gamma(i))\mid
F_3(\gamma(j)))\otimes_{\cc}\mu(F_3(\gamma(j))\mid F_3(\gamma(o)))
\mbox{.}
$$
From (\ref{*}) we have that
$$
\mu(F_3(\gamma(j))\mid F_3(\gamma(o)) = \mu(F_2(j) \mid
 F_2(o)) \otimes_{\cc}  \mu(F_1(\veps(o))\mid F_1(\veps(j))).
$$
Therefore
$$
\mu(F_1(\veps(j))\mid V_1)\otimes_{\cc}\mu(F_3(\gamma(j))\mid
F_3(\gamma(o)))=\mu(F_2(j) \mid
F_2(o))\otimes_{\cc}\mu(F_1(\veps(o))\mid V_1).
$$

Thus we have a well-defined map
$$
(\beta_{ij})_*\otimes\mu(F_3(\gamma(j))\mid F_3(\gamma(o)))  \quad
\colon
$$
$$
\D(F_2(i)/F_2(j))\otimes_{\cc}\mu(F_2(j) \mid
F_2(o))\otimes_{\cc}\mu(F_1(\veps(o))\mid V_1)\arrow{e}  \qquad
\qquad \qquad \qquad \qquad
$$
$$
\qquad \qquad \qquad \qquad \qquad \qquad \arrow{e}
\D(F_3(\gamma(i)) /
F_3(\gamma(j)))\otimes_{\cc}\mu(F_3(\gamma(j))\mid F_3(\gamma(o)))
\mbox{.}
$$
Now we take the projective limits which are used to construct spaces
$\D_{F_2(o)}(E_2)$ and $D_{F_3(\gamma(o))}(E_3)$. Then we obtain a
well-defined map $\beta_*$.
\vspace{0.3cm}

\begin{nt} \em{
Dually to the maps above, we obtain a well-defined map
$$
\beta^{*} \quad \colon \quad
\D'_{F_3(\gamma(o))}(E_3)\otimes_{\cc}\mu(F_3(\veps(o))\mid
V_1)\arrow{e} \D'_{F_2(o)}(E_2)
$$
such that the maps $\beta_*$ and $\beta^*$ are conjugate maps with
respect to the pairings
\linebreak
$<\cdot,\cdot>_{\D_{F_2(o)}(E_2)}$ and
$<\cdot,\cdot>_{\D_{F_3(\gamma(o))}(E_3)}$.}
\end{nt}

\subsubsection{}
Let
$$
0\arrow{e} E_1\arrow{e,t}{\alpha} E_2\arrow{e,t}{\beta}E_3\arrow{e}0
$$
be an admissible triple of $C_2$-spaces over a finite field $\df_q$,
where $E_i=(I_i,F_i,V_i)$, \linebreak$1 \le i \le 3$. By definition,
there are order-preserving  functions:
\begin{itemize}
\item[(i)]
$ \gamma \; \colon \; I_2 \arrow{e} I_3 \quad\mbox{such that}\quad
\beta(F_2(i))=F_3(\gamma(i))\quad\mbox{for any} \quad i\in I_2
\mbox{,} $
\item[(ii)]
$
\veps \; \colon \; I_2 \arrow{e} I_1 \quad\mbox{such that}\quad
F_2(i)\cap V_1=F_1(\veps(i))\quad\mbox{for any}\quad i\in I_2
\mbox{.}
$
\end{itemize}

\begin{prop} \label{pr27}
We suppose that $E_3$ is a $dC_2$-space. Let $o \in I_2$. Then there
is the inverse image
$$
\alpha^* \quad \colon \quad
\D_{F_{2}(o)}(E_2)\otimes_{\cc}\mu(F_3(\gamma(o))\mid
\{0\})\arrow{e} \D_{F_{1}(\veps(o))}(E_1) \mbox{,}
$$
where $\{0\}$ is  the zero subspace of $V_3$.
\end{prop}
\proof. For any $i,j\in I_2$ we have canonically
\begin{equation} \label{**}
\mu(F_2(i)\mid F_2(j))=\mu(F_1(\veps(i))\mid
F_1(\veps(j)))\otimes_{\cc} \mu(F_3(\gamma(i))\mid F_3(\gamma(j)))
\mbox{.}
\end{equation}

Since $\beta$ is an admissible epimorphism and $E_3$ is a
$dC_2$-space, there is $k\in I_2$ such that $F_3(\gamma(k))=\{0\}$.
Let $i\ge j\in I_2$ be any such that $i\ge k\ge j$. Then we have an
admissible triple of $C_1$-spaces:
\begin{equation} \label{***}
0\arrow{e} F_1(\veps(i))/F_1(\veps(j))
\arrow{e,t}{\alpha_{ij}}F_2(i)/F_2(j)\arrow{e}
F_3(\gamma(i))\arrow{e}0 \mbox{.}
\end{equation}

We have a well-defined map
$$
\alpha_{ij}^* \quad \colon \quad \D(F_2(i)/F_2(j))\arrow{e}
\D(F_1(\veps(i))/F_1(\veps(j))).
$$

Therefore we have a well-defined map
$$
\alpha_{ij}^* \otimes \mu(F_1(\veps(j))\mid F_1(\veps(o))) \quad
\colon
$$
$$
\D(F_2(i)/F_2(j))\otimes_{\cc}\mu(F_1(\veps(j))\mid
F_1(\veps(o)))\arrow{e} \qquad \qquad \qquad \qquad \qquad \qquad $$
$$
\qquad \qquad \qquad \qquad \qquad \qquad
\arrow{e}
 \D(F_1(\veps(i))\mid
F_1(\veps(j)))\otimes_{\cc}\mu(F_1(\veps(j))\mid
F_1(\veps(o)))\mbox{.}
$$
From (\ref{**}) and (\ref{***}) we have that
$$
\mu(F_2(o) \mid F_2(j))=\mu(F_1(\veps(o))\mid
F_1(\veps(j)))\otimes_{\cc}\mu(F_3(\gamma(o))\mid \{0\}) \mbox{.}
$$
Hence
$$
\mu(F_1(\veps(j))\mid F_1(\veps(o)))=\mu(F_2(j) \mid
F_2(o))\otimes_{\cc}\mu(F_3(\gamma(o))\mid \{0\}) \mbox{.}
$$
Thus we have a well-defined map
$$
\alpha_{ij}^*\otimes\mu(F_1(\veps(j))\mid F_1(\veps(o))) \quad
\colon
$$
$$
\D(F_2(i)/F_2(j))\otimes_{\cc}\mu(F_2(j)/F_2(o))\otimes_{\cc}\mu(F_3(\gamma(o))\mid
\{0\})\arrow{e} \qquad \qquad \qquad \qquad \qquad \qquad
$$
$$
\qquad \qquad \qquad \qquad \qquad \qquad
 \arrow{e} \D(F_1(\veps(i))\mid
F_1(\veps(j)))\otimes_{\cc}\mu(F_1(\veps(j))\mid
F_1(\veps(o)))\mbox{.}
$$
Now we take the projective limits (with respect to $i\ge k\ge j$)
which we used to construct spaces $\D_{F_2(o)}(E_2)$ and
$\D_{F_1(\veps(o))}(E_1)$. Then we obtain a well-defined map
$\alpha^*$.
\vspace{0.3cm}

\begin{nt} {\em
Dually to the maps above, we obtain a well-defined map
$$
\alpha_{*} \quad \colon \quad
\D'_{F_1(\veps(o))}(E_1)\otimes_{\cc}\mu(F_3(\gamma(o))\mid
\{0\})\arrow{e} \D'_{F_2(o)}(E_2)
$$
such that the maps $\alpha^*$ and $\alpha_*$ are conjugate maps with
respect to the pairings
\linebreak
$<\cdot,\cdot>_{\D_{F_2(o)}(E_2)}$ and
$<\cdot,\cdot>_{\D_{F_1(\veps(o))}(E_1)}$. }
\end{nt}

\subsubsection{}
Let
$$
0\arrow{e} E_1\arrow{e,t}{\alpha} E_2\arrow{e,t}{\beta}E_3\arrow{e}0
$$
be an admissible triple of $C_2$-spaces over a finite field $\df_q$,
where $E_i=(I_i,F_i,V_i)$, \linebreak $1 \le i \le 3$. By
definition, there are order-preserving  functions:
\begin{itemize}
\item[(i)]
$ \gamma \; \colon \; I_2\arrow{e} I_3\quad\mbox{such that}\quad
\beta(F_2(i))=F_3(\gamma(i))\quad\mbox{for any}\quad i\in I_2
\mbox{,} $
\item[(ii)]
$ \veps \; \colon \; I_2\arrow{e} I_1\quad\mbox{such that}\quad
F_2(i)\cap V_1=F_1(\veps(i))\quad\mbox{for any}\quad i\in I_2
\mbox{.} $
\end{itemize}

\begin{prop} \label{pr28}
We suppose that $E_1$ is a $cfC_2$-space. Let $o \in I_2$. Then
there is the inverse image
$$
\beta^* \quad \colon \quad \D_{F_{3}(\gamma(o))}(E_3)\arrow{e}
\D_{F_{2}(o)}(E_2).
$$
\end{prop}
\proof. For any $i,j\in I_2$ we have canonically
\begin{equation} \label{sharp}
\mu(F_2(i)\mid F_2(j))=\mu(F_1(\veps(i))\mid
F_1(\veps(j)))\otimes_{\cc} \mu(F_3(\gamma(i))\mid F_3(\gamma(j)))
\mbox{.}
\end{equation}

Since $E_1$ is a $cfC_2$-space, for any $i\ge j\in I_2$ we have a
canonical element $1_{ij}\in \mu(F_1(\veps(i)) / F_1(\veps(j)))$
such that $1_{ij}(\delta_{F_1(\veps(i))\mid F_1(\veps(j))})=1$.

Therefore we have a canonical element $1_{ij}\in \mu(F_1(\veps(i))
\mid F_1(\veps(j)))$ for any \linebreak $i,j\in I_2$ such that
$1_{ij}\otimes 1_{jk}=1_{ik}$ for any $i,j,k\in I_2$. Thus,  we have
canonically that  $\mu(F_1(\veps(i))\mid F_1(\veps(j)))=\cc$ for any
$i,j\in I_2$. (See also remark~\ref{remark9}).

Thus from formula (\ref{sharp}) we have for any $i,j\in I_2$
\begin{equation} \label{sharpsharp}
\mu(F_2(i)\mid F_2(j))=\mu(F_3(\gamma(i))\mid F_3(\gamma(j)))
\mbox{.}
\end{equation}

For any $i\ge j\in I_2$ we have an admissible triple of $C_1$-spaces:
$$
0\arrow{e}
\frac{F_1(\veps(i))}{F_1(\veps(j))}\arrow{e}\frac{F_2(i)}{F_2(j)}
\arrow{e,t}{\beta_{ij}}
\frac{F_3(\gamma(i))}{F_3(\gamma(j))}\arrow{e}0 \mbox{.}
$$

We have a well-defined map:
$$
\beta_{ij}^* \quad \colon \quad \D(F_3(\gamma(i))/F_3(\gamma(j)))
\arrow{e} \D(F_2(i)/F_2(j)) \mbox{.}
$$

From formula (\ref{sharpsharp}) we have that
$$
\mu(F_2(j) \mid F_2(o)) = \mu(F_3(\gamma(j))\mid
F_3(\gamma(o)))\mbox{.}
$$

Therefore we have a map:
$$
\beta_{ij}^* \otimes \mu(F_2(j)\mid F_2(o)) \quad \colon
$$
$$
\D(F_3(\gamma(i)) /
F_3(\gamma(j)))\otimes_{\cc}\mu(F_3(\gamma(j))\mid F_3(\gamma(o)))
\arrow{e} \qquad \qquad \qquad \qquad \qquad \qquad \qquad \qquad
\qquad
$$
$$
\qquad \qquad \qquad \qquad \qquad \qquad \qquad \qquad \qquad
 \arrow{e} \D(F_2(i) / F_2(j))\otimes_{\cc}\mu(F_2(j)\mid
F_2(o)) \mbox{.}
$$
These maps (for various $i\ge j\in I_2$) are  compatible when we
take the projective limits according to formulas which  define
$\D_{F_3(\gamma(o))}(E_3)$ and $\D_{F_2(o)}(E_2)$. Therefore we
obtain the map $\beta^*$.
\vspace{0.3cm}

\begin{nt} {\em
Dually to the maps, which we considered above, we have a
well-defined map
$$
\beta_{*} \quad \colon \quad \D'_{F_2(o)}(E_2)\arrow{e}
D'_{F_3(\gamma(o))}(E_3)
$$
such that the maps $\beta^*$ and $\beta_*$ are conjugate maps with
respect to the pairings \linebreak$<\cdot,\cdot>_{\D_{F_2(o)}(E_2)}$
and $<\cdot,\cdot>_{\D_{F_3(\gamma(o))}(E_3)}$. }
\end{nt}

\subsubsection{}
 Let
$$
0\arrow{e} E_1\arrow{e,t}{\alpha} E_2\arrow{e,t}{\beta}E_3\arrow{e}0
$$
be an admissible triple of $C_2$-spaces over a finite field $\df_q$,
where $E_i=(I_i,F_i,V_i)$, \linebreak $1 \le i \le 3$. By
definition, there are order-preserving functions:
\begin{itemize}
\item[(i)]
$ \gamma \; \colon \; I_2\arrow{e} I_3\quad\mbox{such that}\quad
\beta(F_2(i))=F_3(\gamma(i))\quad\mbox{for any}\quad i\in I_2
\mbox{,} $
\item[(ii)]
$
\veps \; \colon \; I_2\arrow{e} I_1\quad\mbox{such that}\quad
F_2(i)\cap V_1=F_1(\veps(i))\quad\mbox{for any}\quad i\in I_2
\mbox{.}
$
\end{itemize}

\begin{prop} \label{pr29}
We suppose that $E_3$ is a $dfC_2$ space. Let $o\in I_2$. Then there
is the direct image
$$
\alpha_* \quad \colon \quad \D_{F_{1}(\veps(o))}(E_1)\arrow{e}
\D_{F_{2}(o)}(E_2) \mbox{.}
$$
\end{prop}
\proof. For any $i,j\in I_2$ we have canonically
\begin{equation} \label{natural}
\mu(F_2(i)\mid F_2(j))=\mu(F_1(\veps(i))\mid
F_1(\veps(j)))\otimes_{\cc} \mu(F_3(\gamma(i))\mid F_3(\gamma(j)))
\mbox{.}
\end{equation}
Since $E_3$ is a $dfC_2$-space, for any $i \ge j \in I_2$ we have a
canonical element \linebreak $\delta_{ij} \in \mu(F_3(\gamma(i)) /
F_3(\gamma(j)))$ such that $\delta_{ij}(\{0\})=1$. Therefore we have
a canonical ele\-ment  $\delta_{ij}\in \mu(F_3(\gamma(i))\mid
F_3(\gamma(j)))$ for any $i,j\in I_2$ such that $\delta_{ij}\otimes
\delta_{jk}=\delta_{ik}$ for any $i,j,k \in I_2$. Therefore
 we have canonically $\mu(F_3(\gamma(i)) \mid
F_3(\gamma(j)))=\cc$ for any $i,j\in I_2$. (See also
remark~\ref{remark9}).

Thus from formula (\ref{natural}) we have for any $i,j\in I_2$
\begin{equation} \label{naturalnatural}
\mu(F_2(i)\mid F_2(j))=\mu(F_1(\veps(i))\mid F_1(\veps(j))) \mbox{.}
\end{equation}

For any $i\ge j\in I_2$ we have an admissible triple of $C_1$-spaces:
$$
0\arrow{e}
\frac{F_1(\veps(i))}{F_1(\veps(j))}\arrow{e,t}{\alpha_{ij}}\frac{F_2(i)}{F_2(j)}
\arrow{e} \frac{F_3(\gamma(i))}{F_3(\gamma(j))}\arrow{e}0 \mbox{.}
$$

We have a well-defined map:
$$
(\alpha_{ij})_* \quad \colon \quad
\D(F_1(\veps(i))/F_1(\veps(j)))\arrow{e} \D(F_2(i)/F_2(j)),
$$
because $\frac{F_3(\gamma(i))}{F_3(\gamma(j))}$ is a discrete
$C_1$-space.

From formula (\ref{naturalnatural}) we have that
$$
\mu(F_2(j) \mid F_2(o))=\mu(F_1(\veps(j)) \mid F_1(\veps(o))).
$$

Therefore we have a map:
$$
(\alpha_{ij})_*\otimes \mu(F_2(j) \mid F_2(o)) \quad \colon
$$
$$
\D(F_1(\veps(i))/F_1(\veps(j)))\otimes_{\cc}\mu(F_1(\veps(j))\mid
F_1(\veps(o))) \arrow{e} \qquad \qquad \qquad \qquad \qquad \qquad
\qquad \qquad
$$
$$
\qquad \qquad \qquad \qquad \qquad \qquad \qquad \qquad \arrow{e}
\D(F_2(i)/F_2(j)) \otimes_{\cc} \mu(F_2(j)\mid F_2(o)).
$$
These maps (for various $i\ge j\in I_2$) are  compatible when we
take the projective limits according to the formulas which  define
$\D_{F_1(\veps(o))}(E_1)$ and $\D_{F_2(o)}(E_2)$. Therefore we
obtain the map $\alpha_*$.
\vspace{0.3cm}

\begin{nt} {\em
Dually to the maps, which we considered above, we have a
well-defined map
$$
\alpha^{*} \quad \colon \quad \D'_{F_2(o)}(E_2) \arrow{e}
\D'_{F_1(\veps(o))}(E_1)
$$
such that the maps $\alpha_*$ and $\alpha^*$ are conjugate maps with
respect to the pairings \linebreak
$<\cdot,\cdot>_{\D_{F_2(o)}(E_2)}$ and
$<\cdot,\cdot>_{\D_{F_1(\veps(o))}(E_1)}$. }
\end{nt}

\subsection{Composition of maps and base change rules.} \label{bcr2}
In this section we consider  base change rules
(subsection~\ref{ss5.7.1}) and rules of composition of maps
(subsection~\ref{ss5.7.2}).

\subsubsection{} \label{ss5.7.1}
Let
$$
0 \arrow{e} E_1  \arrow{e,t}{\alpha}  E_2 \arrow{e,t}{\beta}   E_3
\arrow{e} 0
$$
be an admissible triple of $C_2$-spaces over the field $\df_q$,
where $E_i = (I_i, F_i, V_i)$, $ 1 \le i \le 3$.

Let
$$
0 \arrow{e} D  \arrow{e,t}{\gamma}  E_3 \arrow{e,t}{\delta}   B
\arrow{e} 0
$$
be an admissible triple of  $C_2$-spaces over the field $\df_q$,
where $D = (R, S, Y)$, $B = (T, U, W)$. Then we have the following
commutative diagram of $C_2$-spaces.

\begin{equation} \label{zvezda}
\begin{diagram}
\node[3]{0} \arrow{s} \node{0} \arrow{s} \\
 \node{0} \arrow{e} \node{E_1}
\arrow{e,t}{\gamma_{\alpha}} \arrow{s,=} \node{E_2
\mathop{\times}_{E_3} D} \arrow{e,t}{\gamma_{\beta}}
\arrow{s,l}{\beta_{\gamma}} \node{D} \arrow{e} \arrow{s,r}{\gamma} \node{0} \\
\node{0} \arrow{e} \node{E_1}  \arrow{e,t}{\alpha}  \node{E_2}
\arrow{e,t}{\beta}  \arrow{s,l}{\beta_{\delta}} \node{E_3} \arrow{e}
\arrow{s,r}{\delta} \node{0} \\
\node[3]{B} \arrow{r,=} \arrow{s} \node{B} \arrow{s} \\
\node[3]{0} \node{0}
\end{diagram}
\end{equation}

Let $X' = E_2 \mathop{\times}\limits_{E_3} D = (N, Q, X)$ as
$C_2$-space. In this diagram two horizontal triples and two vertical
triples are admissible triples of  $C_2$-spaces.

We choose any  $i \le j \in I_2$. Then diagram~(\ref{zvezda})
induces the following commutative diagram of  $C_1$-spaces (see the
construction of fibered product from lemma~\ref{lbc2}).

\begin{equation} \label{2zvezda}
\begin{diagram}
\node[3]{0} \arrow{s} \node{0} \arrow{s} \\
 \node{0} \arrow{e}
 \node{\frac{F_2(j) \cap V_1}{F_2(i) \cap V_1}}
\arrow{e,t}{(\gamma_{\alpha})_{ji}} \arrow{s,=} \node{\frac{F_2(j)
\cap X}{F_2(i) \cap X}} \arrow{e,t}{(\gamma_{\beta})_{ji}}
\arrow{s,l}{(\beta_{\gamma})_{ji}} \node{\frac{\beta(F_2(j)) \cap
Y}{\beta(F_2(i)) \cap Y}}
\arrow{e} \arrow{s,r}{\gamma_{ji}} \node{0} \\
\node{0} \arrow{e} \node{\frac{F_2(j) \cap V_1}{F_2(i) \cap V_1}}
\arrow{e,t}{\alpha_{ji}}  \node{\frac{F_2(j)}{F_2(i)}}
\arrow{e,t}{\beta_{ji}} \arrow{s,l}{(\beta_{\delta})_{ji}}
\node{\frac{\beta(F_2(j))}{\beta(F_2(i))}} \arrow{e}
\arrow{s,r}{\delta_{ji}} \node{0} \\
\node[3]{\frac{\delta \beta (F_2(j))}{\delta \beta (F_2(i))}}
\arrow{r,=} \arrow{s}
\node{\frac{\delta \beta (F_2(j))}{\delta \beta (F_2(i))}} \arrow{s} \\
\node[3]{0} \node{0}
\end{diagram}
\end{equation}

Here two horizontal triples and two vertical triples are admissible
triples of
 $C_1$-spaces. At the same time,
$$
\frac{F_2(j) \cap X}{F_2(i) \cap X} = \frac{F_2(j)}{F_2(i)}
\mathop{\times}_{\frac{\beta(F_2(j))}{\beta(F_2(i))}}
\frac{\beta(F_2(j)) \cap Y}{\beta(F_2(i)) \cap Y}
$$
as  $C_1$-spaces.

\begin{prop}. \label{przam1}
Using notations of diagram~(\ref{zvezda}),  let $E_1$ be a
$cC_2$-space, $B$ be a $dC_2$-space. Let $o \in I_2$, $\mu \in
\mu(F_2(o) \cap V_1 \mid V_1)$, $ \nu \in \mu(\delta \beta (F_2(o))
\mid \{ 0\})$. Then the following formulas are satisfied.
\begin{enumerate}
\item
For any $f \in  \D_{F_2(o)} (E_2)$
\begin{equation} \label{ut1}
\gamma^*( \beta_*(f \otimes \mu) \otimes \nu )= (\gamma_{\beta})_*
(\beta_{\gamma}^* (f \otimes \nu) \otimes \mu) \mbox{.}
\end{equation}
\item
For any  $G \in \D'_{\beta(F_2(o))\cap Y} (D)$
\begin{equation} \label{ut2}
\beta^* ( \gamma_* (G \otimes \nu) \otimes \mu) =
(\beta_{\gamma})_*(\gamma_{\beta}^* (G \otimes \mu) \otimes \nu)
\mbox{.}
\end{equation}
\end{enumerate}
\end{prop}
\begin{nt}. {\em
The formulation of this proposition is well-defined  (with respect
to definitions of direct and inverse images), since $\alpha =
\beta_{\gamma} \circ \gamma_{\alpha}$, $\beta_{\delta} = \delta
\circ \beta $.}
\end{nt}
\proof  of proposition~\ref{przam1}. We prove at first
formula~(\ref{ut1}). Let $j_0 \in I_2$ be such that $F_2(j_0) \cap
V_1 = V_1$. (Such a $j_0$ it is possible to choose, since $E_1$ is a
$cC_2$-space.) Let $i_0 \in I_2$ be such that  $\beta \delta
(F_2(i_0)) = \{ 0\}$. (Such an $i_0$ it is possible to choose, since
 $B$ is a $dC_2$-space.)
We consider arbitrary  $j \ge i \in I_2$ such that
  $j \ge j_0$, $i
\le i_0$. Then  the base change formula on $C_1$-spaces is satisfied
for the maps $\beta_{ji}$, $\gamma_{ji}$, $(\gamma_{\beta})_{ji}$,
$(\beta_{\gamma})_{ji}$ from cartesian square of
diagram~(\ref{2zvezda}) (see formula~(\ref{e''1}) from
proposition~\ref{pr18}). We multiply this formula by corresponding
measure spaces, and then we take the projective limit with respect
to all such $j \ge i \in I_2$ ($j \ge j_0$, $i \le i_0$). Using the
constructions of maps $\beta_*$ and $\gamma^*$ from
propositions~\ref{pr26} and~\ref{pr27},  we obtain the base change
formula~(\ref{ut1}).

Formula~(\ref{ut2}) is a dual formula to formula~(\ref{ut1}), and it
can be obtained by analogous reasonings as in the proof of
formula~(\ref{ut1}), but one has to take the inductive limit with
respect to $j, i \in I_2$, and one has to use formula~(\ref{e''2})
instead of formula~(\ref{e''1}).
\vspace{0.3cm}

\begin{prop}. \label{przam2}
Using notations of diagram~(\ref{zvezda}), let   $E_1$ be a
$cfC_2$-space, $B$ be a $dfC_2$-space. Let $o \in I_2$. Then the
following formulas are satisfied.
\begin{enumerate}
\item For any $f \in \D_{\beta(F_2(o)) \cap Y} (D)$
\begin{equation} \label{ut3}
\beta^* \gamma_* (f) = (\beta_{\gamma})_* \gamma_{\beta}^* (f)
\mbox{.}
\end{equation}
\item
For any $G \in \D'_{F_2(o)}(E_2)$
\begin{equation} \label{ut4}
\gamma^* \beta_* (G) = (\gamma_{\beta})_* \beta_{\gamma}^* (G)
\mbox{.}
\end{equation}
\end{enumerate}
\end{prop}
\proof. We prove formula~(\ref{ut3}). We consider arbitrary  $j \ge
i \in I_2$. Then  the base change formula on $C_1$-spaces is
satisfied for the maps $\beta_{ji}$, $\gamma_{ji}$,
$(\gamma_{\beta})_{ji}$, $(\beta_{\gamma})_{ji}$ from cartesian
square of diagram~(\ref{2zvezda}) (see formula~(\ref{pr18f5}) from
proposition~\ref{pr18}). We take the projective limit with respect
to  such $j \ge i \in I_2$. Using the constructions of maps
$\beta^*$ and
 $\gamma_*$ from propositions~\ref{pr28}
and~\ref{pr29},  we obtain the base change formula~(\ref{ut3}).

Formula~(\ref{ut4}) can be obtained by analogous reasonings, but one
has to take the inductive limit with respect to $ i, j \in I_2$, and
one has to use formula~(\ref{pr18f6}) instead of
formula~(\ref{pr18f5}).
\vspace{0.3cm}

\begin{prop}.
Using notations of diagram~(\ref{zvezda}), let   $E_1$ be a
$cC_2$-space, $B$ be a $dfC_2$-space. Let $o \in I_2$, $\mu \in
\mu(F_2(o) \cap V_1 \mid V_1)$. Then the following formulas are
satisfied.
\begin{enumerate}
\item For any $f \in \D_{F_2(o) \cap X}(X')$, $\mu \in \mu(E_1)$
\begin{equation} \label{ut5}
\beta_* ((\beta_{\gamma})_* (f) \otimes \mu) = \gamma_*
(\gamma_{\beta})_* (f \otimes \mu) \mbox{.}
\end{equation}
\item For any $G \in \D_{\beta(F_2(o))}'(E_3)$,
$\mu \in \mu(E_1)$
\begin{equation} \label{ut6}
\gamma_{\beta}^* (\gamma^* (G) \otimes \mu ) = \beta_{\gamma}^*
\beta^* (G \otimes \mu) \mbox{.}
\end{equation}
\end{enumerate}
\end{prop}
\proof  is analogous to the proof of proposition~\ref{przam1}
and~\ref{przam2}, and it is reduced to corresponding
formula~(\ref{pr18f4c})
 for $C_1$-spaces (to prove formula~(\ref{ut5}))
and to corresponding formula~(\ref{pr18f4d}) for  $C_1$-spaces (to
prove formula~(\ref{ut6})).
\vspace{0.3cm}

\begin{prop}.
Using notations of diagram~(\ref{zvezda}), let   $E_1$ be a
$cfC_2$-space, $B$ be a $dC_2$-space. Let $o \in I_2$, $ \nu \in
\mu(\delta \beta (F_2(o)) \mid \{ 0\})$. Then the following formulas
are satisfied.
\begin{enumerate}
\item
For any $f \in \D_{\beta(F_2(o))}(E_3)$
\begin{equation} \label{ut7}
\gamma_{\beta}^* \gamma^* (f \otimes \nu) = \beta_{\gamma}^*
(\beta^* (f) \otimes \nu) \mbox{.}
\end{equation}
\item For any $G \in \D_{F_2(o) \cap X}'(X')$
\begin{equation} \label{ut8}
\beta_* (\beta_{\gamma})_* (G \otimes \nu) = \gamma_*
((\gamma_{\beta})_* (G) \otimes \nu ) \mbox{.}
\end{equation}
\end{enumerate}
\end{prop}
\proof  is analogous to the proof of proposition~\ref{przam1}
and~\ref{przam2}, and it is reduced to corresponding
formula~(\ref{pr18f4a})
 for $C_1$-spaces (to prove formula~(\ref{ut7}))
and to corresponding formula~(\ref{pr18f4b}) for  $C_1$-spaces (to
prove formula~(\ref{ut8})).
\vspace{0.3cm}

\subsubsection{} \label{ss5.7.2}
We consider again diagram~(\ref{zvezda}). If $E_1$ and $D$ are
$cC_2$-spaces, then $E_2 \mathop{\times}\limits_{E_3} D$ is a
$cC_2$-space. Indeed, it follows from the admissible triple of
$C_2$-spaces:
\begin{equation} \label{diez}
0 \arrow{e} E_1  \arrow{e,t}{\gamma_{\alpha}} E_2
\mathop{\times}\limits_{E_3} D \arrow{e,t}{\gamma_{\beta}} D
\arrow{e} 0 \mbox{.}
\end{equation}

Let $o \in I_2$. Then from~(\ref{diez}) we obtain the canonical
isomorphism
\begin{equation} \label{2diez}
\mu(F_2(o) \cap X \mid X) = \mu (F_2(o) \cap V_1 \mid V_1)
\otimes_{\cc} \mu (\beta (F_2(o)) \cap Y \mid Y) \mbox{.}
\end{equation}
(The subspaces, which  appear in formula~(\ref{2diez}), are elements
of the filtration of corresponding $C_2$-spaces. Therefore the
spaces of virtual measures are well-defined, see section~\ref{vm}.)

We have from diagram~(\ref{zvezda}) the following admissible triple
of  $C_2$-spaces:
$$
0 \arrow{e} E_2 \mathop{\times}\limits_{E_3} D
\arrow{e,t}{\beta_{\gamma}} E_2 \arrow{e,t}{\delta  \beta} B
\arrow{e} 0 \mbox{.}
$$

\begin{prop}. \label{prop34}
Using notations of diagram~(\ref{zvezda}), let   $E_1$ and $D$ be
$cC_2$-spaces. Let $o \in I_2$, $\mu \in \mu(F_2(o) \cap V_1 \mid
V_1)$, $\nu \in \mu (\beta (F_2(o)) \cap Y \mid Y)$. Then the
following formulas are satisfied.
\begin{enumerate}
\item
For any $f \in \D_{F_2(o)}(E_2)$
\begin{equation} \label{ut9}
(\delta \beta)_* (f \otimes (\mu \otimes \nu)) = \delta_* ( \beta_*
(f \otimes \mu) \otimes \nu) \mbox{.}
\end{equation}
\item
For any $G \in \D'_{\delta \beta (F_2(o))}(B)$
\begin{equation} \label{ut10}
(\delta \beta)^* (G \otimes (\mu \otimes \nu))= \beta^* (\delta^* (G
\otimes \nu) \otimes \mu) \mbox{.}
\end{equation}
\end{enumerate}
\end{prop}
\begin{nt}. {\em
Due to~(\ref{2diez}) we have that $\mu \otimes \nu \in \mu(F_2(o)
\cap X \mid X)$. }
\end{nt}
\proof of  proposition~\ref{prop34}. Let $j_0 \in I_2$ be such that
$F_2(j_0) \cap V_1 = V_1$ and $\quad \beta (F_2(j_0)) \cap Y =Y$. We
consider arbitrary  $j \ge i \in I_2$ such that $j \ge j_0$. Then
formula~(\ref{e1}) from proposition~\ref{pr16} is satisfied for the
maps  $\beta_{ji}$, $\delta_{ji}$ and $(\beta_{\delta})_{ji} =
\delta_{ji} \beta_{ji}$ (see diagram~(\ref{2zvezda})). We multiply
this formula by corresponding measure spaces, and then we take the
limit with respect to  $j \ge i \in I_2$ ($j \ge j_0$). Using
explicit construction of direct image from proposition~\ref{pr26},
we obtain formula~(\ref{ut9}).

Formula~(\ref{ut10}) can be obtained by analogous reasonings, but
one has to use formula~(\ref{e2}) from proposition~\ref{pr16}.
\vspace{0.3cm}

Using notations of diagram~(\ref{zvezda}), we suppose that $E_1$ and
$D$ are $cfC_2$-spaces. Then from admissible triple of
$C_2$-spaces~(\ref{diez}) it follows that $E_2
\mathop{\times}\limits_{E_3} D$ is also a $cfC_2$-space.

\begin{prop}. \label{predut}
Using notations of diagram~(\ref{zvezda}), let   $E_1$ and $D$ be
$cfC_2$-spaces. Let $o \in I_2$. Then the following formulas are
satisfied.
\begin{enumerate}
\item
For any $f \in \D_{\delta \beta (F_2(o))}(B)$
\begin{equation} \label{ut11}
(\delta \beta)^* (f )= \beta^* \delta^* (f) \mbox{.}
\end{equation}
\item
For any $G \in \D'_{F_2(o)}(E_2)$
\begin{equation} \label{ut12}
(\delta \beta)_* (G ) = \delta_*  \beta_* (G) \mbox{.}
\end{equation}
\end{enumerate}
\end{prop}
\proof. We consider arbitrary  $j \ge i \in I_2$. Then
formula~(\ref{pr16f5}) from proposition~\ref{pr16} is satisfied for
the maps $\beta_{ji}$, $\delta_{ji}$ and $(\beta_{\delta})_{ji} =
\delta_{ji} \beta_{ji}$ (see diagram~(\ref{2zvezda})). We multiply
this formula by corresponding measure spaces, and then we take the
limit with respect to  $j \ge i \in I_2$. Using explicit
construction of inverse image from proposition~\ref{pr28}, we obtain
formula~(\ref{ut11}).

Formula~(\ref{ut12}) can be obtained by analogous reasonings, but
one has to use formula~(\ref{pr16f6}) from proposition~\ref{pr16}.
\vspace{0.3cm}

A little bit rewriting  (and redenoting) diagram~(\ref{zvezda}), we
obtain the following commutative diagram  of $C_2$-spaces over the
field $\df_q$:
\begin{equation} \label{3zvezda}
\begin{diagram}
\node[3]{0} \arrow{s} \node{0} \arrow{s} \\
 \node{0} \arrow{e} \node{E_1}
\arrow{e,t}{\alpha} \arrow{s,=} \node{E_2} \arrow{e,t}{\beta}
\arrow{s,l}{\alpha'} \node{E_3} \arrow{e} \arrow{s,r}{\beta_{\alpha'}} \node{0} \\
\node{0} \arrow{e} \node{E_1}  \arrow{e,t}{\alpha' \alpha} \node{H'}
\arrow{e,t}{\alpha'_{\beta}}  \arrow{s,l}{\theta} \node{E_3
\mathop{\amalg}\limits_{E_2} H'} \arrow{e}
\arrow{s,r}{\theta'} \node{0} \\
\node[3]{L'} \arrow{r,=} \arrow{s} \node{L'} \arrow{s} \\
\node[3]{0} \node{0}
\end{diagram}
\end{equation}
In this diagram two horizontal and two vertical triples are
admissible triples of  $C_2$-spaces.

Let $E_i = (I_i, F_i, V_i)$ ($1 \le i \le 3$), $H' = (J', T', W')$
as $C_2$-spaces. From admissible triple of  $C_2$-spaces
\begin{equation} \label{zv}
0 \arrow{e} E_3  \arrow{e,t}{\beta_{\alpha'}} E_3
\mathop{\amalg}\limits_{E_2} H' \arrow{e,t}{\theta'} L' \arrow{e} 0
\end{equation}
we obtain that if  $E_3$ and $L'$ are $dC_2$-spaces, then $ E_3
\mathop{\amalg}\limits_{E_2} H'$ is also a $dC_2$-space.

Let $ o \in J'$ be any, then from~(\ref{zv}) we canonically have
\begin{equation} \label{2zv}
\mu (\alpha_{\beta}' (T'(o)) \mid \{ 0 \} ) = \mu(\beta(T'(o) \cap
V_2) \mid \{ 0 \}) \otimes_{\cc} \mu(\theta(T'(o)) \mid \{ 0\})
\mbox{.}
\end{equation}
(The subspaces, which  appear in formula~(\ref{2zv}), are elements
of the filtration of corresponding $C_2$-spaces. Therefore the
spaces of virtual measures are well-defined, see section~\ref{vm}.)

From diagram~(\ref{3zvezda}) we have an admissible triple
$$
0 \arrow{e} E_1  \arrow{e,t}{\alpha' \alpha} H'
\arrow{e,t}{\alpha_{\beta}'}  E_3 \mathop{\amalg}\limits_{E_2} H'
\arrow{e} 0 \mbox{.}
$$

\begin{prop}.
Using notations of diagram~(\ref{3zvezda}), let   $E_3$ and $L'$ be
$dC_2$-spaces. Let $o \in J'$, $\mu \in \mu(\beta(T'(o) \cap V_2)
\mid \{ 0 \}) $, $\nu \in \mu(\theta(T'(o)) \mid \{ 0\})$. Then the
following formulas are satisfied.
\begin{enumerate}
\item
For any $f \in \D_{T'(o)} (H')$
\begin{equation} \label{ut13}
(\alpha' \alpha)^* (f \otimes (\mu \otimes \nu))= \alpha^*
((\alpha')^*(f \otimes \nu) \otimes \mu ) \mbox{.}
\end{equation}
\item
For any $G \in \D'_{T'(o) \cap V_1} (E_1)$
\begin{equation} \label{ut14}
(\alpha' \alpha)_* (G \otimes (\mu \otimes \nu)) = (\alpha')_*
(\alpha_*(G \otimes \mu) \otimes \nu) \mbox{.}
\end{equation}
\end{enumerate}
\end{prop}
\begin{nt}. {\em
Due to~(\ref{2zv}) we have that $\mu \otimes \nu \in
\mu(\alpha_{\beta}' (T'(o)) \mid \{ 0 \})$. }
\end{nt}
\proof of proposition. Formula~(\ref{ut13}) can be reduced to
formula~(\ref{e'1}) from proposition~\ref{pr17} after the taking the
projective limit with respect to  $j \ge i \in J'$, where $i \le
i_0$, and $i_0$ is chosen such that  $\theta(T'(i_0)) = \{ 0\}$ and
$\beta (T'(i_0) \cap V_2) = \{ 0\}$. At that we use the explicit
construction of inverse image from proposition~\ref{pr27}. (Compare
 with the proof of proposition~\ref{predut}.)

Formula~(\ref{ut14}) can be proved analogously, but one has to use
formula~(\ref{e'2}) from proposition~\ref{pr17} and take the
inductive limit.
\vspace{0.3cm}

Using notations of diagram~(\ref{3zvezda}) we suppose that  $E_3$
and $L'$ are $dfC_2$-spaces. Then from admissible triple of
$C_2$-spaces~(\ref{zv}) it follows that $E_3
\mathop{\amalg}\limits_{E_2} H'$ is also a  $dfC_2$-space.

\begin{prop}. Using notations of diagram~(\ref{3zvezda}), let  $E_3$ and $L'$
are $dfC_2$-spaces. Let $o \in J'$. Then the following formulas are
satisfied.
\begin{enumerate}
\item
For any $f \in \D_{T'(o) \cap V_1} (E_1)$
\begin{equation} \label{ut15}
(\alpha' \alpha)_* (f ) = (\alpha')_* \alpha_*(f ) \mbox{.}
\end{equation}
\item
For any $G \in \D'_{T'(o)} (H')$
\begin{equation} \label{ut16}
(\alpha' \alpha)^* (G )= \alpha^* ((\alpha')^*(G ) \mbox{.}
\end{equation}
\end{enumerate}
\end{prop}
\proof of formula~(\ref{ut15}) can be reduced to
formula~(\ref{pr17f5}) of proposition~\ref{pr17} after the choice of
 $j \ge i \in J'$ and taking the projective limit with respect to this  $j \ge
i \in J'$.  At that we use the explicit construction of direct image
from proposition~\ref{pr29}.

The proof of formula~(\ref{ut16}) is analogous, but one has to take
the inductice limit with respect to   $j \ge i \in J'$ and use
formula~(\ref{pr17f6}) instead of formula~(\ref{pr17f5}).
\vspace{0.3cm}

\subsection{Fourier transform and direct and inverse images}
\label{ftdi}
Let
$$
0\arrow{e} E_1\arrow{e,t}{\alpha} E_2\arrow{e,t}{\beta}E_3\arrow{e}0
$$
be an admissible triple of $C_2$-spaces over a finite field $\df_q$,
where $E_i=(I_i,F_i,V_i)$, \linebreak $1 \le i \le 3$. We suppose
that $E_i$ ($1 \le i\le 3$) is a {\it complete} $C_2$-space. Then
$\check{E}_i=E_i$, $1 \le i \le 3$.

By definition, there are order-preserving  functions:
\begin{itemize}
\item[(i)]
$ \gamma \; \colon \; I_2\arrow{e} I_3\quad\mbox{such that}\quad
\beta(F_2(i))=F_3(\gamma(i))\quad\mbox{for any}\quad i\in I_2
\mbox{;} $
\item[(ii)]
$ \veps \; \colon \; I_2\arrow{e} I_1\quad\mbox{such that}\quad
F_2(i)\cap V_1=F_1(\veps(i))\quad\mbox{for any}\quad i\in
I_2\mbox{.} $
\end{itemize}
We consider the following admissible triple:
$$
0\arrow{e} \check{E_3} \arrow{e,t}{\check\beta} \check{ E_2}
\arrow{e,t}{\check\alpha} \check{E_1} \arrow{e}0 \mbox{.}
$$
We have $\check{E_i}=(I^0_i,F^0_i,\check{V_i})$, $1 \le i \le 3$. We
recall that the partially ordered set $I^{0}_i$ ($1 \le i \le 3$) is
equal to $I_i$ as a set, but has inverse partial order comparing
with $I_i$ ($1 \le i \le 3$). We have that for any $i,j \in I_2$:
\begin{equation} \label{bigstar}
\begin{array}{c}
\check\alpha (F^{0}_2(i))=F^{0}_1(\veps(i)) \\[0.2cm]
F_2^0(i)\cap \check{V_3}=F^{0}_3(\gamma(i)) \mbox{.}
\end{array}
\end{equation}

For any $i,j\in I_2$ we have from formula~(\ref{fff}) that
\begin{equation} \label{bigstarbigstar}
\begin{array}{c}
\mu(F_1(\veps(i))\mid F_1(\veps(j)))=\mu(F^{0}_1(\veps(i))\mid
F^{0}_1(\veps(j))) \\[0.2cm]
\mu(F_3(\gamma(i))\mid F_3(\gamma(j)))=\mu(F^{0}_3(\gamma(i))\mid
F^{0}_3(\gamma(j))) \mbox{.}
\end{array}
\end{equation}

We have the following statements.
\begin{prop} \label{pppp}
Let $o \in I_2$. The following diagrams are
commutative.
\begin{enumerate}
\item If $E_1$ is a $cC_2$-space. Then
$$
\begin{diagram}
\node{\D_{F_2(o)}(E_2)\otimes_{\cc}\mu(F_1(\veps(o))\mid V_1)}
\arrow{e,t}{\beta_*} \arrow{s,l}{\F\otimes\mu(F_1(\veps(o))\mid
V_1)}
\node{\D_{F_3(\gamma(o))}(E_3)}\arrow{s,r}{\F}\\
\node{\D_{F^{0}_2(o)}(\check{
E_2})\otimes_{\cc}\mu(F^{0}_1(\veps(o))\mid \{0\})}
\arrow{e,t}{(\check\beta)^*} \node{\D_{F^{0}_3(\gamma(o))}(\check
E_3) \mbox{.}}
\end{diagram}
$$
\item If $E_3$ is a $dC_2$-space. Then
$$
\begin{diagram}
\node{\D_{F_2(o)}(E_2)\otimes_{\cc}\mu(F_3(\gamma(o))\mid\{0\})}
\arrow{e,t}{\alpha^*}
\arrow{s,l}{\F\otimes\mu(F_3(\gamma(o))\mid\{0\})}
\node{\D_{F_1(\veps(o))}(E_1)}\arrow{s,r}{\F}\\
\node{\D_{F^{0}_2(o)}(\check{
E_2})\otimes_{\cc}\mu(F^{0}_3(\gamma(o))\mid \check{ V_3})}
\arrow{e,t}{(\check\alpha)_*}
\node{\D_{F^{0}_1(\veps(o))}(\check{E_1}) \mbox{.}}
\end{diagram}
$$
\item If $E_1$ is a $cC_2$-space. Then
$$
\begin{diagram}
\node{\D'_{F_3(\gamma(o))}(E_3)\otimes_{\cc}\mu(F_1(\veps(o))\mid
V_1)} \arrow{e,t}{\beta^*}\arrow{s,l}{\F\otimes\mu(F_1(\veps(o))\mid
V_1)}
\node{\D'_{F_2(o)}(E_2)}\arrow{s,r}{\F}\\
\node{\D'_{F^{0}_3(\gamma(o))}(\check{
E_3})\otimes_{\cc}\mu(F^{0}_1(\veps(o))\mid \{0\})}
\arrow{e,t}{(\check\beta)_*} \node{\D'_{F^{0}_2(o)}(\check{E_2})
\mbox{.}}
\end{diagram}
$$
\item If $E_3$ is a $dC_2$-space. Then
$$
\begin{diagram}
\node{\D'_{F_1(\veps(o))}(E_1)\otimes_{\cc}\mu(F_3(\gamma(o))\mid
\{0\})}
\arrow{e,t}{\alpha_*}\arrow{s,l}{\F\otimes\mu(F_3(\gamma(o))\mid
\{0\})}
\node{\D'_{F_2(o)}(E_2)}\arrow{s,r}{\F}\\
\node{\D'_{F^{0}_1(\veps(o))}(\check{
E_1})\otimes_{\cc}\mu(F^{0}_3(\gamma(o))\mid \check{V_3})}
\arrow{e,t}{(\check\alpha)^*} \node{\D'_{F^{0}_2(o)}(\check{ E_2})
\mbox{.}}
\end{diagram}
$$
\item If $E_1$ is a $cfC_2$-space. Then
$$
\begin{diagram}
\node{\D_{F_3(\gamma(o))}(E_3)} \arrow{e,t}{\beta^*}\arrow{s,l}{\F}
\node{\D_{F_2(o)}(E_2)}\arrow{s,r}{\F}\\
\node{\D_{F^{0}_3(\gamma(o))}(\check{E_3})}
\arrow{e,t}{(\check\beta)_*} \node{\D_{F^{0}_2(o)}(\check{E_2})
\mbox{.}}
\end{diagram}
$$
\item If $E_3$ is a $dfC_2$-space. Then
$$
\begin{diagram}
\node{\D_{F_1(\veps(o))}(E_1)}\arrow{e,t}{\alpha_*}\arrow{s,l}{\F}
\node{\D_{F_2(o)}(E_2)}\arrow{s,r}{\F}\\
\node{\D_{F^{0}_1(\veps(o))}(\check{E_1})}
\arrow{e,t}{(\check\alpha)^*} \node{\D_{F^{0}_2(o)}(\check{E_2})
\mbox{.} }
\end{diagram}
$$
\item If $E_1$ is a $cfC_2$-space. Then
$$
\begin{diagram}
\node{\D'_{F_2(o)}(E_2)} \arrow{e,t}{\beta_*}\arrow{s,l}{\F}
\node{\D'_{F_3(\gamma(o))}(E_3)}\arrow{s,r}{\F}\\
\node{\D'_{F^{0}_2(o)}(\check{E_2})} \arrow{e,t}{(\check\beta)^*}
\node{\D'_{F^{0}_3(\gamma(o))}(\check{E_3}) \mbox{.}}
\end{diagram}
$$
\item If $E_3$ is a $dfC_2$-space. Then
$$
\begin{diagram}
\node{\D'_{F_2(o)}(E_2)} \arrow{e,t}{\alpha^*}\arrow{s,l}{\F}
\node{\D'_{F_1(\veps(o))}(E_1)}\arrow{s,r}{\F}\\
\node{\D'_{F^{0}_2(o)}(\check{E_2})} \arrow{e,t}{(\check\alpha)_*}
\node{\D'_{F^{0}_1(\veps(o))}(\check{E_1}) \mbox{.}}
\end{diagram}
$$
\end{enumerate}
\end{prop}
\proof  follows from formulas (\ref{bigstar}) and
(\ref{bigstarbigstar}), constructions of two-dimensional direct and
inverse images from section~\ref{secdi}, two-dimensional Fourier
transform, and corresponding properties of direct and inverse images
and of Fourier transform on $C_1$-spaces from section~\ref{dif}.
\vspace{0.3cm}

\subsection{Two-dimensional Poisson formulas} \label{2PF}
In this section we obtain two-dimensional Poisson formulas: the
Poisson formula I (subsection~\ref{s5.9.1}) and the Poisson formula
II (subsection~\ref{s5.9.2}).
\subsubsection{}
\label{s5.9.1} Let
$$
0\arrow{e} E_1\arrow{e,t}{\alpha} E_2\arrow{e,t}{\beta}E_3\arrow{e}0
$$
be an admissible triple of complete $C_2$-spaces over a finite field
$\df_q$, where \linebreak $E_i=(I_i,F_i,V_i)$, $1 \le i \le 3$. By
definition, there are order-preserving  functions:
\begin{itemize}
\item[(i)]
$ \gamma \; \colon \; I_2\arrow{e} I_3\quad\mbox{such that}\quad
\beta(F_2(i))=F_3(\gamma(i))\quad\mbox{for any}\quad i\in I_2
\mbox{;} $
\item[(ii)]
$ \veps \; \colon \; I_2\arrow{e} I_1\quad\mbox{such that}\quad
F_2(i)\cap V_1=F_1(\veps(i))\quad\mbox{for any}\quad i\in I_2
\mbox{.} $
\end{itemize}

Let $o\in I_2$. If we suppose that $E_1$ is a $cC_2$-space, then
from any \linebreak $\mu\in \mu(F_1(\veps(o))\mid V_1)$ we
canonically construct $\znakneponyaten_{\mu} \in
\D'_{F_1(\veps(o))}(E_1)$  in the following way. By definition, we
have
\begin{equation} \label{nabla}
\D'_{F_1(\veps(o))}(E_1)=\mathop{\Lim_{\lto}}_{k\in I_1}
\D'(V_1/F_1(k))\otimes_{\cc}\mu(F_1(\veps(o))\mid k) \mbox{.}
\end{equation}
By definition, $\mu(F_1(\veps(o))\mid V_1)\subset
\D'(V_1/F_1(\veps(o)))$. Therefore $\mu\in \D'(V_1/F_1(\veps(o)))$.
In formula (\ref{nabla}) we take $k=\veps(o)$. Then
$\mu(F_1(\veps(o))\mid k)=\cc$ and we put \linebreak $1 \in
\mu(F_1(\veps(o))\mid \veps(o))$. Therefore $\mu$ gives
$\znakneponyaten_\mu\in \D'_{F_1(\veps(o))}(E_1)$ by formula
(\ref{nabla}).

Let $o \in I_2$. If we suppose that $E_3$ is a $dC_2$-space, then
from any \linebreak $\nu \in \mu(F_3(\gamma(o))\mid \{0\})$ we
canonically construct $\delta_\nu \in \D'_{F_3(\gamma(o))}(E_3)$ in
the following way. By definition, we have
\begin{equation} \label{nablanabla}
\D'_{F_3(\gamma(o))}(E_3)=\mathop{\Lim_{\lto}}_{k\in I_2}
\D'(F_3(k))\otimes_{\cc}\mu(F_3(\gamma(o))\mid \{0\}) \mbox{.}
\end{equation}
In formula (\ref{nablanabla}) we take  $k=\gamma(o)$. Then we put
$\delta_0 \in \D'(F_3(\gamma(o)))$, where \linebreak
$\delta_0(f)=f(0)$, $f\in \D(F_3(\gamma(o)))$. We put $\nu \in
\mu(F_3(\gamma(o))\mid \{0\})$ in formula (\ref{nablanabla}). Thus
we have constructed $\delta_\nu\in \D'_{F_3(\gamma(o))}(E_3)$.

We note that
\begin{equation} \label{triangle}
\F(\znakneponyaten_\mu)=\delta_\mu\quad\mbox{and}\quad
\F(\delta_\nu)=\znakneponyaten_\nu
\end{equation}
under the maps
$$\F \quad \colon \quad \D'_{F_1(\veps(o))}(E_1)\arrow{e}
\D'_{F^{0}_1(\veps(o))}(\check{E_1}) \quad  \mbox{and}
$$
$$
\F \quad \colon \quad \D'_{F_3(\gamma(o))}(E_3)\arrow{e}
\D'_{F^{0}_3(\gamma(o))}(\check{ E_3}) \mbox{.}
$$
(We used that
$$
\mu(F_1(\veps(o))\mid V_1)=\mu(F^{0}_1(\veps(o))\mid
\{0\})\quad\mbox{and}\quad \mu(F_3(\gamma(o))\mid
\{0\})=\mu(F^{0}_3(\gamma(o))\mid \check{V_3}) \mbox{.})
$$

Now we suppose that simultaneously $E_1$ is a $cC_2$-space and $E_3$
is a $dC_2$-space. Let $o \in I_2$. Then by any $\mu\in
\mu(F_1(\veps(o))\mid V_1)$, $\nu\in \mu(F_3(\gamma(o))\mid \{0\})$
it is well-defined {\em the characteristic function}
\begin{equation}
\delta_{E_1,\mu\otimes \nu}\eqdef\alpha_*(\znakneponyaten_\mu\otimes
\nu)\in \D'_{F_2(o)}(E_2).
\end{equation}

We have the following statement.
\begin{lemma} \label{lem5}
$$
\delta_{E_1,\mu\otimes \nu}=\beta^*(\delta_\nu\otimes \mu).
$$
\end{lemma}
\proof follows from constructions and corresponding statements for
the following admissible triple of $C_1$-spaces:
$$
0\arrow{e} V_1/F_1(\veps(k_2))\arrow{e}\frac{F_2(k_1)}{F_2(k_2)}\arrow{e}
F_3(\gamma(k_1))\arrow{e}0,
$$
where $k_1\ge o \ge k_2\in I_2$ and $F_1(\veps(k_1))=V_1$,
$F_3(\gamma(k_2))=\{0\}$.
\vspace{0.3cm}

Now we have the following first two-dimensional Poisson formula.
\begin{Th}[\it Poisson formula I]. \label{tpf1}
Let $E_1$ be a $cC_2$-space, Let $E_3$ be a $dC_2$-space. Let $o \in
I_2$. Then for any $\mu\in \mu(F_1(\veps(o))\mid V_1)$, $\nu\in
\mu(F_3(\gamma(o))\mid \{0\})$ we have
$$
\F(\delta_{E_1,\mu\otimes \nu})=\delta_{\check{E_3},\nu\otimes \mu}
\mbox{.}
$$
\end{Th}
\proof follows from lemma~\ref{lem5}, from
formulas~(\ref{triangle}), and from proposition~\ref{pppp} about
connection of direct and inverse images with two-dimensional Fourier
transform. The theorem is proved.
 \vspace{0.3cm}

\begin{nt} \label{rem23} {\em
Let $E=(I, F, V)$ be a $C_2$-space over a field $k$. We construct a
$C_2$-space $\overline{E}$, which dominates a $C_2$-space $E$. We
define $\overline{E}= (I, F, V)$, but for every $i \le j \in I$ the
$C_1$-space $\overline{E}_{ij} $ is the maximal $C_1$-space which
dominates the $C_1$-space $E_{ij}$, where $\overline{E}_{ij}$ and
$E_{ij}$ are $C_1$-spaces with underlying $k$-vector space
$F(j)/F(i)$. (To construct for a $C_1$-space $(K, G, W)$ its maximal
dominating $C_1$-space $(K', G', W)$ we denote by $Gr(W)$ the set of
$k$-vector subspaces $U \subset W$ such that there are $k \le l \in
K $ with the property $G(k) \subset U \subset G(l)$. It is clear
that $Gr(W)$ is a partially ordered set which is ordered by
embeddings of $k$-subspaces. Then we put $K'= Gr(W)$, and the
function $G'$ is an embedding of the corresponding $k$-subspace to
$W$.)

If
$$ 0\arrow{e} E_1\arrow{e}
E_2\arrow{e} E_3 \arrow{e} 0
$$
is an admissible triple of $C_2$-spaces, then
$$ 0 \arrow{e} \overline{E_1} \arrow{e}
\overline{E_2} \arrow{e} \overline{E_3} \arrow{e}0
$$
is also an admissible triple of $C_2$-spaces. }
\end{nt}

We consider any admissible triple of $C_2$-spaces over a finite
field $\df_q$
$$
0\arrow{e} E_1\arrow{e,t}{\alpha}
E_2\arrow{e,t}{\beta}E_3\arrow{e}0,
$$
where $E_i=(I_i,F_i,V_i)$, $1 \le i \le 3$, and $E_1$ is a
$cC_2$-space,  $E_3$ is a $dC_2$-space. Let
  $g\in \Aut_{C_2}(E_2)$ such that for any $i \in I_2$ there is $j \in I_2$ with the property
$g F_2(i)= F_2(j)$. Then
   we have the following
admissible triple of $C_2$-spaces
\begin{equation} \label{kvadrat}
0\arrow{e} g\overline{E_1} \arrow{e,t}{g(\alpha)} g \overline{ E_2}
\arrow{e,t}{g(\beta)} g \overline{E_3} \arrow{e} 0 \mbox{,}
\end{equation}
where the $C_2$-space $\overline{E_i}$  dominates the $C_2$-space
$E_i$ ($1 \le i \le 3$) (see remark~\ref{rem23}), $g \overline{
E_2}=\overline{E_2}$, $g \overline{E_3}= \overline{ E_2}/g
\overline{E_1}$.

From~(\ref{kvadrat}) we have canonically
$$
\mu(F_2(o)\mid gF_2(o))=\mu(F_2(o) \cap gV_1 \mid
gF_1(\veps(o)))\otimes_{\cc}\mu(g(\beta)(F_2(o)) \mid
gF_3(\gamma(o))) \mbox{.}
$$
Then for any $a\in\mu(F_2(0)\mid gF_2(0))$ we consider $a=b\otimes
c$, where
$$
b\in\mu(F_2(o) \cap gV_1 \mid gF_1(\veps(o))) \mbox{,} \quad
c\in\mu(g(\beta)(F_2(o)) \mid gF_3(\gamma(o))).
$$
Under conditions and notations of theorem~\ref{tpf1}  we have
$$
g\mu\in\mu(gF_1(\veps(o))\mid gV_1) \mbox{,} \qquad
g\nu\in\mu(gF_3(\gamma(o))\mid \{0\}) \mbox{.}
$$
We have $\widetilde g=(g,a)\in\widehat{\Aut_{C_2}(E_2)}_{F_2(o)}$.
Then
$$
R'_{\widetilde g}(\delta_{E_1,\mu\otimes\nu})=\delta_{gE_1,
(g\mu\otimes b)\otimes(g\nu\otimes c)}=\delta_{gE_1, g\mu\otimes
g\nu\otimes a}.
$$

\begin{cons}
Let $(g,a)\in\widehat{\Aut_{C_2}(E_2)}_{F_2(o)}$ such that for any
$i \in I_2$ there is $j \in I_2$ with the property $g F_2(i)=
F_2(j)$. Then we have
$$
\F(\delta_{gE_1, g\mu\otimes g\nu\otimes a})=\delta_{\check
g^{-1}(\check{E_3}),\check g^{-1}(\nu)\otimes\check
g^{-1}(\mu)\otimes a}.
$$
\end{cons}
\proof follows from theorem~\ref{tpf1} (Poisson formula I)  and
section~\ref{sec5.5.5}.
\vspace{0.3cm}

\subsubsection{}
\label{s5.9.2}
 Let
$$
0\arrow{e} E_1\arrow{e,t}{\alpha} E_2\arrow{e,t}{\beta}E_3\arrow{e}0
$$
be an admissible triple of complete $C_2$-spaces over a finite field
$\df_q$, where \linebreak $E_i=(I_i,F_i,V_i)$, $1 \le i \le 3$. By
definition, there are order-preserving  functions:
\begin{itemize}
\item[(i)]
$ \gamma \; \colon \; I_2\arrow{e} I_3\quad\mbox{such that}\quad
\beta(F_2(i))=F_3(\gamma(i))\quad\mbox{for any}\quad i\in I_2
\mbox{;} $
\item[(ii)]
$
\veps \; \colon \; I_2\arrow{e} I_1\quad\mbox{such that}\quad
F_2(i)\cap V_1=F_1(\veps(i))\quad\mbox{for any}\quad i\in I_2
\mbox{.}
$
\end{itemize}

If $E_1$ is a $cfC_2$-space, then using remark~\ref{remark9} we have
$$
\D(E_1)=\mathop{\Lim_{\longleftarrow}}_{i\ge j\in I_1}
\D(F_1(i)/F_1(j)) \mbox{.}
$$
$F_1(i)/F_1(j)$ is a compact $C_1$-space (for any $i\ge j\in I_1$).
Therefore \linebreak $1\in \D(F_1(i)/F_1(j))$. Taking the projective
limits, we obtain that in this case the following element is
well-defined
\begin{equation}
\znakneponyaten \in \D(E_1) \mbox{.}
\end{equation}

If $E_3$ is a $dfC_2$-space, then using remark~\ref{remark9} we have
$$
\D(E_3)=\mathop{\Lim_{\longleftarrow}}_{i\ge j\in I_3}
\D(F_3(i)/F_3(j)) \mbox{.}
$$
$F_3(i)/F_3(j)$ is a discrete $C_1$-space (for any $i\ge j\in I_3$).
Therefore \linebreak $\delta_0 \in \D(F_3(i)/F_3(j))$, where
$\delta_0(0)=1$. Taking the projective limits, we obtain in this
case a well-defined element
\begin{equation}
\delta_0 \in \D(E_3) \mbox{.}
\end{equation}

We note that
\begin{equation} \label{triangletriangle}
\F(\znakneponyaten)=\delta_0\quad\mbox{and}\quad
\F(\delta_0)=\znakneponyaten
\end{equation}
under the maps
$$\F \quad \colon \quad \D(E_1)\arrow{e} \D(\check E_1) \quad \mbox{and} $$
$$\F \quad \colon \quad \D(E_3)\arrow{e}
\D(\check E_3) \mbox{. }$$

Now we suppose that simultaneously $E_1$ is a $cfC_2$-space and
$E_3$ is a $dfC_2$-space. Let $o \in I_2$. Then an element
\begin{equation}
\delta_{E_1} \eqdef \alpha_*(\znakneponyaten) \in \D_{F_2(o)}(E_2)
\end{equation}
is well-defined.

We have the following statement.
\begin{lemma} \label{lem6}
$$
\delta_{E_1}=\beta^*(\delta_0) \mbox{.}
$$
\end{lemma}
\proof follows from constructions and corresponding statements for
the following admissible triples of $C_1$-spaces  (for any $i\ge
j\in I_2$)
$$
0\arrow{e}\frac{F_1(\veps(i))}{F_1(\veps(j))}\arrow{e} \frac{F_2(i)}{F_2(j)}\arrow{e}
\frac{F_3(\gamma(i))}{F_3(\gamma(j))}\arrow{e}0,
$$
where $\frac{F_1(\veps(i))}{F_1(\veps(j))}$ is a compact $C_1$-space
and $\frac{F_3(\gamma(i))}{F_3(\gamma(j))}$ is a discrete
$C_1$-space.
\vspace{0.3cm}

Now we have the following second two-dimensional Poisson formula.
\begin{Th}[\it Poisson formula II]. \label{th3}
Let $E_1$ be a $cfC_2$-space, Let $E_3$ be a $dfC_2$-space. Let $o
\in I_2$. Then
$$
\F(\delta_{E_1})=\delta_{\check{E_3}} \mbox{.}
$$
\end{Th}
\proof follows from lemma~\ref{lem6}, from
formulas~(\ref{triangletriangle}), and from proposition~\ref{pppp}
about direct and inverse images and two-dimensional Fourier
transform. The theorem is proved.
\vspace{0.3cm}

We consider any admissible triple of $C_2$-spaces over a finite
field $\df_q$
$$
0\arrow{e} E_1\arrow{e,t}{\alpha}
E_2\arrow{e,t}{\beta}E_3\arrow{e}0,
$$
where $E_i=(I_i,F_i,V_i)$, $1 \le i \le 3$, and $E_1$ is a
$cfC_2$-space,  $E_3$ is a $dfC_2$-space. Let
  $g\in \Aut_{C_2}(E_2)$ such that for any $i \in I_2$ there is $j \in I_2$ with the property
$g F_2(i)= F_2(j)$. Then
   we have the following
admissible triple of $C_2$-spaces
\begin{equation} \label{kvadrat}
0\arrow{e} g\overline{E_1} \arrow{e,t}{g(\alpha)} g \overline{ E_2}
\arrow{e,t}{g(\beta)} g \overline{E_3} \arrow{e} 0 \mbox{,}
\end{equation}
where the $C_2$-space $\overline{E_i}$  dominates the $C_2$-space
$E_i$ ($1 \le i \le 3$) (see remark~\ref{rem23}), $g \overline{
E_2}=\overline{E_2}$, $g \overline{E_3}= \overline{ E_2}/g
\overline{E_1}$.

\begin{cons}
Let $g \in \Aut_{C_2}(E_2)$ such that for any $i \in I_2$ there is
$j \in I_2$ with the property $g F_2(i)= F_2(j)$. Then we have
$$
\F(\delta_{gE_1})=\delta_{\check g^{-1}(\check{E_3})}.
$$
\end{cons}
\proof follows from theorem~\ref{th3} (Poisson formula II)  and
section~\ref{sec5.5.5}.

\vspace{0.3cm}

\noindent D.V. Osipov\\
Steklov Mathematical Institute RAS \\
{\it E-mail:}  ${d}_{-} osipov@mi.ras.ru$ \\

\vspace{0.3cm}

\noindent A.N. Parshin\\
Steklov Mathematical Institute RAS \\
{\it E-mail:}  $parshin@mi.ras.ru$
\end{document}